\RequirePackage{fix-cm}
\documentclass[a4paper,12pt]{amsart}
\usepackage{amssymb}
\usepackage{amsfonts}
\usepackage{amsmath}
\usepackage{mathrsfs}
\usepackage{enumerate}

\numberwithin{equation}{section}

\def\qed{\hfill $\Box$} 
\makeatletter
    
    \@addtoreset{equation}{section}
  \makeatother

\newcommand{\K}{{\mathbb{K}}}
\newcommand{\R}{{\mathbb{R}}}
\newcommand{\C}{{\mathbb{C}}}
\newcommand{\Z}{{\mathbb{Z}}}

\newcommand{\Si}{\Sigma}
\newcommand{\SiAp}{\Sigma_{A}^{+}}
\newcommand{\SiBp}{\Sigma_{B}^{+}}
\newcommand{\SiM}{\Sigma_{M}}
\newcommand{\SiMp}{\Sigma_{M^{\prime}}}
\newcommand{\SiMpp}{\Sigma_{M^{\prime\prime}}}
\newcommand{\si}{\sigma}
\newcommand{\Ph}{\varPhi}
\newcommand{\Phe}{\varPhi(\epsilon,\cdot)}
\newcommand{\tPhe}{\tilde{\varPhi}(\epsilon,\cdot)}
\newcommand{\ph}{\varphi}
\newcommand{\phe}{\varphi(\epsilon,\cdot)}
\newcommand{\lam}{\lambda}
\newcommand{\e}{\epsilon}
\newcommand{\ep}{\epsilon^{\prime}}

\newcommand{\p}{{\prime}}
\newcommand{\pp}{{\prime\prime}}
\newcommand{\ppp}{{\prime\prime\prime}}
\newcommand{\LR}{\mathcal{L}}
\newcommand{\tLR}{\tilde{\mathcal{L}}}
\newcommand{\bT}{\mathbb{T}}
\newcommand{\ali}[1]{\begin{align*}#1\end{align*}}
\newcommand{\alil}[1]{\begin{align}#1\end{align}}
\newcommand{\ite}[1]{\begin{enumerate}[(1)]#1\end{enumerate}}
\newcommand{\prop}[1]{\begin{proposition2}#1\end{proposition2}}
\newcommand{\thm}[1]{\begin{theorem2}#1\end{theorem2}}
\newcommand{\cor}[1]{\begin{corollary2}#1\end{corollary2}}
\newcommand{\lem}[1]{\begin{lemma2}#1\end{lemma2}}

\newcommand{\rems}{\begin{remark2}\normalfont}
\newcommand{\reme}{\end{remark2}}
\newcommand{\pros}{\begin{proof}}
\newcommand{\proe}{\end{proof}}
\newcommand{\case}[1]{\begin{cases}#1\end{cases}}
\newcommand{\cd}{\cdot}

\newcommand{\om}{\omega}
\newcommand{\omp}{\omega^{\prime}}

\newcommand{\bi}{\mathbf{i}}

\newcommand{\bbJ}{\mathfrak{J}}

\newcommand{\mV}{\mathscr{V}}
\newcommand{\bV}{\mathbf{V}}
\newcommand{\bU}{\mathbf{U}}
\newcommand{\mG}{\mathscr{G}}
\newcommand{\bG}{\mathbf{G}}

\newcommand{\tmV}{\tilde{\mathscr{V}}}
\newcommand{\tbV}{\tilde{\mathbf{V}}}
\newcommand{\tbU}{\tilde{\mathbf{U}}}
\newcommand{\tmG}{\tilde{\mathscr{G}}}
\newcommand{\tbG}{\tilde{\mathbf{G}}}

\newcommand{\SP}{\mathrm{SP}}
\newcommand{\bM}{\mathbb{M}}
\newcommand{\bbM}{\mathfrak{M}}
\newcommand{\IR}{\mathcal{I}}
\newcommand{\te}{t_{\epsilon}}
\newcommand{\tte}{\tilde{t}_{\epsilon}}
\newcommand{\ttte}{\hat{t}_{\epsilon}}
\newcommand{\GTe}{\mathrm{GSP}_{\epsilon}}
\newcommand{\VTe}{\mathrm{VSP}_{\epsilon}}
\newcommand{\LT}{\mathcal{L}}
\newcommand{\RLR}{\mathrm{RL}}
\newcommand{\LRR}{\mathrm{LR}}
\newcommand{\LRLR}{\mathrm{LRL}}
\newcommand{\RR}{\mathcal{R}}
\newcommand{\XR}{\mathcal{X}}
\newcommand{\mM}{L}
\newcommand{\ti}[1]{\tilde{#1}}
\renewcommand{\MR}{$\mathrm{SMR}$ }

\newcommand{\matII}[1]{\left(\begin{array}{cc}#1\end{array}\right)}

\newcommand{\matIV}[1]{\left(\begin{array}{cccc}#1\end{array}\right)}

\newcommand{\matVI}[1]{\left(\begin{array}{cccccc}#1\end{array}\right)}

\newcommand{\matVIII}[1]{\left(\begin{array}{cccccccc}#1\end{array}\right)}
\newcommand{\dia}{\mathrm{diag}}
\newcommand{\supp}{\mathrm{supp}\,}

\makeatletter
    
    \@addtoreset{equation}{section}
  \makeatother

\begin{document}

\newtheorem{theorem2}{\sc \bf Theorem}[section]
\newtheorem{lemma2}[theorem2]{\sc \bf Lemma}
\newtheorem{corollary2}[theorem2]{\sc \bf Corollary}
\newtheorem{proposition2}[theorem2]{\sc \bf Proposition}
\newtheorem{claim2}[theorem2]{\sc \bf Claim}
\newtheorem{definition2}[theorem2]{\sc \bf Definition}
\newtheorem{remark2}[theorem2]{\sc \bf Remark}
\newtheorem{example2}[theorem2]{\sc \bf Example}
\newtheorem{notation2}[theorem2]{\sc \bf Notation}
\newtheorem{assertion2}[theorem2]{\sc \bf Assertion}

\author{Haruyoshi Tanaka}
\date{June 4, 2016}
\title[Perturbation analysis in thermodynamic]{Perturbation analysis in thermodynamics using matrix representations of Ruelle transfer operators}

\address{
{\rm Haruyoshi Tanaka}\\
Department of Mathematics and Statistics,\\
Wakayama Medical University, \\
580, Mikazura, Wakayama-city, Wakayama, 641-0011, Japan\\
}
\email{htanaka@wakayama-med.ac.jp}

\setcounter{tocdepth}{1}

\maketitle

\begin{abstract}
We study perturbations of topological pressures, Gibbs measures and measure-theoretic entropies of these measures concerning perturbed potentials defined on topologically transitive subshift of finite type. The subshift with respect to non-perturbed system is assumed to be no topologically transitive. Therefore, the subshift of the perturbed systems and the subshift of the unperturbed system are different. We reduce this situation to a perturbation problem of certain irreducible nonnegative matrices generated by Ruelle transfer operators. Consequently, under suitable conditions of potentials, we characterize the limit points of those thermodynamics and give a necessary and sufficient condition for convergence of Gibbs measures and the measure-theoretic entropy of this measure when the subshift of the non-perturbed system has $2$ or $3$ transitive components with the maximal pressure. Finally, we illustrate the relation between potentials and convergence of Gibbs measures by using asymptotic expansion techniques for eigenvalues of Ruelle transfer operators.
\end{abstract}

\tableofcontents
\section{Introduction and main results}\label{sec:intro}
We study perturbations of thermodynamic features (topological pressures, Gibbs measures and measure-theoretic entropies of these measures) concerning perturbed potentials defined on topologically transitive subshift of finite type. Our perturbed potential which is given by (\ref{eq:Phe}) tends to $-\infty$ partially. This implies that the subshift with respect to the perturbed systems and the subshift with respect to the unperturbed system are different in general. Such a situation was first considered by Ikawa \cite{Ikawa,Ikawa2} in the study of billiard dynamics problems. This was developed by Morita and Tanaka \cite{MT} who also applied it to the study of degeneration of one-dimensional Markov maps and of asymptotic variance concerning with the central limit theorem. They only considered the case that the subshift with respect to non-perturbed system is mixing essentially. Therefore, it was difficult to consider a more natural perturbation. In this paper, we treat the case where the subshift with respect to the unperturbed system is no topologically transitive, and aim to give a necessary and sufficient condition of convergence of thermodynamic features and to characterize these limit points.

For details, we introduce some notation and notions of thermodynamic formalism below.
Let $d\geq 1$ be an integer and $S=\{1,2,\dots, d\}$ a finite state space with the discrete topology. Denoted by $\Z^{+}$ the all of nonnegative integers and by $S^{\Z^{+}}$ the direct sum of $S$ endowed with the normal product topology. If $\om\in S^{\Z^{+}}$ then $\om$ is written by $\om_{0}\om_{1}\cdots$. For a $d\times d$ one-zero matrix $M=M(ij)$ indexed by $S\times S$, we define $\Sigma_{M}^{+}=\{\om=(\om_{n})_{n=0}^{\infty}\in S^{\Z^{+}}\,:\, M(\om_{n}\om_{n+1})=1\text{ for any } n\geq 0\}$. In general, $\Si_{M}^{+}$ may be an empty set. Assume $\Si_{M}^{+}\neq \emptyset$. The shift transformation $\sigma$ on $S^{\Z^{+}}$ is given by $(\sigma \om)_{n}=\om_{n+1}$ for any $n\geq 0$ and $\om\in S^{\Z^{+}}$. For $\sigma_{M}=\sigma|_{\Si_{M}^{+}}$, the dynamics $(\Sigma_{M}^{+}, \sigma_{M})$ is called by subshift of finite type with transition matrix $M$. In particular, the dynamics $(S^{\Z^{+}}, \sigma)$ is called by a full shift.
A word $w\in S^{n}$ is $M$-admissible if $M(w_{0}w_{1})M(w_{1}w_{2})\cdots M(w_{n-2}w_{n-1})=1$ is satisfied. We set $W_{n}(M)=\{w\in S^{n}\,:\, w \text{ is }M\text{-admissible}\}$ for $n\geq 1$, $W_{0}(M)=\emptyset$, and $W_{*}(M)=\bigcup_{n=0}^{\infty}W_{n}(M)$.
For integers $k\geq 0$ and $n\geq 1$ and a word $w\in S^{n}$, we write ${}_{k}[w]=\{\om\in S^{\Z^{+}}\,;\,\om_{k}\cdots \om_{k+n-1}=w\}$ and ${}_{k}[w]^{M}={}_{k}[w]\cap \Sigma_{M}^{+}$. Such sets are called by cylinder sets on $\Sigma_{M}^{+}$.
For $\theta\in (0,1)$, we define a metric $d_{\theta}$ on $S^{\Z^{+}}$ by $d_{\theta}(\omega,\omega^{\prime})=\theta^{n^\p}$ with
$n^\p=\min\{n\geq 0\,:\,\omega_{n}\neq \omega^{\prime}_{n}\}$.
For $\K=\R$ or $\K=\C$, 
we denote by $C(\Sigma_{M},\K)$ the totally of all $\K$-valued continuous functions on $\Sigma_{M}^{+}$ and by $F_{\theta}(\Sigma_{M}^{+},\K)$ the totally of all $\K$-valued $d_{\theta}$-Lipschitz continuous functions on $\Si_{M}^{+}$. For simplicity, we write $C(\Sigma_{M}^{+})$ as $C(\Sigma_{M},\C)$ and $F_{\theta}(\Sigma_{M}^{+})$ as $F_{\theta}(\Sigma_{M}^{+},\C)$. For $f\in F_{\theta}(\Sigma_{M}^{+})$ and $n\geq 1$, we let $[f]_{\theta,n}=\sup\{|f(\om)-f(\upsilon)|/d_{\theta}(\om,\upsilon)\,:\, \om,\upsilon\in \Sigma_{M}^{+}, \om\neq \upsilon \text{ and }\om\in {}_{0}[\upsilon_{0}\cdots\upsilon_{n-1}]^{M}\}$. For  convenience sake, we set $[f]_{\theta}=[f]_{\theta,1}$ and $[[f]]_{\theta}=[f]_{\theta,2}$. We see that $C(\Sigma_{M}^{+},\K)$ is a Banach space endowed with the supremum norm $\|f\|_{\infty}=\sup_{\om\in \Sigma_{M}^{+}}|f(\om)|$, and $F_{\theta}(\Sigma_{M}^{+},\K)$ a Banach space endowed with the Lipschitz norm $\|f\|_{\theta}=\|f\|_{\infty}+[f]_{\theta}$. We denote by $M(\Sigma_{M}^{+})$ the totally of all Borel probability measures on $\Sigma_{M}^{+}$ and by $M_{\sigma}(\Sigma_{M}^{+})$ the totally of all $\sigma_{M}$-invariant measures belonging to $M(\Sigma_{M}^{+})$.

For $f\in C(\Sigma_{M}^{+},\R)$, it is known in \cite{Bowen} that the limit
$$P(\sigma_{M},f)=\lim_{n\to \infty}\frac{1}{n}\log \sum_{w\in W_{n}(M)}\exp(\sup_{\om\in {}_{0}[w]^{M}}S_{n}f(\om))$$
exists in $\R$ and is called the topological pressure of $f$. Assume that the matrix $M$ is irreducible, namely for any indexes $i,j$ of $M$, there exists an integer $n\geq 0$ such that $M^{n}(ij)>0$. For $\varphi\in C(\Sigma_{M}^{+}, \R)$, a measure $\mu\in M_{\sigma}(\Sigma_{M}^{+})$ is called a Gibbs measure of the potential $\varphi$ if there exist constants $c\geq 1$ and $P\in \R$ such that for all $n\geq 1$ and $\om\in \Sigma_{M}^{+}$,
$$c^{-1}\leq \frac{\mu({}_{0}[\om_{0}\om_{1}\cdots \om_{n-1}]^{A})}{\exp(-nP+S_{n}\varphi(\om)}\leq c$$
holds, where $S_{n}\varphi(\om)=\sum_{k=0}^{n-1}\varphi(\sigma_{M}^{k}\om)$. It is well known that if $\varphi$ is in $F_{\theta}(\Sigma_{M}^{+},\R)$, then there exists an unique ergodic Gibbs measure $\mu_{\ph}$ of the potential $\varphi$. Moreover, the number $P$ is equals to $P(\sigma_{M},\varphi)$. For $m\in M_{\sigma}(\Sigma_{M}^{+})$, the measure-theoretic entropy of $m$ is defined by
$$h(\sigma_{M},m)=\lim_{n\to \infty}\frac{1}{n}\sum_{w\in W_{n}(M)}m({}_{0}[w]^{M})\log m({}_{0}[w]^{M}).$$
By virtue of Variational Principle \cite{Bowen}, for $f\in C(\SiAp,\R)$, the topological pressure $P(\si_{M},f)$ satisfies the equation
\ali
{
P(\si_{M},f)=\sup\{\int f\,dm+h(\si_{M},m)\,:\,m\in M_{\si}(\SiAp)\}.
}
When $f$ is in $F_{\theta}(\Sigma_{M}^{+},\R)$, this supremum is attained at $m=\mu_{f}$
\ali
{
P(\sigma_{M},f)=\int \varphi d\mu+h(\sigma_{M},\mu_{f}).
}

We mention our formulation and main results. Let $A=(A(ij))$ and $B=(B(ij))$ be $d\times d$ zero-one matrices indexed by $S\times S$.
Throughout this paper, we assume that $\SiAp$ is not empty.
We introduce the following conditions for $A$ and $B$:
\begin{itemize}
\item[$(\Si.1)$] $A$ is irreducible, i.e. the dynamics $(\SiAp,\si_{A})$ is topologically transitive.
\item[$(\Si.2)$] $B(ij)=1$ implies $A(ij)=1$.
\item[$(\Si.3)$] $\SiBp$ is not empty.
\end{itemize}
Assume $(\Si.2)$. General theory of nonnegative matrices implies that there is a suitable permutation matrix $P$ so that
\begin{equation}
P^{-1}BP=\left(
\begin{array}{cccc}
B_{11}&B_{12}&\cdots &B_{1m}\\
O      &B_{22}&\ddots & \vdots\\
\vdots &\ddots &\ddots & B_{m-1 m}\\
O     &\cdots &O &B_{mm}
\end{array}\label{eq:tPBP=...2}
\right)
\end{equation}
with an unique integer $m\geq 1$, where each submatrix $B_{kk}$ of $B$ is square and irreducible.
For a square matrix $M$, we define a subset $S_{M}$ by the set of indexes of $M$. Put $\Sigma_{M}=\bigcup_{i\in S_{M}}{}_{0}[i]^{A}$ for a square submatrix $M$ of $A$.
We give the notation
\alil{
\bT=\bT(B)=\{B_{kk}\,:1\leq k\leq m,\,\Si_{B_{kk}}\neq \emptyset\}.\label{eq:bT}
}
Let $\ph\in F_{\theta}(\SiAp,\R)$. We write $\lam(M)=\exp(P(\si_{M},\ph_{M}))$ for each $M\in \bT$ with $\Si_{M}^{+}\neq \emptyset$, where $\ph_{M}=\ph |_{\Si_{M}^{+}}$. If $\Si_{M}^{+}=\emptyset$ then we set $\lam(M)=0$.
Put $\lam=\lam(B,\varphi)=\max_{M\in \bT}\lam(M)$,
\alil{
\bT_{0}&=\bT_{0}(B,\ph)=\{M\in \bT(B)\,:\,\lam(M)=\lam(B,\ph)\}\text{ and}\label{eq:T02}\\
\bT_{1}&=\bT_{1}(B,\ph)=\bT(B)\setminus \bT_{0}(B,\ph).\label{eq:T1}
}
We sometimes call $\Si_{M}^{+}$ ($M\in \bT_{0}$) {\it the transitive component of $\SiBp$ with maximal pressure}.
Let $N=\bigcup_{ij\,:\, B(ij)=0}{}_{0}[ij]^{A}$.
We give the following three conditions for two real-valued functions $\phe,\psi(\e,\cd)\in F_{\theta}(\SiAp,\R)$ with a small parameter $\e>0$:
\begin{itemize}
\item[$(\Ph.1)$] There exists a function $\varphi\in F_{\theta}(\SiAp,\R)$ such that $\|\phe-\varphi\|_{\infty}\to 0$ as $\e\to 0$.
\item[$(\Ph.2)$] $\max_{\om\in N}\psi(\e,\om)\to -\infty$ as $\e\to 0$.
\item[$(\Ph.3)$] $\sup_{\e>0}[[\phe]]_{\theta}<\infty$ and $\sup_{\e>0}[[\psi(\e,\cd)]]_{\theta}<\infty$.
\end{itemize}
We define 
\alil
{
\Phe=\phe+\chi_{N}\psi(\e,\cd),\label{eq:Phe}
}
where $\chi_{N}$ is the identify map of $N$. We obtain one of the main theorems below:
\thm{
\label{th:conv_pressure}
Assume that the conditions $(\Si.1)$-$(\Si.2)$ and $(\Ph.1)$-$(\Ph.2)$ are satisfied. Then if $(\Si.3)$ also holds then $P(\sigma_{A},\Phe)$ converges to $P(\sigma_{B},\varphi_{B})$ as $\e\to 0$. If $(\Si.3)$ does not valid then $P(\sigma_{A},\Phe)$ tends to $-\infty$.
}
Denoted by $\mu(M,\cd)$ the Gibbs measure of the potential $\ph_{M}$ for each $M\in \bT$ with $\Si_{M}^{+}\neq \emptyset$. Next we have the forms of limit points of the Gibbs measure $\mu(\e,\cd)$ of the potentical $\Phe$ and the measure-theoretic entropy $h(\si_{A},\mu(\e,\cd))$ of $\mu(\e,\cd)$:
\thm{\label{th:limpoint_Gibbs_entropy}
Assume that the conditions $(\Si.1)$-$(\Si.3)$ and $(\Ph.1)$-$(\Ph.3)$ are satisfied. Then any accumulation point of the Gibbs measure $\mu(\e,\cd)$ in $M(\SiAp)$ with weak${}^{*}$ topology has the form $\sum_{M\in \bT_{0}}\delta(M)\mu(M,\cd)$ for some nonnegative constant $\delta(M)$ with $\sum_{M\in \bT_{0}}\delta(M)=1$. Moreover, if $\mu(\e,\cd)$ converges to $\sum_{M\in \bT_{0}}\delta(M)\mu(M,\cd)$ weakly, then the entropy $h(\si_{A},\mu(\e,\cd))$ converges to $\sum_{M\in \bT_{0}}\delta(M)h(\si_{M},\mu(M,\cd))$.
}
By virtue of this theorem, when $\sharp \bT_{0}=1$ (write $\bT_{0}=\{M\}$), the Gibbs measure $\mu(\e,\cd)$ converges to $\mu(M,\cd)$ weakly and the entropy $h(\si_{A},\mu(\e,\cd))$ converges to $h(\si_{M},\mu(M,\cd))$.

When $\sharp \bT_{0}\geq 2$, $\mu(\e,\cd)$ and $h(\si_{A},\mu(\e,\cd))$ do not converge as $\e\to 0$ in general. In fact, the relation between the potential $\Phe$ and convergence of $\mu(\e,\cd)$ is very complex and difficult. In this paper, we find a simple relation between convergence of $\mu(\e,\cd)$ and convergence of an expression composed of Perron eigenvalues of generalized Ruelle operators in the case when $\sharp \bT_{0}=2,3$.  To describe the precise statement, we give some notation as follows.
For $M,M^\p\in \bT$, $A_{MM^\p}$ denotes the submatrix of $A$ indexed by $S_{M}\times S_{M^\p}$. For a non-empty subset $\bM\subset \bT$, $A(\bM)$ is defined by the submatrix of $A$ indexed by $\left(\bigcup_{M\in \bM}S_{M}\right)^{2}$. Similarly, $B(\bM)$ is given as a submatrix of $B$. When we write $\bT_{0}=\{M_{1},M_{2},\dots, M_{m_{0}}\}$, we let $\bM(k)=\{M_{k}\}\cup \bT_{1}$ and $\bM(k,k^\p)=\{M_{k},M_{k^\p}\}\cup \bT_{1}$ for $k,k^\p=1,2,\dots, m_{0}$ with $k\neq k^\p$.
Set
\ali{
\lam(\e)&=\exp(P(\si_{A},\Phe)) \text{ for } \e>0,\\
\lam(\bM,\e)&=\exp(P(\si_{A(\bM)},\Phe_{A(\bM)})) \text{ for } \bM\subset \bT \text{ with }\Si_{A(\bM)}^{+}\neq \emptyset.
}
If $\Si_{A(\bM)}^{+}= \emptyset$ then we put $\lam(\bM,\e)=0$.
Note that the number $\lam(\e)$ becomes the eigenvalue of the Ruelle operator of $\Phe$, and $\lam(\bM,\e)$ coincides with the Perron eigenvalue of the generalized Ruelle operator of $\{A(\bM),\Phe\}$ (see Section \ref{sec:Gen.Rue}).
In the case when $\sharp \bT_{0}=2,3$, we can equate convergence of $\mu(\e,\cd)$ with convergence of the number of a representation by differences of these eigenvalues:
\thm{\label{th:m0=2}
Assume that the conditions $(\Si.1)$-$(\Si.3)$ and $(\Ph.1)$-$(\Ph.3)$ are satisfied and $\sharp \bT_{0}=2$, i.e. $\bT_{0}=\{M_{1},M_{2}\}$. Then the number $\delta_{\e}(k)=(\lam(\e)-\lam(\bM(k^\p),\e))/(\sum_{l=1}^{2}(\lam(\e)-\lam(\bM(l),\e))$ converges to a number $\delta(k)$ for each $k=1,2$ with $\{k,k^\p\}=\{1,2\}$ if and only if the Gibbs measure $\mu(\e,\cd)$ converges to a measure $\mu$. In these cases, $\mu$ has the form $\mu=\sum_{k=1}^{2}\delta(k)\mu(M_{k},\cd)$.
}
For $k=1,2,3$ and $\{k,k^\p,k^\pp\}=\{1,2,3\}$, we define
\alil
{
\delta_{\e}(k)&=\frac{\delta^{0}_{\e}(k)}{\sum_{l=1}^{3}\delta^{0}_{\e}(l)} \text{ with}\label{eq:dk_m=3}\\
\delta^{0}_{\e}(k)&=(\lam(\e)-\lam(\bM(k^\p,k^\pp),\e))\times\nonumber\\
&\qquad\qquad(\lam(\e)+\lam(\bM(k^\p,k^\pp),\e)-\lam(\bM(k^\p),\e)-\lam(\bM(k^\pp),\e)).\nonumber
}
\thm
{\label{th:m0=3}
Assume that the conditions $(\Si.1)$-$(\Si.3)$ and $(\Ph.1)$-$(\Ph.3)$ are satisfied and $\sharp \bT_{0}=3$, i.e. $\bT_{0}=\{M_{1},M_{2},M_{3}\}$. Then 
the vector $\delta_{\e}(k)$ defined by (\ref{eq:dk_m=3}) converges to a number $\delta(k)$ for all $k=1,2,3$ if and only if
the Gibbs measure $\mu(\e,\cd)$ converges to a measure $\mu$. In these cases, $\mu$ has the form $\sum_{k=1}^{3}\delta(k)\mu(M_{k},\cd)$.
}
\rems
\item[(1)] Theorem \ref{th:m0=2} and Theorem \ref{th:m0=3} are proven by using the fact from which convergence of the Gibbs measure $\mu(\e,\cd)$ is reduced to convergence of the right Perron eigenvector of a $\sharp \bT_{0}\times \sharp \bT_{0}$ nonnegative irreducible matrix $\tmV_{\bbM,\e}$ (Theorem \ref{th:conv_gibbs_iff_mVM}). In addition, it plays the essential role in the proof of Theorem \ref{th:m0=3} that the difference of the eigenvalues of generalized Ruelle operators and the difference of the eigenvalues of the submatrices of $\tmV_{\bbM,\e}$ going to $0$ at the same speed (Lemma \ref{lem:lam_v-l/lam-l->1_m0>=3}(2) and Lemma \ref{lem:le-lvie/le-lje->1_m0=3}).

On the other hand, it is a natural question whether similar theorems follow when $\sharp \bT_{0}\geq 4$. Unfortunately, such key lemmas are not satisfied at $\sharp \bT_{0}=4$ (Section \ref{sec:m0>=4}) in general. Therefore, another approaches are necessary for knowing whether similar assertions are valid under the case $\sharp \bT_{0}\geq 4$.
\item[(2)] By a general theory of measure-theoretic entropies, if a parametrized $\sigma$-invariant measure $m_{e}\in M_{\si}(\SiAp)$ converges to a measure $m$ in $M(\SiAp)$ weakly, then $\limsup_{\e\to 0}h(\si_{A},m_{e})\leq h(\si_{A},m)$ is valid (e.g. Theorem 8.2 in \cite{Walters}) and it doesn't mean convergence of this entropy. 
Theorem \ref{th:limpoint_Gibbs_entropy} states that if $\mu(\e,\cd)$ converges, then so is for this entropy $h(\si_{A},\mu(\e,\cd))$.
\item[(3)] Morita and Tanaka \cite{MT} gave the following conditions stronger than $(\Si.1)$, $(\Si.3)$ and $(\Ph.3)$:
\ite{
\item[$(\Si.1)^\p$] $A^{n_{0}}>0$ for some $n_{0}\geq 1$, i.e. $(\SiAp,\si_{A})$ is topologically mixing.
\item[$(\Si.3)^\p$] The set $\Sigma_{B}^{+}$ is not empty and  if  $\Sigma^{+}$ is the maximal $\sigma$-invariant subset of $\Sigma_{B}^{+}$, then the subshift $(\Sigma^{+},\sigma|_{\Sigma^+})$ is topologically mixing. 
\item[$(\Ph.3)^\p$] $\sup_{\e>0}[\ph(\e,\cd)]_{\theta}<\infty$ and $\sup_{\e>0}[\psi(\e,\cd)]_{\theta}<\infty$.
}
Note that the condition $(\Si.3)^{\p}$ is satisfied if and only if $\bT_{0}$ consists of only one element $M$, the subshift $(\Si_{M}^{+},\si_{M})$ is topologically mixing, and for each $M^\p\in \bT_{1}$, $M^{\p}$ is a $1\times 1$ zero matrix.
They showed that under the six conditions $(\Si.1)^\p,(\Si.2),(\Si.3)^\p,(\Ph.1),(\Ph.2)$ and $(\Ph.3)^\p$, the three thermodynamic futures $P(\si_{A},\Phe)$, $\mu(\e,\cd)$, $h(\si_{A},\mu(\e,\cd))$ converge (Theorem 1.1 in \cite{MT} ).
\item[(4)] In \cite{T2012}, we showed that in addition to the conditions $(\Si.1)$-$(\Si.3)$ and $(\Ph.1)$-$(\Ph.3)$, when we assume differentiability conditions for the potential $\Phe$ at $\e=0$, a semisimplicity for $\lam(B,\ph)$ and some strong conditions for $B$ and $\ph$, the Gibbs measure $\mu(\e,\cd)$ and the entropy of this measure converge.
\item[(5)] The general analytic perturbation theory provides with a (Puiseux) series expansion of an isolated eigenvalue with finite multiplicity of a bounded linear operator and a (Puiseux) series expansion of the corresponding eigenvector \cite{Kato}. 
Therefore, when the potential $\Phe$ is analytically expanded in $(F_{\theta}(\SiAp,\R), \|\cdot\|_{\theta})$, so is for the Ruelle operator $\LR_{A,\Phe}\in \LR(F_{\theta}(\SiAp))$ and consequently thermodynamic futures are continuous at $\e=0$, where these notation appear in Section \ref{sec:RPFth}. 
\item[(6)] The condition $(\varPhi.3)$ implies Lasota-Yorke type inequality for the Ruelle operator of $\Phe$ uniformly in $\e>0$. Indeed, assume $(\Si.1)$ and $(\Ph.3)$ (other conditions do not need). Recall that the Ruelle operator $\LR_{A,\Phe}$ of $\Phe$ is a bounded linear operator acting on $C(\SiAp)$ or $F_{\theta}(\SiAp)$ which is defined as (\ref{eq:LR=}) by putting $M=A$ and $\ph=\Phe$.
Let 
\alil
{
\tPhe&=\varPhi(\e,\cdot)-\log g(\e,\cd)\circ \sigma_{A}+\log g(\e,\cd)\label{eq:tPhe=}\\
\hat{\Ph}(\e,\cd)&=\tPhe-\log \lam(\e),\label{eq:hatPhe=}
}
where $g(\e,\cd)\in F_{\theta}(\Sigma_{A}^{+})$ is the Perron eigenfunction of $\LR_{A,\Phe}$ with $\|g(\e,\cd)\|_{\infty}=1$. 
In this setting, the equality
\alil{
\LR_{A,\hat{\Ph}(\e,\cd)}1=1\label{eq:LR1=1}
}
holds for any $\e>0$. For $c>0$, we define
\alil{
\Lambda_{c}=&\{f\in C(\SiAp)\,:\,0\leq f\leq 1\text{ and }f(\om)\leq e^{cd_{\theta}(\om,\upsilon)}f(\upsilon)\nonumber\\
&\qquad\qquad\qquad\qquad\qquad\qquad \text{ for }\om,\upsilon\in \SiAp \text{ with }\om_{0}=\upsilon_{0}\}.\label{eq:Lc=}
}
Note that this set is a subset of $F_{\theta}(\SiAp,\R)$ and a compact subset of $C(\SiAp)$ for each $c>0$. 
It is well-known that $g(\e,\cd)$ belongs to $\Lambda_{c}$ with the constant $c=\theta\sup_{\e^\p>0}[[\Phi(\e^\p,\cd)]]_{\theta}/(1-\theta)$ for any $\e>0$ \cite{MT}. Therefore, $[\log g(\e,\cd)]_{\theta}\leq c^\p$ for some constant $c^\p>0$. Moreover, it is not hard to see that there exist constants $c^\pp>0$ such that for any $\e>0, n\geq 1$ and $f\in F_{\theta}(\Sigma_{A}^{+})$,
\alil{
[\LR_{A,\hat{\Ph}(\e,\cd)}^{n}f]_{\theta}\leq \theta^{n}[f]_{\theta}+c^\pp\|f\|_{\infty}.\label{eq:LY}
}
by using a similar argument of Proposition 2.1 in \cite{PP}.
On the other hand,  $\LR_{A,\hat{\Ph}(\e,\cd)}\in \LR(C(\SiAp))$ has a convergence subsequence $(\e_{n})$ and a limit point $\hat{\LR}$ from Proposition \ref{prop:conv_nor_TSC_M}(1) and the equation (\ref{eq:LR1=1}). Thus it follows from these results that for a sufficiently small enough $\eta>0$, the eigenprojection
\ali{
\int_{\partial B(1,\eta)}(z\IR-\LR_{A,\hat{\Ph}(\e,\cd)})^{-1} dz
}
converges to the eigenprojection
\ali{
\int_{\partial B(1,\eta)}(z\IR-\ti{\LR})^{-1} dz
}
as $\e\to 0$ running through $(\e_{n})$ \cite{KL}, where $B(a,\eta)$ is the open ball in $\C$ with a center $a$ and a radius $\eta$, and $\IR$ is the identity operator belonging to $\LR(C(\SiAp))\cap \LR(F_{\theta}(\SiAp))$.
However even if $\LR_{A,\hat{\Ph}(\e,\cd)}$ converge, this fact does not imply convergence of the (geometric) eigenfunction $g(\e,\cd)$, convergence of the Perron eigenvector of the dual $\LR_{A,\Phe}^{*}$ of $\LR_{A,\Phe}$ and convergence of the Gibbs measure $\mu(\e,\cd)$.
\medskip
\par
In Section \ref{sec:pre}, we introduce the notion which is often used throughout this paper. 
In the beginning of this section, we give a notation of asymptotic relation. After that, we show an expansion formula of M-matrices by using the eigenvalues of submatrices ( Proposition \ref{prop:decom_eigen_Mmat}). Furthermore, we recall spectral properties of the generalized Ruelle operators which are introduced in \cite{T2009} and are basic tools to prove the main theorems. Note that the topological pressure, the Gibbs measure and the measure-theoretic entropy are expressible in terms of the eigenvalues and eigenspaces of this operators. In Section \ref{sec:mat_bdd}, we mention an abstract frame that reduce the perturbation problem of the (generalized) Ruelle operators to the perturbation problem of certain irreducible nonnegative matrices. In particular, we shall reduce convergence of the Gibbs measure $\mu(\e,\cd)$ to convergence of the Perron eigenvector of a nonnegative irreducible matrix (Theorem \ref{th:conv_gibbs_iff_mVM}). 
Section \ref{sec:proof} is devoted to proofs of main results. Under the case when $\sharp \bT_{0}=3$, we show that the ratio between the difference of the eigenvalues of generalized Ruelle operators and the difference of the eigenvalues of the submatrices of $\tmV_{\bbM,\e}$ converges to $1$ as $\e\to 0$ (Lemma \ref{lem:lam_v-l/lam-l->1_m0>=3}(2) and Lemma \ref{lem:le-lvie/le-lje->1_m0=3}). These facts play an important role in the proof of Theorem \ref{th:m0=3}. In the final section, we demonstrate that such facts do not follow in general under the case $\sharp \bT_{0}=4$ (Section \ref{sec:m0>=4}). Moreover, we illustrate the relation between the potential $\Phe$ and convergence of the Gibbs measure $\mu(\e,\cd)$ under the case $\sharp \bT_{0}=2$ by using asymptotic expansion techniques for eigenvalues of Ruelle transfer operators (Section \ref{sec:ex_asymp}).
\reme
\noindent
Acknowledgements:\\
The author is heartily grateful to Katsukuni Nakagawa, Hiroshima University, for many valuable conversations.
\section{Preliminaries}\label{sec:pre}
\subsection{Asymptotic relation}
Let $f_{\e}$ and $g_{\e}$ be two nonnegative valued functions defined on a set $X$ with a small parameter $\e>0$. We consider the condition that there exist constants $c\geq 1$ and $\epsilon_{0}>0$ such that for any $x\in X$ and for any $0<\e\leq \epsilon_{0}$
$$c^{-1}g_{\e}(x)\leq f_{\e}(x)\leq c g_{\e}(x)$$
holds. We then write $f_{\e} \asymp g_{\e}\ (\e \to 0)$, or simply,  $f_{\e} \asymp g_{\e}$. We see that $\asymp$ is an equivalent relation. The following basic proposition is useful in some proofs.
\prop{\label{prop:a+b/c+d}
Assume that $a_{\e}(i)$ and $b_{\e}(i)$ are nonnegative numbers with a small parameter $\e>0$ for $i=1,2,\dots, n$ satisfying that there exists $c_{\e}(i)$ such that $a_{\e}(i)=c_{\e}(i)b_{\e}(i)$ and $c_{\e}(i)\asymp 1$ as $\e\to 0$ for each $i=1,2,\dots, n$. Assume also $a_{\e}(n)>0$ and $b_{\e}(n)>0$. Then
\begin{align*}
d_{\e}=\frac{a_{\e}(1)+a_{\e}(2)+\cdots+a_{\e}(n)}{b_{\e}(1)+b_{\e}(2)+\cdots+b_{\e}(n)}\asymp 1
\end{align*}
 as $\e\to 0$. In particular, if $c_{\e}(i)\to 1$ for each $i$, then $d_{\e}\to 1$.
}
\pros
Let $\underline{c}_{\e}=\min_{i}c_{\e}(i)$ and $\overline{c}_{\e}=\max_{i}c_{\e}(i)$. Since $d_{\e}$ is a number between $\underline{c}_{\e}$ and $\overline{c}_{\e}$, we obtain the assertion.
\proe
\subsection{Some properties of nonnegative irreducible matrices}\label{sec:tensor}
We start with elementary results for M-matrices.
\subsubsection{Relation between M-matrices and eigenvectors}\label{sec:re_Mmat_evec}
Let $T$ be a finite set and $M=(M(ij))$ a $\sharp T\times \sharp T$ nonnegative irreducible matrix indexed by $T\times T$. Perron-Frobenius Theorem provides the right Perron eigenvector ${}^{t}(b(i))$ and the left Perron eigenvector $(c(i))$ with $\sum_{i\in T}b(i)=\sum_{i\in T}b(i)c(i)=1$. Denoted by $I$ an identity matrix. For a square nonnegative matrix $M_{0}$, let $\eta(M_{0})$ be the Perron eigenvalue of $M_{0}$.
We write the adjoint matrix of $\eta(M) I-M$ by $(M_{\eta(M)}(ji))$ whose elements are given by
$M_{\eta(M)}(ji)=(-1)^{j+i}\det D_{ji}$ and $D_{ji}$ is a minor of order $n-1$ of $\eta(M) I-M$ is defined to be a submatrix of $\eta(M) I-M$ obtained by striking out $j$th row and $i$th column.

We often call an element $\bi=(i_{j})_{j=1}^{p}=i_{1}i_{2}\dots i_{p}\in T^{p}$ a {\it path} on $T$ from $i_{1}$ to $i_{p}$. 
This path $\bi$ is a {\it simple path} on $T$ from $i_{1}$ to $i_{p}$ if $i_{1}$, $i_{2}$,\dots, $i_{p}$ are distinct. We give
\begin{align}
\SP(i, j : T^{\prime})=\{\bi\in \bigcup_{n=1}^{\infty}T^{n} :  \bi\text{ is a simple path on }T^{\p}\cup\{i,j\} \text{ from }i \text{ to }j\}\label{eq:SP}
\end{align}
for a subset $T^{\prime}\subset T$ and for each $i,j\in T$.
\prop{[\cite{Seneta}]\label{prop:cm=cb}
\ite{
\item $c(j)b(i)=M_{\eta(M)}(ji)/\sum_{k\in T}M_{\eta(M)}(kk)$ for each $i,j\in T$.
\item $b(i)/b(k)=M_{\eta(M)}(ji)/M_{\eta(M)}(jk)$ for each $i,k,j\in T$.
\item $c(j)/c(k)=M_{\eta(M)}(ji)/M_{\eta(M)}(ki)$ for each $j,k,i\in T$.
}
}
\pros
The assertion (1) follows from \cite{Seneta} (Corollary 2, p 9). The assertions (2) and (3) are obtained by (1).
\proe
For a simple path $\bi=(i_{j})_{j=1}^{p}$ on $T$, we define a submatrix $M(\bi)$ of $M$ by
\alil{
M(\bi)=M(i_{1},i_{2},\dots, i_{p})=\left(\begin{array}{cccc}
M(i_{1}i_{1}) &M(i_{1}i_{2})  &\cdots  &M(i_{1}i_{p})\\
M(i_{2}i_{1}) &M(i_{2}i_{2})  &\cdots  &M(i_{2}i_{p})\\
\vdots &\vdots  &\ddots  &\vdots \\
M(i_{p}i_{1}) &M(i_{p}i_{2})  &\cdots  &M(i_{p}i_{p})\\
\end{array}\right).\label{eq:M(bi)=}
}
We have the following:
\prop{\label{prop:decom_eigen_Mmat}
For each $k\in T$, there exist integers $q(k)\geq 1$ and $p(k,i)\geq 1$ $(1\leq i\leq q(k))$, and simple paths $\bi(k,i,j), \bi^{\prime}(k,i,j)$ on $T$ such that $\bi^{\prime}(k,i,j)\subsetneq \bi(k,i,j)$ by regarding as sets, 
and $M_{\eta(M)}(kk)$ has the form
$$M_{\eta(M)}(kk)=\sum_{i=1}^{q(k)}\prod_{j=1}^{p(k,i)}(\eta(M(\bi(k,i,j)))-\eta(M(\bi^{\prime}(k,i,j)))).$$
}
\pros
We write $T=\{i_{1},i_{2},\dots, i_{m}\}$. Put $M^{(k)}=M(i_{1},\dots, i_{k-1},i_{k+1},\dots, i_{m})$ for $k=1,2,\dots, m$. For $\xi=(\xi_{j})\in \R^{p}$, we denote by $\dia (\xi)$ the diagnostic matrix whose diagonal elements are $\xi_{1},\xi_{2},\dots, \xi_{p}$.
Consider the following steps (a) and (b) by the first defining as $D:=M^{(1)}$, $q:=m-1$ and $\xi:=(\xi_{1},\xi_{2},\dots, \xi_{q})$ with $\xi_{j}=\eta(M)$ for $j=1,2,\dots q$:
\ite{
\item[(a)] 
In the case when $q=1$,
$$\det(\dia(\xi)-D)=\xi_{1}-D(11)=\xi_{1}-\eta(D(i_{1}))$$
is satisfied.

In the case when $q\geq 2$, let ${\bf m}_{1}, {\bf m}_{2}, \dots, {\bf m}_{q}$ be the columns of $D$. By using the unit vectors ${\bf e}_{1}, {\bf e}_{2}, \dots, {\bf e}_{q}$, we have
\begin{align*}
&\ \ \ \det(\dia(\xi)-D)=\det(\xi_{1}{\bf e}_{1}-{\bf m}_{1}\ \ \xi_{2}{\bf e}_{2}-{\bf m}_{2}\cdots\xi_{q}{\bf e}_{q}-{\bf m}_{q})\\
&=\sum_{i=1}^{q}\Bigg(\det\Big(\underbrace{\eta(D){\bf e}_{1}-{\bf m}_{1}\cdots \eta(D){\bf e}_{i-1}-{\bf m}_{i-1}}_{i-1}\ \ \underbrace{\xi_{i}{\bf e}_{i}-{\bf m}_{i}\cdots \xi_{q}{\bf e}_{q}-{\bf m}_{q}}_{q-i+1}\Big)\\
&\ \ -\det\Big(\underbrace{\eta(D){\bf e}_{1}-{\bf m}_{1}\cdots\eta(D){\bf e}_{i}-{\bf m}_{i}}_{i}\ \ \underbrace{\xi_{i+1}{\bf e}_{i+1}-{\bf m}_{i+1}\cdots \xi_{q}{\bf e}_{q}-{\bf m}_{q}}_{q-i}\Big)\Bigg)\\
&=\sum_{i=1}^{q}(\xi_{i}-\eta(D))\det\Big(\underbrace{\eta(M_{j}){\bf e}_{1}-{\bf m}_{1}\cdots \eta(M_{j}){\bf e}_{i-1}-{\bf m}_{i-1}}_{i-1}\\
&\qquad\qquad\qquad\qquad\qquad\qquad\qquad\qquad\underbrace{\xi_{i+1}{\bf e}_{i+1}-{\bf m}_{i+1}\cdots \xi_{d_{j}}{\bf e}_{d_{j}}-{\bf m}_{q}}_{q-i}\Big)\\
&=\sum_{i=1}^{q}(\xi_{i}-\eta(D))\det(\dia({\bf \xi}^{\p}(i))-D^{\p}(i)),
\end{align*}
where each ${\bf \xi}^{\p}(i)$ is in $\R^{q-1}$ and each $D^{\p}(i)$ is a $(q-1)\times (q-1)$ nonnegative matrix.
\item[(b)] We repeat the step (a) as $D:=D^{\p}(i)$ and ${\bf \xi}:={\bf \xi}^{\p}(i)$. 
}
By using the steps $(a)$ and $(b)$ repeatedly, we see that the assertion is valid.
\proe
\rems
\label{rem:cb=}
Denote by $T=\{i_{1},i_{2},\dots, i_{m}\}$.
\item[(1)] In the case when $M$ is a $2\times 2$ matrix, we have
\begin{align}
M_{\eta(M)}(i_{1}i_{1})=\eta(M)-\eta(M(i_{2})),\quad M_{\eta(M)}(i_{2}i_{2})=\eta(M)-\eta(M(i_{1}))
\end{align}
and $b(i_{k})c(i_{k})=M_{\eta(M)}(i_{k}i_{k})/\sum_{l=1}^{2}M_{\eta(M)}(i_{l}i_{l})$ for $k=1,2$, by using Proposition \ref{prop:cm=cb}(1).
\item[(2)] In the case when  $M$ is a $3\times 3$ matrix, we obtain
\ali
{
&M_{\eta(M)}(i_{1}i_{1})\\
=&(\eta(M)-\eta(M(i_{2},i_{3})))(\eta(M)-\eta(M(i_{2}))+\eta(M(i_{2},i_{3}))-\eta(M(i_{3})))
}
for $T=\{i_{1},i_{2},i_{3}\}$
and $b(i_{k})c(i_{k})=M_{\eta(M)}(i_{k}i_{k})/\sum_{l=1}^{3}M_{\eta(M)}(i_{l}i_{l})$ for $k=1,2,3$.

\reme
\subsubsection{Behaviours of Perron eigenvectors of perturbed nonnegative matrices}\label{sec:behav_evec_per_mat}
Let $T$ be a finite totally ordered set and
$M_{\e}=(M_{\e}(ij))$ be nonnegative irreducible matrix indexed by $T\times T$ with a small parameter $\e>0$ satisfying the following (M.1)-(M.4):
\begin{itemize}
\item[(M.1)] For any $i,j\in T$, either $M_{\e}(ij)>0$ for all $\e>0$ or $M_{\e}(ij)=0$ for all $\e>0$.
\item[(M.2)] There exists a constant $c>0$ such that $M_{\e}(ij)\leq c$ for any $\e>0$ and $i<j$.
\item[(M.3)] $M_{\e}(ij)$ vanishes as $\e\to 0$ for any $i>j$.
\item[(M.4)] $M_{\e}(ii)$ converges to a number $\eta_{i}\geq 0$ as $\e\to 0$ for any $i$.
\end{itemize}
Put $\eta_{\e}(i)=M_{\e}(ii)$ for any $\e>0$ and $i\in T$. Give
\begin{align}
T_{0}=\{i\in T\,:\,\eta_{i}=\max_{j\in T}\eta_{j}\},\quad T_{1}=T\setminus T_{0}.
\end{align}
Furthermore, there exists $\e_{0}>0$ such that for any $0<\e\leq \e_{0}$, the inequality $\eta_{\e}(i)<\max_{j\in T}\eta_{j}$ holds for $i\in T_{1}$.

Denoted by $\eta(\e)$ the Perron eigenvalue of $M_{\e}$, by ${}^{t}(b_{\e}(i))$ the right Perron eigenvector of $M_{\e}$ with $\sum_{i\in T}b_{\e}(i)=1$ and by $(c_{\e}(i))$ the left Perron eigenvector of $M_{\e}$ with $\sum_{i\in T}b_{\e}(i)c_{\e}(i)=1$. Then we have the following which will be used in Section \ref{sec:mat_bdd}.
\prop{\label{prop:convspd_evec_T1}
Under the above notation, for any $i\in T_{1}$
\alil
{
b_{\e}(i)&\asymp \sum_{j\in T_{0}}\sum_{w\in \SP(i,j:T_{1})}M_{\e}(w)b_{\e}(j)\label{eq:convspd_evec_T1_b}\\
c_{\e}(i)&\asymp \sum_{j\in T_{0}}\sum_{w\in \SP(j,i:T_{1})}M_{\e}(w)c_{\e}(j)\label{eq:convspd_evec_T1_c}.
}
}
\pros
Put Let $L$ be a $\sharp T_{1}\times \sharp T_{1}$ zero-one matrix indexed by $T_{1}\times T_{1}$ so that $L(jk)=1$ for $j\neq k$ and $L(jj)=0$ for $j$. For $i\in T_{1}$,
\alil
{
\eta(\e)b_{\e}(i)&=\sum_{j\in T_{0}}M_{\e}(ij)b_{\e}(j)+\sum_{j\in T_{1}\setminus\{i\}}M_{\e}(ij)b_{\e}(j)+\eta_{i}b_{\e}(i).\label{eq:etaebe=}
}
Put $a_{\e}=(a_{\e}(i))_{i\in T_{1}}$ with $a_{\e}(i)=\sum_{j\in T_{0}}(\eta(\e)-\eta_{i})^{-1}M_{\e}(ij)b_{\e}(j)$ and $L_{\e}=(L_{\e}(ij))$ be a matrix indexed by $T_{1}\times T_{1}$ so that $L_{\e}(ij)=(\eta(\e)-\eta_{i})^{-1}M_{\e}(ij)$ for $i\neq j$ and $L_{\e}(ii)=0$. The equation (\ref{eq:etaebe=}) yields
$b_{\e}(i)=a_{\e}(i)+L_{\e}b_{\e}(i)$.
By using this equation repeatedly, we obtain
\begin{align*}
b_{\e}(\cd)=(I+L_{\e}+\cdots+L_{\e}^{nm_{1}-1})a_{\e}(\cd)+L_{\e}^{nm_{1}}b_{\e}(\cd)
\end{align*}
with $m_{1}=\sharp T_{1}$. Note that for any $w=w_{1}\cdots w_{m_{1}+1}\in W_{m_{1}+1}(L)$, $w$ contains a cycle, i.e. $w_{j}=w_{p}$ for some $1\leq j< p\leq m_{1}+1$. Therefore $L_{\e}(w)=\prod_{k=1}^{m_{1}}L_{\e}(w_{k}w_{k+1})$ tends to $0$ by the conditions (M.2)-(M.4).
For a sufficiently small $\e>0$, $\|L_{\e}^{m_{1}}\|<1$ is satisfied, where $\|\cd\|$ is a matrix norm. We see 
\begin{align}
|L_{\e}^{nm_{1}}b_{\e}(\cd)|\leq \|L_{\e}^{m_{1}}\|^{n}|b_{\e}(\cd)|\to 0
\end{align}
as $n\to \infty$.
Thus we have the form
\alil{
b_{\e}(i)=a_{\e}(i)+\sum_{n=0}^{\infty}\sum_{j\in T_{1}}\sum_{w\in W_{n}(L)}L_{\e}(i\cd w\cd j)a_{\e}(j).\label{eq:convspd_evec_T1}
}
Let $d_{\e}=(d_{\e}(i))$ with $d_{\e}(i)=\sum_{j\in T_{0}}\sum_{w\in \SP(i,j:T_{1})}M_{\e}(w)b_{\e}(j)$. Note that for $w\in \SP(i,j:T_{1})$, $w=i\cd u\cd j$ for some $u\in W_{k}(L)$ for some $0\leq k\leq m_{1}$. Now we show $b_{\e}(i)\asymp d_{\e}(i)$ for $i\in T_{1}$. Since $a_{\e}$ and $M_{\e}$ are nonnegative,
\begin{align*}
b_{\e}(i)&\geq (I+L_{\e}+\cdots+L_{\e}^{m_{1}-1})a_{\e}(i)\\
&=a_{\e}(i)+\sum_{l\in T_{1}}\sum_{j=1}^{m_{1}}\sum_{u\in W_{j-1}(L)}L_{\e}(i\cd u\cd l)\sum_{s\in T_{0}}(\eta(\e)-\eta_{l})^{-1}M_{\e}(ls)b_{\e}(s)\\
&\geq c\sum_{s\in T_{0}}\sum_{j=0}^{m_{1}}\sum_{u\in W_{j}(L)}M_{\e}(i\cd u\cd s)b_{\e}(s)\geq cd_{\e}(i)
\end{align*}
is satisfied for any small $\e>0$ with $c=\min_{0\leq k\leq m_{1}}(\max_{j\in T_{1}}(\max_{k}\eta_{k}-\eta_{j}))^{-k}/2>0$.

To prove the converse inequality,
let
$$r(\e)=\max_{v\in W_{*}(L)\,:\,\text{simple cycle}}M_{\e}(v).$$
Then we see $r(\e)\to 0$ as $\e\to 0$. Take $\e_{0}>0$ and $r\in (0,1)$ so that $m_{1}r(\e)^{1/m_{1}}\leq r$ for $0<\e<\e_{0}$.
Recall the form (\ref{eq:convspd_evec_T1}).
For any $w\in W_{n}(L)$ with $i\cd w\cd j\in W_{n+2}(L)$, there exist $u=u(w)\in W_{*}(L)$ and simple cycles $v_{1},\dots,v_{k}\in W_{*}(L)$ such that $|u|\leq m_{1}-2$, $i\cd u\cd j$ is simple path, and 
\ali{
L_{\e}(i\cd w\cd j)=L_{\e}(i\cd u\cd j)\prod_{s=1}^{k}L_{\e}(v_{s})
}
with $n=|w|=|u|+\sum_{s=1}^{k}(|v_{s}|-1)\leq (k+1)m_{1}$, where $|w|$ is the length of $w$. Therefore, 
\begin{align*}
\sum_{j\in T_{1}}L_{\e}(i\cd w\cd j)a_{\e}(j)=&\sum_{j\in T_{1}}L_{\e}(i\cd u\cd j)a_{\e}(j)\prod_{s=1}^{k}L_{\e}(v_{s})\\
\leq& d_{\e}(i)r(\e)^{k}\leq d_{\e}(i)r(\e)^{n/m_{1}-1}.
\end{align*}
Thus we see by (\ref{eq:convspd_evec_T1})
\begin{align*}
b_{\e}(i)\leq a_{\e}(i)+d_{\e}(i)\sum_{n=0}^{\infty}(m_{1})^{n}r(\e)^{n/m_{1}}\leq d_{\e}(i)\left(1+\frac{1}{1-r}\right).
\end{align*}
The proof of the assertion (\ref{eq:convspd_evec_T1_b}) is complete. We also obtain the assertion (\ref{eq:convspd_evec_T1_c}) by considering the transposed matrix of $M_{\e}$.
\proe
\subsection{Generalized Ruelle operators}\label{sec:Gen.Rue}
\subsubsection{Ruelle-Perron-Frobenius type theorem}\label{sec:RPFth}
Denoted by $\LR(\mathcal{X})$ the totality of all bounded linear operators acting on a Banach space $\mathcal{X}$.
In this section, we recall an analogous Ruelle operators acting on $C(\SiAp)$ or on $F_{\theta}(\SiAp)$ (\cite{T2009}). Let $S^{\prime}, S^{\prime\prime}\subset S$ be non-empty subsets and $M=(M(ij))$ a $\sharp S^{\prime}\times \sharp S^{\prime\prime}$ zero-one matrix indexed by $S^{\prime}\times S^{\prime\prime}$ satisfying that $M(ij)=1$ implies $A(ij)=1$.
For $\varphi\in C(\SiAp,\R)$, we define $\LR_{M,\varphi}\in \LR(C(\SiAp))$ by
\alil{
\LR_{M,\varphi}f(\om)=\sum_{i\in S\,:\,M(i\om_{0})=1}e^{\varphi(i\cd \om)}f(i\cd \om)\label{eq:LR=},
}
where $i\cd\om$ is the concatenation of $i$ and $\om$, i.e. $i\cd\om=i\om_{0}\om_{1}\cdots$, and where if $i$ or $\om_{0}$ is not an index of $M$, then $M(i\om_{0})$ regards as $0$. For $\varphi\in F_{\theta}(\SiAp,\R)$, $\LR_{M,\varphi}\in \LR(F_{\theta}(\SiAp))$ is well-defined. We call this operator a generalized Ruelle operator of $\{M,\ph\}$. 
The dual operator $\LR_{M,\varphi}^{*}\in \LR(M(\SiAp))$ of $\LR_{M,\ph}$ is defined by $\LR_{M,\varphi}^{*}m(f)=m(\LR_{M,\varphi}f)$ for $m\in M(\SiAp)$ and $f\in C(\SiAp)$. 

The following is an analogues of Ruelle-Perron-Frobenius Theorem:
\thm{[\cite{MT},\cite{T2009}]\label{th:aRPF}
Assume that the conditions $(\Si.1)$-$(\Si.3)$ are satisfied. Let $\varphi\in F_{\theta}(\SiAp,\R)$ and $\bM$ a subset of $\bT$ with $\sharp \bT_{0}(B(\bM),\ph)=1$. Then there exists an unique triplet $(\lam,h,\nu)\in \R\times F_{\theta}(\SiAp)\times M(\SiAp)$ such that $\lam$ is the maximal simple eigenvalue of the operator $\LR_{B(\bM),\varphi}\in \LR(C(\SiAp))\cap \LR(F_{\theta}(\SiAp))$, $h$ is a nonnegative eigenfunction corresponding to $\lam$, and $\nu$ is an eigenvector corresponding to $\lam$ of the dual $\LR_{B(\bM),\varphi}^{*}$ with $\nu(h)=1$. In particular, $h/\|h\|_{\infty}$ is in $\Lambda_{c}$ with $c=\theta[[\varphi]]_{\theta}/(1-\theta)$ is satisfied.
}
Denote $\bT_{0}(B(\bM),\ph)=\{M\}$. We also see
\alil
{
\supp h&=\bigcup\{{}_{0}[j]^{A}\,:\,B(\bM)^{n}(ij)>0\text{ for some }n\geq 0,\ i\in S_{M}\}\nonumber\\
\supp \nu&=\Si_{B(\bM)}^{+}\cap\bigcup\{{}_{0}[j]^{A}\,:\,B(\bM)^{n}(ji)>0\text{ for some }n\geq 0,\ i\in S_{M}\}.\label{eq:suppnu=}
}
It is known that the topological pressure $P(\sigma_{B(\bM)},\varphi_{B(\bM)})$ of $\varphi_{B(\bM)}=\varphi|_{\Sigma_{B(\bM)}^{+}}$ is equals to $\log\lam$, $\mu=h\nu$ is in $M_{\si}(\SiAp)$ with $\mathrm{supp}\,\mu=\Sigma_{M}^{+}$, and the restriction $\mu|_{\Sigma_{M}^{+}}$ of this measure becomes the Gibbs measure of $\varphi_{M}$. For the sake of convenience, we call the triplet $(\lam,h,\nu)$ the thermodynamic spectral characteristic (or TSC for a short) of the operator $\LR_{B(\bM),\varphi}$. We sometimes write $(\lam,h,\nu)$ by
\alil{
(\lam_{B(\bM),\ph},h_{B(\bM),\ph},\nu_{B(\bM),\ph}).\label{eq:TSC}
}
Put
\alil{
g_{B(\bM),\ph}=h_{B(\bM),\ph}/\|h_{B(\bM),\ph}\|_{\infty}.\label{eq:g}
}

Finally, we consider a behaviour of the operator $\LR_{A_{MM^\p},\Phe}$ for each $M,M^\p\in \bT$. Assume the conditions $(\Si.2)$ and $(\Ph.1)$-$(\Ph.3)$.
Let 
\begin{align}
t_{\e}(MM^\p)=\max_{\om\in \Sigma_{MM^{\prime}}}e^{\varPhi(\e,\om)}\label{eq:tkkp=}
\end{align}
with 
\alil
{
\Sigma_{MM^{\prime}}=\{\om\in \SiAp\,:\,\om_{0}\in S_{M},\ \om_{1}\in S_{M^{\prime}}\}\label{eq:SiMMp=}
}
for $M,M^\p\in \bT$. We note that $t_{\e}(MM^\p)=0$ iff $\Si_{MM^\p}=\emptyset$ iff $A_{MM^\p}=O$.
\prop{
\label{prop:conv_TSC_kkp}
Assume that the conditions $(\Si.2)$ and $(\Ph.1)$-$(\Ph.3)$ are satisfied. Then for any $M,M^\p\in \bT$ with $A_{MM^\p}\neq O$, $\LR_{A_{MM^\p},\Phe}/t_{\e}(MM^\p)$ has a convergent subsequence in $\LR(C(\SiAp))$ running through any positive sequence $(\epsilon_{n})$ with $\lim_{n\to \infty}\epsilon_{n}=0$. This limit point is not zero operator. Moreover, $\|\LR_{A_{MM^\p},\Phe}1\|_{\infty}\asymp t_{\e}(MM^\p)$.
}
\pros
For $\om,\omp\in \Sigma_{MM^\p}$ with $\om_{0}\om_{1}=\omp_{0}\omp_{1}$, we have
\begin{align*}
\left|\frac{e^{\varPhi(\e,\om)}}{t_{\e}(MM^\p)}-\frac{e^{\varPhi(\e,\omp)}}{t_{\e}(MM^\p)}\right|&=\frac{e^{\varPhi(\e,\omp)}}{t_{\e}(MM^\p)}\left|e^{\varPhi(\e,\om)-\varPhi(\e,\omp)}-1\right|
\leq e^{c\theta^{2}}c d_{\theta}(\om,\omp)
\end{align*}
with $c=\sup_{\e>0}[[\Phe]]_{\theta}$. Note that $e^{\Phe}/t_{\e}(MM^\p)$ on $\Sigma_{MM^\p}$ is bounded by $1$.
Therefore Ascoli-Arzel\`a Theorem implies that 
there exists a subsequence $(\e^\p_{n})$ of $(\e_{n})$ such that $\chi_{\Sigma_{MM^\p}}e^{\varPhi(\e^\p_{n},\cd)}/t_{\e^\p_{n}}(MM^\p)$ converges to a function $\psi_{MM^\p}$ in $C(\SiAp)$.
For $f\in C(\SiAp)$ with $\|f\|_{\infty}\leq 1$, we obtain
\begin{align*}
&\left\|\frac{\LR_{A_{MM^\p},\Phe}f}{t_{\e}(MM^\p)}-\LR_{A_{MM^\p},0}(\psi_{MM^\p}f)\right\|_{\infty}\\
=&\sup_{\om\in \SiMp}\left|\sum_{i\,:\,A_{MM^\p}(i\om_{0})=1}\left(\frac{e^{\varPhi(\e,i\cd \om)}}{t_{\e}(MM^\p)}-\psi_{MM^\p}(i\cd \om)\right)f(i\cd \om)\right|\\
\leq& d\left\|\frac{\chi_{\Sigma_{MM^\p}}e^{\Phe}}{t_{\e}(MM^\p)}-\psi_{MM^\p}\right\|_{\infty}\to 0
\end{align*}
as $\e\to 0$ running through $(\e^\p_{n})$. Thus $\LR_{A_{MM^\p},\Phe}/t_{\e}(MM^\p)$ has a limit point $\LR_{A_{MM^\p},0}(\psi_{MM^\p}\cd)$ in $\LR(C(\SiAp))$.

We finally show $\|\LR_{A_{MM^\p},\Phe}1\|_{\infty}\asymp t_{\e}(MM^\p)$. We see
\begin{align*}
\|\LR_{A_{MM^\p},\Phe}1\|_{\infty}\leq \|\LR_{A_{MM^{\prime}},0}1\|_{\infty}t_{\e}(MM^\p)\leq d t_{\e}(MM^\p).
\end{align*}
On the other hand, take $\om(\e)\in \Si_{MM^\p}$ so that $t_{\e}(MM^\p)=e^{\varPhi(\e,\om(\e))}$. Thus
\ali
{
t_{\e}(MM^\p)=e^{\varPhi(\e,\om(\e))}\leq \LR_{A_{MM^\p},\Phe}1(\sigma_{A}\om(\e))\leq \|\LR_{A_{MM^\p},\Phe}1\|_{\infty}.
}
Hence we obtain the final assertion.
\proe
\subsubsection{Perturbation of thermodynamic spectral characteristics}\label{sec:perturb_TSC}
Assume that the conditions $(\Si.1)$-$(\Si.3)$ and $(\Ph.1)$-$(\Ph.2)$ are satisfied. Let $\bM$ be a subset $\bT$ so that the following condition holds:
\ite
{
\item[$(\Si.\bM)$] There exists $\e_{0}>0$ such that for any $0<\e\leq\e_{0}$, $\bT_{0}(A(\bM),\Phe)$ consists of only one element $M_{0}$ and this matrix do not depend on $\e$.
}
When no confusion can arise, we always assume $0<\e\leq \e_{0}$ under this condition. We note that if $\bM\cap \bT_{0}\neq \emptyset$ then $\bT_{k}(B(\bM),\ph)=\bM\cap \bT_{k}$ for $k=0,1$ are satisfied.
We also see that $M_{0}$ is a $1\times 1$ zero matrix $(0)$ iff $A(\bM)=(0)$ iff $\bM$ consists of one element $M$ and $A_{MM}=(0)$.
By using the notion (\ref{eq:TSC}), we denote
\alil
{
&(\lam(\bM,\e),h(\bM,\e,\cd),\nu(\bM,\e,\cd))\nonumber\\
=&
\case
{
(\lam_{A(\bM),\Phe},h_{A(\bM),\Phe},\nu_{A(\bM),\Phe}),&A(\bM)\neq O\\
(0,\chi_{\Sigma_{M}},\delta_{\om(M)}),&A(\bM)=O,
}\label{eq:tsc_L_AikPhe}\\
g(\bM,\e,\cd)=&h(\bM,\e,\cd)/\|h(\bM,\e,\cd)\|_{\infty},\label{eq:tsc_L_AikPhe_g}
}
where $\om(M)\in \Sigma_{M}$ is a fixed element.
For simplicity, if $\bM=\bT$ then we write these by
\begin{align} (\lam(\e),h(\e,\cd),\nu(\e,\cd)) \text{ and }g(\e,\cd),\label{eq:ge=}
\end{align}
and if $\bM$ consists of only one element $M$, then we denote by
\alil{
(\lam(M,\e),h(M,\e,\cd),\nu(M,\e,\cd)) \text{ and }g(M,\e,\cd).\label{eq:gMe=}
}
Let $\ti{h}(\bM,\e,\cd),\ti{g}(\bM,\e,\cd)\in C(\SiAp,\R)$ and $\ti{\nu}(\bM,\e,\cd)\in M(\SiAp)$ be
\alil
{
\ti{h}(\bM,\e,\cd)=\frac{\nu(\bM,\e,g(\e,\cd))}{g(\e,\cd)}h(\bM,\e,\cd),&\quad
\ti{g}(\bM,\e,\cd)=\ti{h}(\bM,\e,\cd)/\|\ti{h}(\bM,\e,\cd)\|_{\infty},\label{eq:tsc_L_AiktPhe_g}\\
\ti{\nu}(\bM,\e,f)=\frac{\nu(\bM,\e,g(\e,\cd)f)}{\nu(\bM,\e,g(\e,\cd))}&\text{ for }f\in C(\SiAp).\label{eq:tsc_L_AiktPhe}
}
In the setting, we see that the triplet $(\lam(\bM,\e),\ti{h}(\bM,\e,\cd),\ti{\nu}(\bM,\e,\cd))$ satisfies the property of the TSC of $\LR_{A(\bM),\tPhe}$, i.e.
\ali
{
&\LR_{A(\bM),\tPhe}\ti{h}(\bM,\e,\cd)=\lam(\bM,\e)\ti{h}(\bM,\e,\cd),\\
&\LR_{A(\bM),\tPhe}^{*}\ti{\nu}(\bM,\e,\cd)=\lam(\bM,\e)\ti{\nu}(\bM,\e,\cd),\quad
\ti{\nu}(\bM,\e,\ti{h}(\bM,\e,\cd))=1,
}
where $\tPhe$ is defined in (\ref{eq:tPhe=}). When $(\Ph.3)$ is imposed, Theorem \ref{th:aRPF} implies that the functions $g(\e,\cd)$ and $g(\bM,\e,\cd)$ are in $\Lambda_{c}$ with the constant $c=\theta\sup_{\e>0}[[\Phe]]_{\theta}/(1-\theta)$ and therefore so is for $\ti{g}(\bM,\e,\cd)$. Furthermore, $\supp h(\bM,\e,\cd)$ and $\supp \nu(\bM,\e,\cd)$ do not change by $\e$ from (\ref{eq:suppnu=}), and we see $\supp \ti{h}(\bM,\e,\cd)=\supp h(\bM,\e,\cd)$ and $\supp \ti{\nu}(\bM,\e,\cd)=\supp \nu(\bM,\e,\cd)$ from these definitions.

We begin with convergence of the Perron eigenvalue $\lam(\bM,\e)$ of $\LR_{A(\bM),\Phe}$.
\prop
{\label{prop:conv_pre_M}
Assume that the conditions $(\Si.1)$-$(\Si.2)$, $(\Ph.1)$-$(\Ph.2)$ and $(\Si.\bM)$ are satisfied. Then $\lam(\bM,\e)$ converges to the Perron eigenvalue $\lam(\bM)$ of $\LR_{B(\bM),\ph}\in \LR(C(\SiAp))$.
}
\pros
Note that $\lam(\bM,\cd)$ is bounded by $\|\LR_{A(\bM),\Phe}1\|_{\infty}$. There exists a sequence $(\epsilon_{n})$ with $\lim_{n\to \infty}\epsilon_{n}=0$ such that $\lam(\bM,\epsilon_{n})$ converges to a number $\lam_{0}\in \R$ and $\nu(\bM,\e_{n},\cd)$ converges weakly to a measure $\nu_{0}\in M(\SiAp)$ as $n\to \infty$. We notice that the equality $\LR_{B(\bM),\ph}^{*}\nu_{0}=\lam_{0}\nu_{0}$ is satisfied.

In the case when $B(\bM)=O$, $\LR_{B(\bM),\ph}^{d}$ becomes a zero operator and therefore $\lam_{0}=0=\lam(\bM)$ is fulfilled. 

In the case when $B(\bM)\neq O$, we see
\begin{align*}
\lam_{0}=\nu_{0}(\LR_{B(\bM),\varphi}^{n}1)^{1/n}\leq \|\LR_{B(\bM),\varphi}^{n}1\|_{\infty}^{1/n}\to \lam(\bM)
\end{align*}
as $n\to \infty$ by using Proposition 5.2 in \cite{T2009}.  The inequality $\lam_{0}\leq \lam(\bM)$ is obtained. To see the opposite inequality, we recall the equations $\log(\lam(\bM,\e))=P(\si_{A(\bM)},\Phe_{A(\bM)})$ and $\log(\lam(\bM))=P(\si_{B(\bM)},\ph_{B(\bM)})$. 
For $B(\bM)$-admissible word $w$, the inclusion ${}_{0}[w]^{B(\bM)}\subset {}_{0}[w]^{A(\bM)}$ holds and then $w$ is $A(\bM)$-admissible. Therefore we have that for any $n\geq 1$
\ali
{
\sum_{w\in W_{n}(B(\bM))}\sup_{\om\in {}_{0}[w]^{B(\bM)}}\exp (S_{n}\varphi(\e,\om))&= \sum_{w\in W_{n}(B(\bM))}\sup_{\om\in {}_{0}[w]^{B(\bM)}}\exp (S_{n}\varPhi(\e,\om))\\
&\leq \sum_{w\in W_{n}(A(\bM))}\sup_{\om\in {}_{0}[w]^{A(\bM)}}\exp (S_{n}\varPhi(\e,\om)).
}
Thus $P(\sigma_{B(\bM)}, \phe_{B(\bM)})\leq P(\sigma_{A(\bM)},\Phe_{A(\bM)})$ for any $\e>0$. We see the inequality $P(\sigma_{B(\bM)},\varphi_{B(\bM)})\leq \liminf_{\e\to \infty}P(\sigma_{A(\bM)},\Phe_{A(\bM)})$. Hence $\lam(\bM)\leq \lam_{0}$. Consequently, $\lam(\bM,\e)$ converges to $\lam(\bM)$.
\proe
\prop
{\label{prop:gM_asymp_g}
Assume that $(\Si.1)$-$(\Si.3)$, $(\Ph.1)$-$(\Ph.3)$ and $(\Si.\bM)$ are satisfied. Then there exist constants $c=c(B(\bM),\varphi)\geq 1$ and $\epsilon_{0}>0$ such that for any $i,j\in S$ with $(B(\bM))^{n}(ij)>0$ for some $n\geq 0$ and for any $\om\in {}_{0}[i]^{A}$, $\upsilon\in {}_{0}[j]^{A}$ and $0<\e<\epsilon_{0}$,
\item[(1)] $g(\bM,\e,\om)\leq c g(\bM,\e,\upsilon)$ and
\item[(2)] $\nu(\bM,\e,{}_{0}[j]^{A})\leq c \nu(\bM,\e,{}_{0}[i]^{A})$.
}
\pros
Without loss of generality, we may assume $d> n$.
\item[(1)] Assume $i\neq j$. Then there exists $w\in W_{n-1}(B(\bM))$ such that $i\cd w\cd j\in W_{n+1}(B(\bM))$. Choose any $\epsilon_{0}>0$ so that for any $0<\e<\epsilon_{0}$, $0<\lam(\bM,\e)\leq \lam(\bM)+1$ and $\|\phe\|_{\infty}\leq\|\varphi\|_{\infty}+1$. We have
\begin{align*}
g(\bM,\e,\upsilon)
&=\lam(\bM,\e)^{-n}\LR_{A(\bM),\Phe}^{n}g(\bM,\e,\upsilon)\\
&\geq \lam(\bM,\e)^{-n}\LR_{B(\bM),\phe}^{n}g(\bM,\e,\upsilon)\\
&\geq \lam(\bM,\e)^{-n}\exp(S_{n}\varphi(\e,\om_{0}\cd w\cd \upsilon))g(\bM,\e,\om_{0}\cd w\cd \upsilon)\\
&\geq \lam(\bM,\e)^{-n}\exp(-n\|\phe\|_{\infty}-c_{0}\theta)g(\bM,\e,\om)\geq c^{-1}g(\bM,\e,\om)
\end{align*}
for $0<\e<\epsilon_{0}$ by putting $c_{0}=\theta\sup_{\e>0}[[\Phe]]_{\theta}/(1-\theta)$ and $c=(\lam(\bM)+1)^{d}\exp(d(\|\varphi\|_{\infty}+1)+c_{0}\theta)$. 

On the other hand, if $i=j$ then we have $\om_{0}=\upsilon_{0}$ and $g(\bM,\e,\om)=\exp(c_{0}d_{\theta}(\om,\upsilon))g(\bM,\e,\upsilon)\leq \exp(c_{0}\theta)g(\bM,\e,\upsilon)$.
Thus $g(\bM,\e,\om)\leq c g(\bM,\e,\upsilon)$ is satisfied.
\smallskip
\item[(2)] The case $i=j$ is trivial. Assume $i\neq j$. By taking $w$ and $\e_{0}$ given in (1), we obtain
\begin{align*}
\nu(\bM,\e,{}_{0}[i]^{A})
&=\lam(\bM,\e)^{-n}\nu(\bM,\e,\LR_{A(\bM),\Phe}^{n}\chi_{{}_{0}[i]^{A}})\\
&\geq \lam(\bM,\e)^{-n}\nu(\bM,\e,\chi_{{}_{0}[j]^{A}}\LR_{B(\bM),\phe}^{n}\chi_{{}_{0}[i]^{A}})\\
&\geq \lam(\bM,\e)^{-n}\int_{{}_{0}[j]^{A}}\exp(S_{n}\varphi(\e,i\cd w\cd \om)) d \nu(\bM,\e,\om)\\
&\geq \lam(\bM,\e)^{-n}\exp(-n\|\phe\|_{\infty})\nu(\bM,\e,{}_{0}[j]^{A})\geq c^{-1}\nu(\bM,\e,{}_{0}[j]^{A})
\end{align*}
for any $0<\e<\epsilon_{0}$.
Hence the assertion holds.
\proe
By virtue of the above proposition, we have the asymptotic relation:
\alil
{
g(\bM,\e,\cd)\chi_{\SiM}\asymp g(\bM,\e,\om)\chi_{\SiM},\quad \nu(\bM,\e,\Si_{M})\asymp \nu(\bM,\e,{}_{0}[i]^{A})\label{eq:gasympg_nuasympnu}
}
for each $M\in \bT$, $\om\in \SiM$ and $i\in S_{M}$.
Now we obtain an important perturbation result for the TSC of $\LR_{A(\bM),\Phe}$:
\prop{
\label{prop:conv_nor_TSC_M}
Assume that the conditions $(\Si.1)$-$(\Si.3)$, $(\Ph.1)$-$(\Ph.3)$ and $(\Si.\bM)$ with $\bM\cap \bT_{0}\neq \emptyset$ are satisfied. Then
\item[(1)] For $M\in \bM\cap \bT_{0}$ with $g(\bM,\e,\cd)\neq 0$ on $\SiM$, $g(\bM,\e,\cd)\chi_{\SiM}/\nu(M,\e,g(\bM,\e,\cd))$ converges to $h(M,\cd)$ in $C(\SiAp)$.
\item[(2)] For $M\in \bM\cap \bT_{0}$ with $\nu(\bM,\e,\SiM)>0$, $\nu(\bM,\e,f\chi_{\SiM})/\nu(\bM,\e,\SiM)$ converges to $\nu(M,f)$ for each $f\in C(\SiAp)$.
}
\pros
\item[(1)] By the former of (\ref{eq:gasympg_nuasympnu}), we have $g(\bM,\e,\cd)\chi_{\SiM}/\nu(M,\e,g(\bM,\e,\cd))\asymp \Si_{M}$. We can take a sequence $(\epsilon_{n})$ and $g_{\infty}\in C(\SiAp)$ with $g_{\infty}\neq 0$ so that the function $g(\bM,\e_{n},\cd)\chi_{\SiM}/\nu(M,\e_{n},g(\bM,\e_{n},\cd))$ converges to $g_{\infty}$ in $C(\SiAp)$ as $n\to \infty$. Denoted by $\mathcal{I}$ the identity operator belonging to $\LR(C(\SiAp))\cap \LR(F_{\theta}(\SiAp))$. For $\om\in \SiM$,
\begin{align*}
&(\lam(\bM,\e)\mathcal{I}-\LR_{A_{MM},\Phe})\left(g(\bM,\e,\cd)\chi_{\SiM}\right)(\om)\\
=&\lam(\bM,\e)g(\bM,\e,\om)-\LR_{A_{MM},\Phe}g(\bM,\e,\om)\\
=&\LR_{A(\bM)-A_{MM},\Phe}g(\bM,\e,\om)\\
=&\sum_{M^\p\in \bM\,:\,M\neq M^{\prime}}\LR_{A_{M^\p M},\Phe}g(\bM,\e,\om)\\
\asymp& \sum_{M^\p\in \bM\,:\,M\neq M^\p}t_{\e}(M^\p M)\nu(M^\p,\e,g(\bM,\e,\cd)).
\end{align*}
by Proposition \ref{prop:conv_TSC_kkp}. Here the last expression has the relation
\begin{align*}
&(\lam(\bM,\e)-\lam(M,\e))\nu(M,\e,g(\bM,\e,\cd))\\
=&\nu(M,\e,\left(\LR_{A(\bM),\Phe}-\LR_{A_{MM},\Phe}\right)g(\bM,\e,\cd))\\
=&\nu(M,\e,\sum_{M^{\prime}\in \bM\,:\,M\neq M^{\prime}}\LR_{A_{M^{\prime}M},\Phe}g(\bM,\e,\om))\\
\asymp& \sum_{M^{\prime}\in \bM\,:\,M\neq M^{\prime}}t_{\e}(M^{\prime}M)\nu(M^{\prime},\e,g(\bM,\e,\cd)).
\end{align*}
Therefore we obtain the relation
\ali{
\left\|(\lam(\bM,\e)\mathcal{I}-\LR_{A_{MM},\Phe})\left(\frac{g(\bM,\e,\cd)\chi_{\SiM}}{\nu(M,\e,g(\bM,\e,\cd))}\right)\right\|_{\infty}\asymp \lam(\bM,\e)-\lam(M,\e).
}
Since $\bM\cap \bT_{0}$ is not empty and $M$ is in $\bT_{0}$, the eigenvalues $\lam(\bM,\e)$ and $\lam(M,\e)$ converge to both $\lam=\lam(B,\ph)$. This implies that the left hand side of the above relation tends to $0$. Thus we have the equation $\LR_{B_{MM},\varphi}g_{\infty}=\LR_{M,\varphi}g_{\infty}=\lam g_{\infty}$. This yields $g_{\infty}=\beta g(M,\cd)$ for some $\beta>0$ from $\lam$ is a simple eigenvalue of $\LR_{M,\ph}\in \LR(C(\SiAp))$. By
\begin{align*}
1=\frac{\nu(M,\e,g(\bM,\e,\cd))}{\nu(M,\e,g(\bM,\e,\cd))}&=\nu(M,\e,\frac{g(\bM,\e,\cd)\chi_{\SiM}}{\nu(M,\e,g(\bM,\e,\cd))})\\
&\to \nu(M,\beta g(M,\cd))=\beta \nu(M,g(M,\cd)),
\end{align*}
as $\e\to 0$ running through $(\epsilon_{n})$, $\beta=1/\nu(M,g(M,\cd))$ is satisfied. Consequently, we obtain $g_{\infty}=g(M,\cd)/\nu(M,g(M,\cd))=h(M,\cd)$. 
\item[(2)] We have
\begin{align}
&|(\lam(\bM,\e)\mathcal{I}^{*}-\LR_{A_{MM},\Phe}^{*})\nu(\bM,\e,f\chi_{\SiM})|\nonumber\\
=&|\nu(\bM,\e,(\lam(\bM,\e)\mathcal{I}-\LR_{A_{MM},\Phe})(f\chi_{\SiM}))|\nonumber\\
=&\left|\nu(\bM,\e,\LR_{A(\bM),\Phe}(f\chi_{\SiM})-\LR_{A_{MM},\Phe}(f\chi_{\SiM}))\right|\nonumber\\
=&\left|\nu(\bM,\e,\sum_{M^{\prime}\in \bM\,:\,M^{\prime}\neq M}\LR_{A_{MM^{\prime}},\Phe}(f\chi_{\SiM}))\right|\label{eq:conv_nor_TSC_1}\\
\leq &\nu(\bM,\e,\sum_{M^{\prime}\in \bM\,:\,M^{\prime}\neq M}\LR_{A_{MM^{\prime}},\Phe}\chi_{\SiM})\|f\|_{\infty}\nonumber\\
\leq&\sum_{M^{\prime}\in \bM\,:\,M^{\prime}\neq M}\nu(\bM,\e,\Sigma_{M^{\prime}})t_{\e}(MM^{\prime})\|f\|_{\infty}.\nonumber
\end{align}
On the other hand,
\alil
{
\lam(\bM,\e)-\lam(M,\e)=
&\frac{\nu(\bM,\e,\left(\LR_{A(\bM),\Phe}-\LR_{A_{MM},\Phe}\right)h(M,\e,\cd))}{\nu(\bM,\e,h(M,\e,\cd))}\nonumber\\
\asymp&\frac{\sum_{M^{\prime}\in \bM\,:\,M^{\prime}\neq M}\nu(\bM,\e,h(M^{\prime},\e,\cd))t_{\e}(MM^{\prime})}{\nu(\bM,\e,h(M,\e,\cd))}\label{eq:conv_nor_TSC_2}\\
\asymp& \frac{\sum_{M^{\prime}\in \bM\,:\,M^{\prime}\neq M}\nu(\bM,\e,\Sigma_{M^{\prime}})t_{\e}(MM^{\prime})}{\nu(\bM,\e,\SiM)}.\nonumber
}
from $h(M^\p,\e,\cd)\asymp \chi_{\SiMp}$. 
Let $m(M,\e,f)=\nu(\bM,\e,f\chi_{\SiM})/\nu(\bM,\e,\SiM)$ for $f\in C(\SiAp)$. There exist a sequence $(\ep_{n})$ and $m_{\infty}\in M(\SiAp)$ such that $m(M,\ep_{n},\cd)\to m_{\infty}$ as $n\to \infty$ weakly. 
Thus the inequalities (\ref{eq:conv_nor_TSC_1}) and (\ref{eq:conv_nor_TSC_2}) yield
$\LR_{M,\varphi}^{*}m_{\infty}=\lam m_{\infty}$.
This implies $m_{\infty}=\gamma \nu(M,\cd)$ for some constant $\gamma>0$ by the simplicity of $\lam$ of $\LR_{M,\ph}^{*}$. Hence $m_{\infty}=\nu(M,\cd)$ is yielded by $m_{\infty}(\SiAp)=\nu(M,\SiAp)=1$.
\proe
Finally, we consider a special case when $T_{0}(B(\bM),\ph)$ consists of only one element.
\prop{
\label{prop:tsc_LAMPhe_T0=1}
Assume that $(\Si.1)$-$(\Si.3)$, and $(\Ph.1)$-$(\Ph.3)$ are satisfied. Assume also that $T_{0}(B(\bM),\ph)$ consists of only one element $M$. Then we have the following:
\item[(1)] The condition $(\Si.\bM)$ is valid.
\item[(2)] The TSC $(\lam(\bM,\e),h(\bM,\e,\cd),\nu(\bM,\e,\cd))$ of $\LR_{A(\bM),\Phe}$ converges to the TSC $(\lam(\bM),h(\bM,\cd),\nu(\bM,\cd))$ of  $\LR_{B(\bM),\varphi}$, and $g(\bM,\e,\cd)$ converges to $g(\bM,\cd)$ in $C(\SiAp)$.
\item[(3)] $h(\bM,\cd)\chi_{\SiM}=h(M,\cd)/\nu(\bM,\SiM)$ and $\nu(\bM,f\chi_{\SiM})=\nu(\bM,\SiM)\nu(M,f)$ for $f\in C(\SiAp)$.
}
\pros
\item[(1)] We see that there exists $M_{0}\in \bT(A(\bM))$ such that $S_{M}\subset S_{M_{0}}$. Since $\lim_{\e\to 0}\lam(M_{0},\e)=\lam$ and $\lim_{\e\to 0}\lam(M_{1},\e)<\lam$ for any $M_{1}\in\bT(A(\bM),\Phe)\setminus \{M_{0}\}$, $\bT_{0}(A(\bM),\Phe)=\{M_{0}\}$ is satisfied for a small $\e>0$.
\item[(2)] Take a limit point $g\in C(\SiAp)$ of $g(\bM,\e,\cd)$ and a limit point $\nu\in M(\SiAp)$ of $\nu(\bM,\e,\cd)$. Since $\LR_{A(\bM),\Phe}$ converges to $\LR_{B(\bM),\varphi}$ in $\LR(C(\SiAp))$ and since $\lam(\bM,\e)$ converges to $\lam(\bM)$ by Proposition \ref{prop:conv_pre_M}, we have the equations $\LR_{B(\bM),\ph}g=\lam(\bM) g$ and $\LR_{B(\bM),\ph}^{*}\nu=\lam(\bM) \nu$. By virtue of Theorem \ref{th:aRPF}, we obtain $g=g(M,\cd)$ and $\nu=\nu(M,\cd)$. In particular, $\nu(g)$ is positive and $h(\bM,\e,\cd)$ converges to $g/\nu(g)=h(M,\cd)$.
\item[(3)] Since $\supp \nu(\bM,\cd) \cap \supp h(\bM,\cd)=\Si_{M}^{+}$, we see
\ali
{
1=\nu(\bM,h(\bM,\cd))=\nu(\bM,h(\bM,\cd)\chi_{\SiM}).
} 
By Proposition 5.3(3) in \cite{T2009}, the equation $h(\bM,\cd)\chi_{\SiM}=c h(M,\cd)$ holds for some constant $c>0$. Similarity, Proposition 5.4(3) in \cite{T2009} yields $\nu(\bM,f\chi_{\SiM})=c^{\p}\nu(M,f)$ for some $c^\p>0$. When $f=1$, we have $c^\p=\nu(\bM,\SiM)$. Furthermore, 
\ali
{
1=\nu(\bM,h(\bM,\cd)\chi_{\SiM})=c^\p c\nu(M,h(M,\cd))=c^\p c.
}
Hence the assertion is valid by $c=1/\nu(\bM,\SiM)$.
\proe

\section{Matrix representations of linear operators}\label{sec:mat_bdd}
The eigenvalue and the eigenvector of bounded linear operators will be able to be reduced to the eigenvalue and the eigenvector of the guided finite matrix under suitable conditions. This idea might be useful for the analysis of the convergence of the eigenfunction and the eigenvector of transfer Ruelle operators. In this section, the abstract formulation is described in the first half, and applications to generalized Ruelle operators are mentioned in the latter half. Note that the way of making our matrix in this section resembles the matrix which appears to block numerical ranges of block operator matrices (e.g. \cite{Tretter}).
\subsection{An abstract formulation}\label{sec:mat_abst}
Let $T$ be a finite set and $\mathcal{X}$ a Banach algebra over $\K=\C$ or $\R$. Take $\LR_{k}\in \LR(\mathcal{X})$ for each $k\in T$. Denoted by $1_{\mathcal{X}}$ the unit element of $\mathcal{X}$. 
We assume the following:
\begin{itemize}
\item[(MR.1)] Pairs $(\eta,\nu)\in \K\times \mathcal{X}^{*}$ with $\LR^{*}\nu=\eta \nu$ and $(\eta_{k},h_{k})\in \K\times \mathcal{X}$ with $\LR_{k}h_{k}=\eta_{k}h_{k}$ for $k\in T$ satisfy $\nu(h_{k})\neq 0$ for some $k\in T$.
\item[(MR.2)] There exists  a subset $\{f_{kk^\p}\in \mathcal{X}\,:\,k,k^{\prime}\in T \text{ s.t. }\nu(h_{k^{\prime}})\neq 0\}$ of $\mathcal{X}$ such that $\nu((\LR-\LR_{k})h_{k})=\sum_{k^{\prime}\,:\,\nu(h_{k^{\prime}})\neq 0}\nu(f_{kk^\p}(\LR-\LR_{k})h_{k})$ for each $k\in T$.
\end{itemize}
Note that if $1_{\mathcal{X}}=\sum_{k^{\prime}\,:\,\nu(h_{k^{\prime}})\neq 0}f_{kk^\p}$ for $k\in T$ then (MR.2) is valid.
For $k\in T$, we have
\begin{align*}
\eta\nu(h_{k})=\nu(\LR h_{k})&=\nu((\LR-\LR_{k})h_{k})+\nu(\LR_{k}h_{k})=\sum_{k^{\prime}\in T}\mV(kk^\p)\bV(k^\p)
\end{align*}
by putting
\alil
{
\bV(k)&=\nu(h_{k})
,\quad \mV(kk^\p)=\delta_{kk^\p}\eta_{k}+\begin{cases}
\frac{\nu(f_{kk^\p}(\LR-\LR_{k})h_{k})}{\nu(h_{k^{\prime}})},& \nu(h_{k^{\prime}})\neq 0\\
0,& \nu(h_{k^{\prime}})= 0\label{eq:MR1=}
\end{cases}
}
for $k,k^\p\in T$, where $\delta_{kk^\p}=1$ if $k=k^\p$ and $\delta_{kk^\p}=0$ otherwise. Therefore $\eta$ is a eigenvalue of the matrix $\mV=(\mV(kk^\p))$ and the corresponding right eigenvector is $\bV=(\bV(k))$.

Next we give a transformed operator $\tilde{\LR}\in \LR(\XR)$ defined by
\begin{align}
\tilde{\LR}f=\xi^{-1}\LR(\xi f),\label{eq:tL=}
\end{align}
where $\xi\in \mathcal{X}$ satisfies that $\xi^{-1}\in \mathcal{X}$ exists. In addition to the conditions (MR.1) and (MR.2), we introduce the following:
\begin{itemize}
\item[(MR.3)] $\xi$ and $f_{kk^\p}$ are commutation for each $k,k^{\prime}$, i.e. $\xi f_{kk^\p}=f_{kk^\p}\xi$.
\item[(MR.4)] $\nu(\xi)\neq 0$. Further,  for each $k\in T$, there exists $\nu_{k}\in \XR^{*}$ with $\LR_{k}^{*}\nu_{k}=\eta_{k}\nu_{k}$ such that $\nu_{k}(h_{k})\neq 0$ and $\nu_{k}(\xi)\neq 0$.
\end{itemize}
We define $\ti{\nu}, \ti{\nu}_{k} \in \XR^{*}$ and $\ti{h}_{k}\in \XR$ by
\alil
{
\tilde{\nu}(f)=\frac{\nu(\xi f)}{\nu(\xi)},\quad
\tilde{\nu}_{k}(f)&=\frac{\nu_{k}(\xi f)}{\nu_{k}(\xi)},\quad \tilde{h}_{k}=\frac{\nu_{k}(\xi)\xi^{-1}h_{k}}{\nu_{k}(h_{k})}.\label{eq:tnuk=}
}
\prop{\label{prop:mr_prop_nor}
Assume (MR.1)-(MR.4). Then 
\ite{
\item The pairs $(\lam,\ti{\nu})$ and $(\lam_{k},\ti{h}_{k})$ satisfy the conditions (MR.1), (MR.2), and $\ti{\nu}_{k}(\ti{h}_{k})=1$ for $k\in T$.
\item We let
\alil{
\tmV(kk^\p)=\mV(kk^\p)\frac{\nu_{k}(\xi)\nu_{k^{\prime}}(h_{k^{\prime}})}{\nu_{k^{\prime}}(\xi)\nu_{k}(h_{k})},\quad
\tbV(k)=\bV(k)\frac{\nu_{k}(\xi)}{\nu(\xi)\nu_{k}(h_{k})}\label{eq:tmVkkp=}
}
for $k,k^\p\in T$. Then the pair $(\tmV, \tbV)$ has the equations (\ref{eq:MR1=}) by replacing $\LR=\tLR,\nu=\ti{\nu},\LR_{k}=\tLR_{k}$ and $h_{k}=\ti{h}_{k}$.
}
}
\pros
Proofs directly follow.
\proe
On the other hand, we consider a symmetric situation of the above:
\begin{itemize}
\item[(MR.5)] Pairs $(\eta,h)\in \K\times \mathcal{X}$ with $\LR h=\eta h$ and $(\eta_{k},\nu_{k})\in \K\times \mathcal{X}^{*}$ with $\LR_{k}^{*}\nu_{k}=\eta_{k}\nu_{k}$ for $k\in T$ satisfy $\nu_{k}(h)\neq 0$ for some $k$.
\item[(MR.6)] There exists a subset $\{g_{k^{\prime}k}\in \XR\,:\,k,k^{\prime}\in T \text{ s.t. }\nu_{k^{\prime}}(h)\neq 0\}$ of $\XR$ such that $\nu_{k}((\LR-\LR_{k})h)=\sum_{k^{\prime}\,:\,\nu_{k^{\prime}}(h)\neq 0}\nu_{k}((\LR-\LR_{k})(g_{k^{\prime}k}h))$ for $k\in T$.
\end{itemize}
Under (MR.5) and (MR.6), we have
\begin{align*}
\eta\nu_{k}(h)=\nu_{k}(\LR h)&=\nu_{k}((\LR-\LR_{k})h)+\nu_{k}(\LR_{k}h)
=\sum_{k^{\prime}\in T}\mG(k^{\prime}k)\bG(k^{\prime})
\end{align*}
with
\alil
{
\bG(k)=\nu_{k}(h),\quad
\mG(k^{\prime}k)=\eta_{k}\delta_{kk^\p}+\begin{cases}
\frac{\nu_{k}((\LR-\LR_{k})(g_{k^{\prime}k}h))}{\nu_{k^{\prime}}(h)},&\nu_{k^{\prime}}(h)\neq 0\\
0,&\nu_{k^{\prime}}(h)=0\\
\end{cases}\label{eq:MR2=}
}
Then $\eta$ becomes the eigenvalue of the matrix $\mG=(\mG(kk^\p))$ and $\bG=(\bG(k))$ is the corresponding left eigenvector.
Take the operator (\ref{eq:tL=}). We also assume the following:
\begin{itemize}
\item[(MR.7)] $\xi$ and $g_{kk^\p}$ are commutation for each $k,k^\p\in T$.
\item[(MR.8)] An eigenvector $\nu\in \XR^{*}$ with $\LR\nu=\eta\nu$ has the condition $\nu(h)\neq 0$.
\end{itemize}
Put 
\alil
{
\tilde{h}&=\xi^{-1}h\nu(h)^{-1}\nu(\xi).\label{eq:th=}
}
\prop{\label{prop:mr_prop_nor2}
Assume (MR.4)-(MR.8) and $\ti{\nu}_{k}$ is defined in (\ref{eq:tnuk=}). Then
\ite
{
\item[(1)] The pairs $(\lam,\ti{h})$ and $(\lam_{k},\ti{\nu}_{k})$ satisfy the conditions (MR.5), (MR.6), and $\ti{\nu}(\ti{h})=1$.
\item[(2)] We set
\alil{
\tmG(kk^\p)=\mG(kk^\p)\frac{\nu_{k}(\xi)}{\nu_{k^\p}(\xi)},\quad
\tbG(k)=\bG(k)\frac{\nu(\xi)}{\nu_{k}(\xi)\nu(h)}\label{eq:tmGkkp=}
}
}
for $k,k^\p\in T$. Then the pair $(\tmG, \tbG)$ gives the equation (\ref{eq:MR2=}) by replacing $\LR=\tLR,h=\ti{h},\LR_{k}=\tLR_{k}$ and $\nu_{k}=\ti{\nu}_{k}$.
}
\pros
The proofs immediately follow.
\proe
We sometimes call the four set
\alil
{
(\mV,\mG,\bV,\bG)
}
a {\it spectral matrix representation} (\MR for a short) of $((\LR,$ $\eta,$ $h,$ $\nu))$ using $((\LR_{k},$ $\eta_{k},$ $h_{k}, $ $\nu_{k})_{k}$, $(f_{kk^\p})_{kk^\p}$, $(g_{kk^\p})_{kk^\p})$.
\subsection{Matrix representation of Ruelle operators (I)}\label{sec:MR_LAiP}
Assume the conditions $(\Si.1)$-$(\Si.3)$,  $(\Ph.1)$-$(\Ph.3)$ and $(\Si.\bM)$. In this section, we study a spectral matrix representation of the operator $\LR_{A(\bM),\Phe}$ by using $\LR_{A_{MM},\Phe}$ for $M\in \bM$. We use the notation in Section \ref{sec:perturb_TSC}.
Let $f_{MM^\p}=g_{M^\p M}=\chi_{\SiMp}$ for each $M,M^\p\in \bT(B(\bM))=\bM$. We take the \MR
\alil{
(\mV_{\e}(\bM,\cdot),\mG_{\e}(\bM,\cdot), \bV_{\e}(\bM,\cdot), \bG_{\e}(\bM,\cdot))\label{eq:MR1_mVeM}
}
of $(\LR_{A(\bM),\Phe},$ $\lam(\bM,\e),$ $\nu(\bM,\e,\cd),$ $g(\bM,\e,\cd))$ using the sets $((\LR_{A_{MM},\Phe},$ $\lam(M,\e),$ $h(M,\e,\cd),$ $\nu(M,\e,\cd))_{M}$, $(f_{MM^\p})$, $(g_{MM^\p}))$.
Recall the function $\ti{g}(\bM,\e,\cd)$ and the TSC ($\lam(\bM,\e),\ti{h}(\bM,\e,\cd),\ti{\nu}(\bM,\e,\cd))$ of the operator $\LR_{A(\bM),\tPhe}$ defined in (\ref{eq:tsc_L_AiktPhe_g}) and (\ref{eq:tsc_L_AiktPhe}).
We pay attention to the equation
\ali
{
\LR_{A(\bM),\tPhe}f=g(\e,\cd)^{-1}\LR_{A(\bM),\Phe}(g(\e,\cd)f).
}
By putting $\xi=g(\e,\cd)$, we can take
\alil{
(\tmV_{\e}(\bM,\cdot),\tmG_{\e}(\bM,\cdot), \tbV_{\e}(\bM,\cdot), \tbG_{\e}(\bM,\cdot))\label{eq:MR1_tmVeM}
}
by the \MR given by (\ref{eq:tmVkkp=}) and (\ref{eq:tmGkkp=}) formed with $(\LR_{A(\bM),\tPhe},$ $\lam(\bM,\e),$ $\ti{g}(\bM,\e,\cd),$ $\ti{\nu}(\bM,\e,\cd))$, $(\LR_{A_{MM},\tPhe},$ $\lam(M,\e),$ $\ti{h}(M,\e,\cd),$ $\ti{\nu}(M,\e,\cd))_{M}$.
For a simple, if $\bM=\bT$ then we may omit $\bM$ from notation of (\ref{eq:MR1_mVeM}) and (\ref{eq:MR1_tmVeM}), i.e. we may write those as
\ali{
(\mV_{\e},\mG_{\e}, \bV_{\e}, \bG_{\e}),\quad (\tmV_{\e},\tmG_{\e}, \tbV_{\e}, \tbG_{\e}).
}
In those setting, we have the forms
\ali
{
&\bV_{\e}(\bM,M)=\nu(\bM,\e,h(M,\e,\cd))\\
&\mV_{\e}(\bM,MM^\p)\\
=&\begin{cases}
\lam(M,\e),&M=M^\p\\
\displaystyle \frac{\nu(\bM,\e,\LR_{A_{MM^\p},\Phe}h(M,\e,\cd))}{\nu(\bM,\e,h(M^\p,\e,\cd))},&M\neq M^\p \text{ and }\bV_{\e}(\bM,M^\p)>0\\
0,& \text{otherwise}
\end{cases}\\
&\bG_{\e}(\bM,M)=\nu(M,\e,g(\bM, \e,\cd)),\\
&\mG_{\e}(\bM,MM^\p)\\
=&\begin{cases}
\lam(M,\e),&M=M^\p\\
\displaystyle \frac{\nu(M^\p,\e,\LR_{A_{MM^\p},\Phe}g(\bM,\e,\cd))}{\nu(M,\e,g(\bM,\e,\cd))},&M\neq M^\p \text{ and }\bG_{\e}(\bM,M)>0\\
0,& \text{otherwise}.
\end{cases}
}
For $M,M^\p\in \bM$, we define
\alil
{
\tte(\bM,MM^\p)&=
\case
{
\frac{\bG_{\e}(\bM,M)}{\bG_{\e}(\bM,M^\p)}\te(MM^\p),&\bG_{\e}(\bM,M^\p)>0\\
0,&\bG_{\e}(\bM,M^\p)=0
},\\
\ttte(\bM,MM^\p)&=
\case
{
\frac{\bV_{\e}(\bM,M^\p)}{\bV_{\e}(\bM,M)}\te(MM^\p),&\bV_{\e}(\bM,M)>0\\
0,&\bV_{\e}(\bM,M)=0
}
}
and for $W=W_{1}W_{2}\cdots W_{n}\in \bM^{n}$
\alil
{
\te(W)=&\prod_{k=1}^{n-1}\te(W_{k}W_{k+1}),\\ \tte(\bM,W)=\prod_{k=1}^{n-1}\tte(\bM,W_{k}W_{k+1}),&\qquad \ttte(\bM,W)=\prod_{k=1}^{n-1}\ttte(\bM,W_{k}W_{k+1}).
}
Let
\ali
{
\VTe(\bM,M)&=\sum_{M^\p\in \bM\cap \bT_{0}}\sum_{W\in \SP(M, M^\p:\bM\cap \bT_{1})}\te(W)\bV_{\e}(\bM,M^\p)\\
\GTe(\bM,M)&=\sum_{M^\p\in \bM\cap \bT_{0}}\sum_{W\in \SP(M^\p,M: \bM\cap \bT_{1})}\bG_{\e}(\bM,M)\te(W).
}
\prop{\label{prop:prop_mr_LAMPhe3}
Assume that the conditions $(\Si.1)$-$(\Si.3)$, $(\Ph.1)$-$(\Ph.3)$ and $(\Si.\bM)$ with $\bM\cap \bT_{0}\neq \emptyset$ are satisfied. Then we have the following:
\ite
{
\item For all $M,M^\p\in \bM$ with $\bV_{\e}(\bM,M^\p)>0$, $\mV_{\e}(\bM,MM^\p)\asymp t_{\e}(MM^\p)$.
\item For all $M,M^\p\in \bM$ with $\bG_{\e}(\bM,M)>0$, $\mG_{\e}(\bM,MM^\p)\asymp t_{\e}(MM^\p)$.
\item For each $M\in \bT_{1}(B(\bM),\ph)$,
$\bV_{\e}(\bM,M)\asymp \VTe(\bM,M)$.
\item For each $M\in \bT_{1}(B(\bM),\ph)$,
$\bG_{\e}(\bM,M)\asymp \GTe(\bM,M)$.
\item For $M,M^\p\in \bM$, $\tte(\bM,M^\p M)=O(1)$. In particular, if $M\in \bT_{0}(B(\bM),\ph)$ and $M^\p\neq M$, then $\tte(\bM,M^\p M)\to 0$.
\item For $M,M^\p\in \bM$, $\tte(\bM,M M^\p)=O(1)$. In particular, if $M\in \bT_{0}(B(\bM),\ph)$ and $M^\p\in \bM$, then $\ttte(\bM,MM^\p)\to 0$.
}
}
\pros
\item[(1)] We have
\begin{align*}
\mV_{\e}(\bM,MM^\p)&=\frac{\nu(\bM,\e,\LR_{A_{MM^\p},\Phe}h(M,\e,\cd))}{\nu(\bM,\e,h(M^\p,\e,\cd))}\\
&\asymp\frac{\nu(\bM,\e,\LR_{A_{MM^\p},\Phe}\chi_{\SiM})}{\nu(\bM,\e,\chi_{\SiMp})}\asymp t_{\e}(MM^{\p})
\end{align*}
by using Proposition \ref{prop:conv_TSC_kkp} and the fact $h(M,\e,\cd)\asymp \chi_{\SiM}$ on $\SiM$.
\item[(2)] Similarity, we obtain
\begin{align*}
\mG_{\e}(\bM,MM^{\prime})&=\frac{\nu(M^\p,\e,\LR_{A_{MM^\p},\Phe}g(\bM,\e,\cd))}{\nu(M,\e,g(\bM,\e,\cd))}\\
&\asymp \nu(M^\p,\e,\LR_{A_{MM^\p},\Phe}\chi_{\SiM})\asymp t_{\e}(MM^\p)
\end{align*}
from the relation (\ref{eq:gasympg_nuasympnu}).
\item[(3)] By virtue of Proposition \ref{prop:convspd_evec_T1}, the relation
\ali{
\bV_{\e}(\bM,M)&\asymp \sum_{M^\p\in \bT_{0}(B(\bM),\ph)}\sum_{W\in \SP(M,M^\p:\bT_{1}(B(\bM),\ph))}\mV_{\e}(\bM,W)\bV_{\e}(\bM,M^\p)
}
is satisfied, where we define
$\mV_{\e}(\bM,W)=\prod_{i=1}^{k-1}\mV_{\e}(\bM,W_{i}W_{i+1})$ for $W=W_{1}W_{2}\cdots W_{k}\in \bT^{k}$.
By using (1), we see the assertion.
\item[(4)] By a similar argument above (3), the assertion follows from Proposition \ref{prop:convspd_evec_T1} and (2).
\item[(5)] We note the relation
\ali
{
\sum_{M^\p\in \bM\,:\,M^\p\neq M}\tte(\bM,M^\p M)&\asymp \sum_{M^\p\in \bM\,:\,M^\p\neq M,\bG_{\e}(\bM,M)>0}\frac{\mG_{\e}(\bM,M^\p M)\bG_{\e}(\bM,M^\p)}{\bG_{\e}(\bM,M)}\\
&=\lam(\bM,\e)-\lam(M,\e).
}
Therefore $\tte(\bM,M^\p M)$ is bounded. When $M\in \bT_{0}(B(\bM),\ph)$ with $M^\p\neq M$, we have $\tte(\bM,M^\p M)\to 0$ by $\lam(\bM,\e)-\lam(M,\e)\to 0$.
\item[(6)] By a similar argument above (5), we obtain the assertion.
\proe
Finally, we consider the speed of convergence of $g(\bM, \e,\cd)$ on $\SiM$ and of $\nu(\bM,\e,\SiM)$ for each $M\in \bM\cap \bT_{1}$. These results are useful to show main theorems.
\prop{\label{prop:speed_tscM}
Assume the conditions $(\Si.1)$-$(\Si.3)$, $(\Ph.1)$-$(\Ph.3)$ and $(\Si.\bM)$ with $\bM\cap \bT_{0}\neq \emptyset$. Assume also that for each $M,M^\p\in \bM$ with $A_{MM^\p}\neq O$, 
 $\LR_{A_{MM^\p},\Phe}/\te(MM^\p)$ converges to an operator $\LT(MM^\p)$ in $\LR(C(\SiAp))$.
\item[(1)] If
$\bG_{\e}(\bM,M^\p)\te(W)/\GTe(\bM,M)$ converges to a number $c_{1}(W)$ in $[0,1]$ for each $M\in \bM\cap\bT_{1}$, $M^\p\in \bM\cap\bT_{0}$, and $W\in \SP(M^\p,M : \bM\cap\bT_{1})$ as $\e\to 0$, then for each $M\in \bM\cap\bT_{1}$
\ali
{
\frac{g(\bM,\e,\cd)}{\GTe(\bM,M)}\to \sum_{M^\p\in \bM\cap\bT_{0}}\sum_{W\in \SP(M^\p,M : \bM\cap\bT_{1})}c_{1}(W)\RLR(W)h(M^\p,\cd)
}
on $\SiM$ as $\e\to 0$, where $\RLR(W)$ is defined by
\ali
{
\RLR(MM^\p)&=(\lam\IR-\LR_{M^\p,\ph})^{-1}\LT(MM^\p)\\
\RLR(W)&=\RLR(W_{n-1}W_{n})\RLR(W_{n-2}W_{n-1})\cdots \RLR(W_{1}W_{2})
}
for each $M\in \bT$, $M^\p\in \bT_{1}$ and $W=W_{1}W_{2}\cdots W_{n}\in \bT\times \bT_{1}^{n-1}$.
\item[(2)] If $\te(W)\bV_{\e}(\bM,M^\p)/\VTe(\bM,M)$ converges to a number $c_{2}(W)$ in $[0,1]$ for each $M\in \bM\cap\bT_{1}\,\ M^\p\in \bM\cap\bT_{0}$, and $W\in \SP(M,M^\p : \bM\cap\bT_{1})$ as $\e\to 0$, then for each $M\in \bM\cap\bT_{1}$ and $f\in C(\SiAp)$,
\ali
{
\frac{\nu(\bM,\e,\chi_{\SiM}f)}{\VTe(\bM,M)}\to \sum_{M^\p\in \bM\cap\bT_{0}}\sum_{W\in \SP(M,M^\p : \bM\cap \bT_{1})}c_{2}(W)\nu(M^\p,\LRR(W)f)
}
as $\e\to 0$, where $\LRR(W)$ is defined by
\ali
{
\LRR(MM^\p)&=\LT(MM^\p)(\lam\IR-\LR_{M,\ph})^{-1}\\
\LRR(W)&=\LRR(W_{n-1}W_{n})\LRR(W_{n-2}W_{n-1})\cdots \LRR(W_{1}W_{2})
}
for each $M\in \bT_{1}$, $M^\p\in \bT$ and $W=W_{1}W_{2}\cdots W_{n}\in \bT_{1}^{n-1}\times \bT$.
}
\pros
Put $\bM_{0}=\bM\cap \bT_{0}$, $\bM_{1}=\bM\cap \bT_{1}$ and $m_{1}=\sharp \bM_{1}$.
\item[(1)] 
For $\e>0$, $M\in \bT$, $M^\p\in \bM_{1}$ and $W=W_{1}W_{2}\cdots W_{n}\in \bM\times (\bM_{1})^{n-1}$, we give the notation
\ali
{
\RLR_{\e}(MM^\p)&=(\lam(\e)\IR-\LR_{A_{M^\p M^\p},\Phe})^{-1}\LR_{A_{MM^\p},\Phe}\\
\RLR_{\e}(W)&=\RLR_{\e}(W_{n-1}W_{n})\RLR_{\e}(W_{n-2}W_{n-1})\cdots\RLR_{\e}(W_{1}W_{2}).
}
For $M\in \bT_{1}$, let $\mM_{M}=(\mM_{M}(M^\p M^\pp))$ be zero-one matrix indexed by $(\bM_{1}\setminus \{M\})\times(\bM_{1}\setminus \{M\})$ defined as $\mM_{M}(M^\p M^\p)=0$ and $\mM_{M}(M^\p M^\pp)=1$ for $M^\p\neq M^\pp$.
Take $M\in \bM_{1}$ and $\omega\in \SiM$. We consider the decomposition
\ali{
\lam(\bM,\e)g(\bM,\e,\om)=&\LR_{A(\bM),\Phe}g(\bM,\e,\om)\\
=&\sum_{M^\p\in \bM_{0}}\LR_{A_{M^\p M},\Phe}g(\bM,\e,\om)\\
&+\sum_{M^\p\in \bM_{1}\setminus \{M\}}\LR_{A_{M^\p M},\Phe}g(\bM,\e,\om)+\LR_{A_{MM},\Phe}g(\bM,\e,\om).
}
Therefore, we have the equation
\ali{
g(\bM,\e,\om)=&\sum_{M^\p\in \bM_{0}}\RLR_{\e}(M^\p M)g(\bM,\e,\om)+\sum_{M^\p\in \bM_{1}\setminus \{M\}}\RLR_{\e}(M^\p M)g(\bM,\e,\om).
}
By using the above repeatedly, we obtain
\alil{
g(\bM,\e,\om)
=&\sum_{M^\p\in \bM_{0}}\sum_{i=0}^{m_{1}-1}\sum_{W\in W_{i}(\mM_{M})}\RLR_{\e}(M^\p W M)g(\bM,\e,\om)\nonumber\\
&+\sum_{i=0}^{m_{1}-1}\sum_{W\in W_{i}(\mM_{M})}\RLR_{\e}(M W M)g(\bM,\e,\om)\label{eq:speed_tsc_2}\\
&+\sum_{W\in W_{m_{1}}(\mM_{M})}\RLR_{\e}(WM)g(\bM,\e,\om).\nonumber
}
To calculate $\lim_{\e\to 0}g(\bM,\e,\om)/\GTe(\bM,M)$, we consider the following.
For $0\leq i\leq m_{1}-1$ and $W\in W_{i}(\mM_{M})$, if $M^\p W M$ is not simple path, there exist $W_{1},W_{2}\in W_{*}(\mM_{M})$ such that $W_{2}$ is cycle, $M^{\p} W_{1} M$ is in $\SP(M^\p,M:\bM_{1})$, and the relation $\RLR_{\e}(M^\p W M)g(\bM,\e,\cdot)\asymp \te(W_{2})\te(M^\p W_{1} M)\bG_{\e}(\bM,M)$ is satisfied. Thus we see
\ali{
\frac{\RLR_{\e}(M^\p W_{1} M)g(\bM,\e,\om)}{\GTe(\bM,M)}\asymp \frac{\te(W_{2})\te(M^\p W_{1} M)\bG_{\e}(\bM,M)}{\GTe(\bM,M)}\leq \te(W_{2})\to 0.
}
If $M^\p W M$ is a simple path, then we have
\ali{
\frac{\RLR_{\e}(M^\p W M)g(\bM,\e,\om)}{\GTe(\bM,M)}&=\frac{\te(M^\p W M)\bG_{\e}(\bM,M^\p)}{\GTe(\bM,M)}\frac{\RLR_{\e}(M^\p W M)g(\bM,\e,\om)}{\te(M^\p W M)\bG_{\e}(\bM,M^\p)}\\
&\to c_{1}(M^\p W M)\RLR(M^\p W M)h(M^\p,\om)
}
and
\ali{
\RLR_{\e}(M W M)\frac{g(\bM,\e,\om)}{\GTe(\bM,M)}&\asymp \te(M W M)\frac{\bG_{\e}(\bM,M)}{\GTe(\bM,M)}\asymp \te(M W M)\to 0
}
by virtue of Proposition \ref{prop:prop_mr_LAMPhe3}(4).
For $W\in W_{m_{1}}(\mM_{M})$, $W M$ contains of a simple path and therefore there exist $W_{1},W_{2}\in W_{*}(\mM_{M})$ such that $W_{2}$ is cycle, $W_{1} M\in \SP(M,M^\p:\bM_{1})$ and $\te(W M)=\te(W_{2})\te(W_{1} M)$. Thus we see
\ali{
\frac{\RLR_{\e}(W M)g(\bM,\e,\om)}{\GTe(\bM,M)}&\asymp \frac{\te(W M)\bG_{\e}(\bM,W_{1})}{\bG_{\e}(\bM,M)}= \tte(W_{2})\tte(W_{1}M)\to 0
}
from Proposition \ref{prop:prop_mr_LAMPhe3}(5).
Consequently, we obtain the assertion.
\item[(2)]
For $\e>0$, $M\in \bM_{1}$, $M^\p\in \bM$ and $W=W_{1}W_{2}\cdots W_{n}\in \bM_{1}^{n-1}\times \bM$, we put
\ali
{
\LRR_{\e}(MM^\p)&=\LR_{A_{MM^\p},\Phe}(\lam(\e)\IR-\LR_{A_{MM},\Phe})^{-1}\\
\LRR_{\e}(W)&=\RLR_{\e}(W_{n-1}W_{n})\RLR_{\e}(W_{n-2}W_{n-1})\cdots\RLR_{\e}(W_{1}W_{2}).
}
Fix $M\in \bM_{1}$ and $f\in C(\SiAp)$. We consider the decomposition
\ali{
&\lam(\e)\nu(\bM,\e,\chi_{\SiM}g)=\LR_{A,\Phe}^{*}\nu(\bM,\e,\chi_{\SiM}g)\\
&\qquad\qquad=\sum_{M^\p\in \bM_{0}}\nu(\bM,\e,\LR_{A_{M M^\p},\Phe}g)\\
&\qquad\qquad\quad+\sum_{M^\p\in \bM_{1}\setminus \{M\}}\nu(\bM,\e,\LR_{A_{M M^\p},\Phe}g)+\nu(\bM,\e,\LR_{A_{MM},\Phe}g)
}
for $g\in C(\SiAp)$. Therefore, we have the equation
\ali{
\nu(\bM,\e,\chi_{\SiM}f)=&\sum_{M^\p\in \bM_{0}}\nu(\bM,\e,\LRR_{\e}(M M^\p)f)+\sum_{M^\p\in \bM_{1}\setminus \{M\}}\nu(\bM,\e,\LRR_{\e}(M M^\p)f)
}
by putting $g=(\lam(\e)\IR-\LR_{A_{MM},\Phe})^{-1}f$. By using the above repeatedly, we obtain
\alil{
\nu(\bM,\e,\SiM f)
&=\sum_{M^\p\in \bM_{0}}\sum_{i=0}^{m_{1}-1}\sum_{W\in W_{i}(\mM_{M})}\nu(\bM,\e,\LRR_{\e}(M W M^\p)f)\nonumber\\
&\ +\sum_{i=0}^{m_{1}-1}\sum_{W\in W_{i}(\mM_{M})}\nu(\bM,\e,\LRR_{\e}(M W M)f)\label{eq:speed_tsc_3}\\
&\ +\sum_{W\in W_{m_{1}}(\mM_{M})}\nu(\bM,\e,\LRR_{\e}(M W)f).\nonumber
}
We will estimate $\lim_{\e\to 0}\nu(\bM,\e,\SiM f)/\VTe(\bM,M)$. For $0\leq i\leq m_{1}-1$ and $W\in W_{i}(\mM_{M})$, if $M W M^\p$ is not simple path, there exist $W_{1},W_{2}\in W_{*}(\mM_{M})$ such that $W_{2}$ is cycle, $M W_{1} M^{\p}$ is in $\SP(M,M^\p:\bM_{1})$ and the inequality $|\nu(\bM,\e,\LRR_{\e}(M W M^\p)f)| \leq c \te(W_{2})\te(M W_{1} M^\p)\bV_{\e}(\bM,M)\|f\|_{\infty}$ is satisfied for some positive constant $c$. Thus we see
\ali{
\frac{|\nu(\bM,\e,\LRR_{\e}(M W M^\p)f)|}{\VTe(\bM,M)}&\leq  \frac{c\te(W_{2})\te(M W_{1} M^\p)\bV_{\e}(\bM,M^\p)\|f\|_{\infty}}{\VTe(\bM,M)}\\
&\leq c\te(W_{2})\|f\|_{\infty}\to 0.
}
If $M W M^\p$ is a simple path, then we have
\ali{
\frac{\nu(\bM,\e,\LRR_{\e}(M W M^\p)f)}{\VTe(\bM,M)}&=\frac{\te(M W M^\p)\bV_{\e}(\bM,M^\p)}{\VTe(\bM,M)}\frac{\nu(\bM,\e,\LRR_{\e}(M W M^\p)f)}{\te(M W M^\p)\bV_{\e}(\bM,M^\p)}\\
&\to c_{2}(M W M^\p)\nu(M^\p,\RLR(M W M^\p)f),\\
\frac{|\nu(\bM,\e,\LRR_{\e}(M W M)f)|}{\VTe(\bM,M)}&\leq c\te(M W M)\frac{\nu(\bM,\e,\SiM)}{\bV_{\e}(\bM,M)}\|f\|_{\infty}\\
&\leq cc^\p \te(M W M)\|f\|_{\infty}\to 0.
}
For $W\in W_{m_{1}}(\mM_{M})$, $M W$ contains of simple path and therefore, there exist $W_{1},W_{2}\in W_{*}(\mM_{M})$ such that $W_{2}$ is cycle, $W_{1} M\in \SP(M,M^\p:\bM_{1})$ and $\ttte(MW)=\ttte(W_{2})\ttte(M W_{1})$. Thus we see
\ali{
\frac{|\nu(\bM,\e,\LRR_{\e}(M W)f)|}{\VTe(\bM,M)}&\leq c\frac{\te(MW)\bV_{\e}(\bM,W_{1})\|f\|_{\infty}}{\bV_{\e}(\bM,M)}\\
&=c\ttte(W_{2})\ttte(MW_{1})\|f\|_{\infty}\to 0
}
by $\ttte(W_{2})\to 0$ and $\ttte(MW_{1})$ is bounded from Proposition \ref{prop:prop_mr_LAMPhe3}(6). 
Consequently, we obtain the assertion.
\proe
\subsection{Matrix representation of Ruelle operators (II)}\label{sec:matrep_LAi}
Assume the conditions $(\Si.1)$-$(\Si.3)$, $(\Ph.1)$-$(\Ph.3)$ and $(\Si.\bM)$ with $\bM\cap \bT_{0}\neq \emptyset$. For $M\in \bM\cap \bT_{0}$, we denotes $\bM_{M}=\{M\}\cup (\bM\cap \bT_{1})$ and set 
\alil{
\bbM=\bbM(B,\bM,\ph)=\{\bM_{M}\,:\,M\in \bM\cap \bT_{0}\}.\label{eq:bbM=}
}
In this section, we study a spectral matrix representation of $\LR_{A(\bM),\Phe}$ by using $\LR_{A(\bM_{M}),\Phe}$ for $M\in \bM\cap \bT_{0}$.
We notice that each $\bM_{M}$ satisfies $\sharp (\bM_{M}\cap \bT_{0})=1$ and the condition $(\Si.\bM)$ (Proposition \ref{prop:tsc_LAMPhe_T0=1}(1)). Therefore Proposition \ref{prop:tsc_LAMPhe_T0=1}(2) implies that the TSC $(\lam(\bM_{M},\e),$ $h(\bM_{M},\e,\cd),$ $\nu(\bM_{M},\e,\cd))$ converges to $(\lam,h(\bM_{M},\cd),\nu(\bM_{M},\cd))$ in $\R\times C(\SiAp) \times M(\SiAp)$ for $M\in \bM\cap \bT_{0}$.
We put 
\begin{align}
f_{MM^{\prime}}=g_{M^{\prime}M}=
\begin{cases}
\sum_{M^\p\in \bM_{M}}\chi_{\SiMp},&M^\p=M\\
\SiMp,&M^{\prime}\neq M\\
\end{cases}
\end{align}
for $M,M^\p\in \bM\cap \bT_{0}$.
Take a \MR
\alil{
(\mV_{\bbM,\e},\mG_{\bbM,\e},\bV_{\bbM,\e},\bG_{\bbM,\e})
}
of $(\LR_{A(\bM),\Phe},$ $\lam(\bM,\e),$ $g(\bM,\e,\cd),$ $\nu(\bM,\e,\cd))$ using $((\LR_{A(\bM_{M}),\Phe},$ $\lam(\bM_{M},\e),$ $h(\bM_{M},\e,\cd),$ $\nu(\bM_{M},\e,\cd))_{M}$, $(f_{MM^\p})$, $(g_{MM^\p}))$.
Furthermore, we have the normalized version as follows:
Assume that
\alil{
(\tmV_{\bbM,\e},\tmG_{\bbM,\e},\tbV_{\bbM,\e},\tbG_{\bbM,\e})
}
is the \MR given by (\ref{eq:tmVkkp=}) and (\ref{eq:tmGkkp=}) formed with $((\LR_{A(\bM),\tPhe},$ $\lam(\bM,\e),$ $\ti{g}(\bM,\e,\cd),$ $\ti{\nu}(\bM,\e,\cd))$, $(\LR_{A(\bM_{M}),\tPhe},$ $\lam(\bM_{M},\e),$ $\tilde{h}(\bM_{M},\e,\cd),$ $\tilde{\nu}(\bM_{M},\e,\cd))_{M}$. 
Therefore, we have
\ali{
&\bV_{\bbM,\e}(M)=\nu(\bM,\e,h(\bM_{M},\e,\cd))\\
&\mV_{\bbM,\e}(MM^\p)\\
=&\begin{cases}
\lam(\bM_{M},\e),&M=M^\p\\
\displaystyle \frac{\nu(\bM,\e,\sum_{M^\pp\in \bM_{M}}\LR_{A_{M^\pp M^\p},\Phe}h(\bM_{M},\e,\cd))}{\nu(\bM,\e,h(\bM_{M^\p},\e,\cd))},&M\neq M^\p \text{ and }\bV_{\bbM,\e}(M^\p)>0\\
0,& \text{otherwise}
\end{cases}\\
&\bG_{\bbM,\e}(M)=\nu(\bM_{M}, \e,g(\bM,\e,\cd))\\
&\mG_{\bM,\e}(MM^\p)\\
=&\begin{cases}
\lam(\bM_{M},\e),&M=M^\p\\
\displaystyle \frac{\nu(\bM_{M^\p},\e,\sum_{M^\pp\in \bM_{M^\p}}\LR_{A_{MM^\pp},\Phe}g(\bM,\e,\cd))}{\nu(\bM_{M},\e,g(\bM,\e,\cd))},&M\neq M^\p \text{ and }\bG_{\bbM,\e}(M)>0\\
0,& \text{otherwise}.
\end{cases}
}
Similarity, it follows from Proposition \ref{prop:mr_prop_nor} and Proposition \ref{prop:mr_prop_nor2} that
\begin{align}
\tmV_{\bbM, \e}(MM)&
=\tmG_{\bbM, \e}(MM)=\lam(\bM_{M},\e)\\
\tmV_{\bbM, \e}(MM^\p)&
=\mV_{\bbM, \e}(MM^\p)\frac{\nu(\bM_{M},\e,g(\e,\cd))}{\nu(\bM_{M^\p},\e,g(\e,\cd))} \text{ for }M\neq M^\p,\label{eq:tmV_Me=}\\
\tbV_{\bbM, \e}(M)&
=\bV_{\bbM, \e}(M)\frac{\nu(\bM_{M},\e,g(\e,\cd))}{\nu(\bM,\e,g(\e,\cd))}\\
\tmG_{\bbM, \e}(MM^\p)&
=\mG_{\bbM, \e}(MM^\p)\frac{\nu(\bM_{M},\e,g(\e,\cd))}{\nu(\bM_{M^\p},\e,g(\e,\cd))} \text{ for }M\neq M^\p\label{eq:tmG_Me=}\\
\tbG_{\bbM, \e}(M)&
=\bG_{\bbM, \e}(M)\frac{\nu(\bM,\e,g(\e,\cd))}{\nu(\bM_{M},\e,g(\e,\cd))\nu(\bM,g(\bM,\cd))}.
\end{align}
Put
\begin{align}
\te(M, M^{\prime}:\bT^{\prime})=\sum_{W\in \SP(M,M^{\prime}:T^{\prime})}\te(W)\\
\tte(M, M^{\prime}:\bT^{\prime})=\sum_{w\in \SP(M,M^{\prime}:T^{\prime})}\tte(W),\quad \ttte(M, M^{\prime}:\bT^{\prime})=&\sum_{W\in \SP(M,M^{\prime}:\bT^{\prime})}\ttte(W)
\end{align}
for a nonempty subset $\bT^{\prime}\subset \bT$ and $M,M^\p\in \bT$, where $\tte(W)=\tte(\bT,W)$ and $\ttte(W)=\ttte(\bT,W)$.
\prop
{\label{prop:thMM->}
Assume the conditions $(\Si.1)$-$(\Si.3)$, $(\Ph.1)$-$(\Ph.3)$ and $(\Si.\bM)$ with $\bM\cap \bT_{0}\neq \emptyset$. Then we have the following for $M\in \bM\cap \bT_{0}$:
\ite
{
\item $\ti{h}(\bM_{M},\e,\cd)\to 1$ on $\SiM$.
\item $\ti{h}(\bM_{M},\e,\cd)\asymp \tte(M, M^\p:\bM\cap \bT_{1})$ on $\SiMp$ for each $M^\p\in \bM\cap \bT_{1}$.
}
}
\pros
\item[(1)] Note the form
\ali
{
&\ti{h}(\bM_{M},\e,\cd)=\frac{\nu(\bM_{M},\e,g(\e,\cd))}{g(\e,\cd)}h(\bM_{M},\e,\cd)\\
=&\left(\frac{\nu(\bM_{M},\e,g(\e,\cd)\chi_{\SiM})}{\bG_{\e}(M)}+\frac{\sum_{M^\p\in \bM\cap \bT_{1}}\nu(\bM_{M},\e,g(\e,\cd)\chi_{\SiMp})}{\bG_{\e}(M)}\right)\frac{h(\bM_{M},\e,\cd)}{g(\e,\cd)/\bG_{\e}(M)}.
}
We have that for $M^\p\in \bM\cap \bT_{1}$
\alil
{
\frac{\nu(\bM_{M},\e,g(\e,\cd)\chi_{\SiMp})}{\bG_{\e}(M)}
&\asymp \bV_{\e}(\bM_{M},M^\p)\frac{\bG_{\e}(M^\p)}{\bG_{\e}(M)}\nonumber\\
&\asymp \sum_{W\in \SP(M^\p,M:\bT_{1}(B(\bM_{M}),\ph))}\te(W)\bV_{\e}(\bM_{M},M)\frac{\bG_{\e}(M^\p)}{\bG_{\e}(M)}\label{eq:lem:thMM->_1}\\
&\asymp \tte(M^\p, M:\bM\cap \bT_{1})\nonumber
}
by using Proposition \ref{prop:prop_mr_LAMPhe3}(3) and the facts $\bT_{1}(B(\bM_{M}),\ph)=\bM\cap \bT_{1}$ and the relation $\bV_{\e}(\bM_{M},M)\asymp 1$. When $W=W_{1}W_{2}\cdots W_{n}\in \SP(M^\p,M:\bM\cap \bT_{1})$, we have $\tte(W)\to 0$ from $\tte(W_{i-1}W_{i})=O(1)$ $(2\leq i\leq n-1)$ and $\tte(W_{n-1}W_{n})\to 0$ by Proposition \ref{prop:prop_mr_LAMPhe3}(5). Thus we see $\tte(M^\p, M:\bM\cap \bT_{1})\to 0$ for $M^\p \in \bM\cap \bT_{1}$. By virtue of Proposition \ref{prop:tsc_LAMPhe_T0=1}(2)(3), we obtain
\ali
{
\ti{h}(\bM_{M},\e,\cd)&\to \frac{\nu(\bM_{M},h(M,\cd))}{h(M,\cd)}h(\bM_{M},\cd)\\
&=\frac{\nu(\bM_{M},\SiM)\nu(M,h(M,\cd))}{h(M,\cd)}\frac{h(M,\cd)}{\nu(\bM_{M},\SiM)}=1\quad\text{ on }\SiM.
}
\item[(2)] For $M^\p\in \bM\cap \bT_{1}$, on $\SiMp$
\ali
{
&\ti{h}(\bM_{M},\e,\cd)\\
=&\frac{\nu(\bM_{M},\e,g(\e,\cd)\chi_{\SiM})+\sum_{M^\pp\in \bM\cap \bT_{1}}\nu(\bM_{M},\e,g(\e,\cd)\chi_{\SiMpp})}{g(\e,\cd)}h(\bM_{M},\e,\cd)\\
\asymp& \frac{\bG_{\e}(M)+\sum_{M^\pp\in \bM\cap \bT_{1}}\tte(M^\pp,M:\bM\cap \bT_{0})\bG_{\e}(M)}{\bG_{\e}(M^\p)}\te(M, M^\p:\bM\cap\bT_{1})\\
=& \left(1+\sum_{M^\pp\in \bM\cap\bT_{1}}\tte(M^\pp,M:\bM\cap \bT_{0})\right)\tte(M, M^\p:\bM\cap\bT_{1})\\
\asymp&\tte(M, M^\p:\bM\cap\bT_{1})
}
is satisfied by using Proposition \ref{prop:prop_mr_LAMPhe3}(3)(4) and the relation (\ref{eq:lem:thMM->_1}).
\proe
Consequently, convergence of the measure $\ti{\nu}(\bM,\e,\cd)$ is reduced to convergence of the right Perron eigenvector of the matrix $\mV_{\bbM,\e}$ as follows. Note that when $\bM=\bT$, $\ti{\nu}(\bM,\e,\cd)$ is equal to the Gibbs measure $\mu(\e,\cd)$.
\thm{\label{th:conv_gibbs_iff_mVM}
Assume that the conditions $(\Sigma.1)$-$(\Sigma.3)$, $(\varPhi.1)$-$(\varPhi.3)$ and $(\Si.\bM)$ with $\bM\cap \bT_{0}\neq \emptyset$ are satisfied. Then $\ti{\nu}(\bM,\e,\cd)$ converges to a measure $\mu$ if and only if $\tbV_{\bbM,\e}(M)$ converges to a number $\delta(M)$ for all $M\in \bM\cap \bT_{0}$. In these cases, $\mu$ has the form $\mu=\sum_{M\in \bM\cap \bT_{0}}\delta(M)\mu(M,\cd)$. In particular, $\delta(M)=\mu(\SiM)$ for each $M\in \bM\cap \bT_{0}$ and  $\supp \mu\subset \bigcup_{M\in \bM\cap \bT_{0}}\SiM$.
}
\pros
Assume that $\ti{\nu}(\bM,\e,\cd)$ converges to $\mu$. First we prove $\mu(\bigcup_{M\in \bM\cap\bT_{1}}\SiM)=0$. For $M\in \bM\cap\bT_{1}$, we see
\ali
{
\ti{\nu}(\bM,\e,\SiM)&=\frac{\nu(\bM,\e,g(\e,\cd)\chi_{\SiM})}{\nu(\bM,\e,g(\e,\cd))}
\asymp \frac{\bG_{\e}(M)\bV_{\e}(\bM,M)}{\nu(\bM,\e,g(\e,\cd))}\\
&\asymp \frac{\bG_{\e}(M)\sum_{M^\p\in \bM\cap\bT_{0}}\sum_{W\in \SP(M,M^\p:\bM\cap\bT_{1})}t_{\e}(W)\bV_{e}(\bM,M^\p)}{\nu(\bM,\e,g(\e,\cd))}\\
&=\frac{\sum_{M^\p\in \bM\cap\bT_{0}}\sum_{W\in \SP(M,M^\p:\bM\cap\bT_{1})}\tte(W)\bG_{\e}(M^\p)\bV_{e}(\bM,M^\p)}{\nu(\bM,\e,g(\e,\cd))}\\
&\leq \sum_{M^\p\in \bM\cap\bT_{0}:\nu(\bM,\e,\SiM)>0}\tte(M,M^\p:\bM\cap\bT_{1})\frac{\bG_{\e}(M^\p)\bV_{\e}(\bM,M^\p)}{\nu(\bM,\e,g(\e,\cd)\chi_{\SiM})}\\
&\asymp \sum_{M^\p\in \bM\cap\bT_{0}}\tte(M,M^\p:\bM\cap\bT_{1})\to 0
}
from Proposition \ref{prop:prop_mr_LAMPhe3}(3)(5). Therefore $\mu(\bigcup_{M\in \bM\cap\bT_{1}}\SiM)=0$ is valid. Next, we show that $\mu$ has the form $\sum_{M\in \bM\cap\bT_{0}}\delta(M)\mu(M,\cd)$ for some nonnegative constants $\delta(M)$ with $\sum_{M\in \bM\cap\bT_{0}}\delta(M)=1$.
Let $m(M,\e,\cd)\in M(\SiAp)$ and $f(M,\e,\cd)\in C(\SiAp)$ be $m(M,\e,f)=\nu(\bM,\e,\chi_{\SiM}f)/\nu(\bM,\e,\SiM)$ and $f(M,\e,\cd)=g(\e,\cd)\chi_{\SiM}/\nu(M,\e,g(\e,\cd))$. We obtain that for $M\in \bM\cap \bT_{0}$ with $\nu(\bM,\e,\SiM)>0$
\alil
{
\ti{\nu}(\bM,\e,f\chi_{\SiM})=&\frac{\nu(\bM,\e,g(\e,\cd)f\chi_{\SiM})}{\nu(\bM,\e,g(\e,\cd))}\nonumber\\
=&\frac{\nu(\bM,\e,g(\e,\cd)\chi_{\SiM})}{\nu(\bM,\e,g(\e,\cd))}
\frac{\nu(\bM,\e,\SiM)\nu(M,\e,g(\e,\cd)\chi_{\SiM})}{\nu(\bM,\e,g(\e,\cd)\chi_{\SiM})}\times\nonumber\\
&\quad\frac{1}{\nu(\bM,\e,\SiM)}\nu\left(\bM,\e,\frac{g(\e,\cd)\chi_{\SiM}}{\nu(M,\e,g(\e,\cd)\chi_{\SiM})}f\right)\label{eq:conv_gibbs_iff_mVM_2}\\
=&\ti{\nu}(\bM,\e,\chi_{\SiM})\frac{m(M,\e,f(M,\e,\cd)f)}{m(M,\e,f(M,\e,\cd))}.\nonumber
}
Letting as $\e\to 0$, $\mu$ has the form $\mu=\sum_{M\in \bM\cap\bT_{0}}\mu(\SiM)\mu(M,f)$ from Proposition \ref{prop:conv_nor_TSC_M}(1)(2).

Now we will prove convergence of $\tbV_{\bbM,\e}(M)$ under the assumption $\mu(\e,\cd)\to \mu$. Choose any positive sequence $(\e_{n})$ with $\inf_{n}\e_{n}=0$. The function $\ti{h}(\bM_{M},\e,\cd)$ is bounded from Proposition \ref{prop:thMM->}(1)(2), and $\ti{h}(\bM_{M},\e,\cd)/\|\ti{h}(\bM_{M},\e,\cd)\|_{\infty}$ belongs to $\Lambda_{c}$ with some constant $c$. Therefore there is a subsequence $(\e^\p_{n})$ of $(\e_{n})$ such that for each $M\in \bT_{0}$, $\ti{h}(\bM_{M},\e^\p_{n},\cd)$ converges to a function $\ti{h}$ in $C(\SiAp)$. Proposition \ref{prop:thMM->}(1) yields $\ti{h}=1$ on $\SiM$. We obtain
\alil
{
\tbV_{\bbM,\e}(M)=\ti{\nu}(\bM,\e,\ti{h}(\bM_{M},\e,\cd))\to \mu(\ti{h})=\sum_{M\in \bM\cap\bT_{0}}\mu(\ti{h}\chi_{\SiM})=\mu(\SiM)\label{eq:conv_gibbs_iff_mVM_1}
}
as $\e\to 0$ running through $(\e_{n}^{\p})$. By arbitrarily choosing $(\e_{n})$, convergence of $\tbV_{\bbM,\e}(M)$ is fulfilled. 

Conversely, we consider convergence of $\ti{\nu}(\bM,\e,\cd)$ under $\tbV_{\bbM,\e}(M)\to \delta(M)$ for each $M\in \bM\cap\bT_{0}$. Take a limit point $\mu$ of $\ti{\nu}(\bM,\e,\cd)$ in $M(\Sigma_{A}^{+})$. 
By virtue of (\ref{eq:conv_gibbs_iff_mVM_2}), $\mu$ has the form $\mu=\sum_{M\in \bM\cap\bT_{0}}\mu(\SiM)\mu(M,\cd)$. The fact (\ref{eq:conv_gibbs_iff_mVM_1}) implies $\mu(\Sigma_{M})=\delta(M)$. Hence the accumulation points of $\ti{\nu}(\bM,\e,\cd)$ consists of only one element $\sum_{M\in \bM\cap\bT_{0}}\mu(\SiM)\mu(M,\cd)$. The proof is complete.
\proe
\prop{\label{prop:geneform_prop2}
Assume the conditions $(\Si.1)$-$(\Si.3)$, $(\Ph.1)$-$(\Ph.3)$ and $(\Si.\bM)$ with $\bM\cap \bT_{0}\neq \emptyset$. Assume also that for each $M,M^\p\in \bT$, $\LR_{A_{MM^\p},\Phe}/\te(MM^\p)$ converges to an operator $\LR(MM^\p)$ in $\LR(C(\SiAp))$ and for each $M,M^\p\in \bM\cap\bT_{0}$ and $W\in \SP(M,M^\p:\bM\cap\bT_{1})$, $\te(W)/\te(M, M^\p:\bM\cap\bT_{1})$ converges to a number $c_{MM^\p}(W)$ in $[0,1]$.
\item[(1)] For all $M,M^\p\in \bM\cap\bT_{0}$ with $M\neq M^\p$ and $\te(M,M^\p:\bM\cap \bT_{1})>0$,
$$\frac{\mV_{\bbM,\e}(MM^\p)}{\te(M,M^\p:\bM\cap\bT_{1})}\to \sum_{W\in \SP(M,M^\p:\bM\cap\bT_{1})}c_{MM^\p}(W)\nu(\bM_{M^\p},\LRLR(W)h(\bM_{M},\cd)).$$
Here for $W=W_{1}W_{2}\cdots W_{n}\in \bT\times \bT_{1}^{n-2}\times \bT$, 
\ali
{
\LRLR(W)=&\LR(W_{n-1}W_{n})\RR(W_{n-1})\LR(W_{n-2}W_{n-1})\RR(W_{n-2})\cdots\\
&\quad\cdots \LR(W_{2}W_{3})\RR(W_{2})\LR(W_{1}W_{2}),
}
and $\RR(M)=(\lambda\mathcal{I}-\mathcal{L}_{M,\varphi})^{-1}$.
\item[(2)] For all $M,M^\p\in \bM\cap\bT_{0}$ with $M\neq M^\p$ and $\te(M,M^\p:\bM\cap \bT_{1})>0$,
\ali
{
\frac{\mG_{\bbM,\e}(MM^\p)}{\te(M,M^\p:\bM\cap\bT_{1})}\to \sum_{W\in \SP(M,M^\p:\bM\cap\bT_{1})}c_{MM^\p}(W)\nu(\bM_{M^\p},\LRLR(W)h(\bM_{M},\cd)).
}
}
\pros
Put $\bM_{k}=\bM\cap \bT_{k}$ and $m_{k}=\sharp \bM_{k}$ for $k=0,1$.
\item[(1)] Let $M,M^\p\in \bM_{0}$ with $M\neq M^\p$. We Recall
\ali
{
\mV_{\bbM,\e}&(MM^\p)=\frac{\nu(\bM,\e,\sum_{M^\pp\in \bM_{M}}\LR_{A_{M^\pp M^\p},\Phe}h(\bM_{M},\e,\cd))}{\nu(\bM,\e,h(\bM_{M^\p},\e,\cd))}\\
=&\frac{\nu(\bM,\e,\sum_{M^\pp\in \bM_{M}}\LR_{A_{M^\pp M^\p},\Phe}h(\bM_{M},\e,\cd))}{\bV_{\e}(\bM,M^\p)}\times\\
& \left(\frac{\nu(\bM,\e,\chi_{\SiMp}h(\bM_{M^\p},\e,\cd))}{\bV_{\e}(\bM,M^\p)}+\frac{\sum_{M^\pp\in \bM_{1}}\nu(\bM,\e,\chi_{\SiMpp}h(\bM_{M^\p},\e,\cd))}{\bV_{\e}(\bM,M^\p)}\right)^{-1}.
}
We will show that the last expression converges to $\nu(M^\p,h(\bM_{M^\p},\cd))^{-1}$. We have that for $M^\pp\in \bM_{1}$
\ali
{
\frac{\nu(\bM,\e,\chi_{\SiMpp}h(\bM_{M^\p},\e,\cd))}{\bV_{\e}(\bM,M^\p)}\asymp& \frac{\nu(\bM,\e,\chi_{\SiMpp})}{\bV_{\e}(\bM,M^\p)}\te(M^\p,M^\pp:\bM_{1})\\
=&\frac{\nu(\bM,\e,\chi_{\SiMpp})}{\bV_{\e}(\bM,M^\pp)}\ttte(M^\p,M^\pp:\bM_{1})\to 0
}
is satisfied from Proposition \ref{prop:prop_mr_LAMPhe3}(4), convergence $\ttte(M^\p,M^\pp:\bM_{1})\to 0$ with Proposition \ref{prop:prop_mr_LAMPhe3}(6), and the fact $\nu(\bM,\e,\chi_{\SiMpp})\asymp \bV_{\e}(\bM,M^\pp)$. Moreover, we see
\ali
{
\frac{\nu(\bM,\e,\chi_{\SiMp}h(\bM_{M^\p},\e,\cd))}{\bV_{\e}(\bM,M^\p)}\to \nu(M^\p,h(\bM_{M^\p},\cd)).
}
On the other hand, we have
\alil
{
&\sum_{M^\pp\in \bM_{M}}\LR_{A_{M^\pp M^\p},\Phe}h(\bM_{M},\e,\cd)\nonumber\\
=&\LR_{A_{MM^\p},\Phe}h(\bM_{M},\e,\cd)+\sum_{M^\pp\in \bM_{1}}\LR_{A_{M^\pp M^\p},\Phe}h(\bM_{M},\e,\cd)\nonumber\\
=&\LR_{A_{MM^\p},\Phe}h(\bM_{M},\e,\cd)\nonumber\\
&+\sum_{M^\pp\in \bM_{1}}\sum_{W\in W_{m_{1}}(\mM_{M^\pp})}\LR_{A_{M^\pp M^\p},\Phe}\RLR_{\e}(W M^\pp,h(\bM_{M},\e,\cd))\nonumber\\
&+\sum_{M^\pp\in \bM_{1}}\sum_{i=1}^{m_{1}-1}\sum_{W\in W_{i}(\mM_{M^\pp})}\LR_{A_{M^\pp M^\p},\Phe}\RLR_{\e}(M W M^\pp)h(\bM_{M},\e,\om)\label{eq:prop:geneform_prop2_1}\\
&+\sum_{M^\pp\in \bM_{1}}\sum_{i=1}^{m_{1}-1}\sum_{W\in W_{i}(\mM_{M^\pp})}\LR_{A_{M^\pp M^\p},\Phe}\RLR_{\e}(M^\pp W M^\pp)h(\bM_{M},\e,\om)\nonumber\\
&+\sum_{M^\pp\in \bM_{1}}\LR_{A_{M^\pp M^\p},\Phe}\RLR_{\e}(MM^\pp)h(\bM_{M},\e,\om)\nonumber\\
=&I(\e)+II(\e)\nonumber
}
with
\ali
{
I(&\e)=\LR_{A_{MM^\p},\Phe}h(\bM_{M},\e,\cd)\\
&+\sum_{M^\pp\in \bM_{1}}\sum_{i=1}^{m_{1}-1}\sum_{W\in W_{i}(\mM_{M^\pp}):\atop{MWM^\pp \text{ is a simple path}}}\LR_{A_{M^\pp M^\p},\Phe}\RLR_{\e}(M W M^\pp)h(\bM_{M},\e,\om)\\
&+\sum_{M^\pp\in \bM_{1}}\LR_{A_{M^\pp M^\p},\Phe}\RLR_{\e}(MM^\pp)h(\bM_{M},\e,\om)
}
and
\ali
{
II(&\e)=\sum_{M^\pp\in \bM_{1}}\sum_{W\in W_{m_{1}}(\mM_{M^\pp}):WM^\pp}\LR_{A_{M^\pp M^\p},\Phe}\RLR_{\e}(W M^\pp,h(\bM_{M},\e,\cd))\\
&+\sum_{M^\pp\in \bM_{1}}\sum_{i=1}^{m_{1}-1}\sum_{W\in W_{i}(\mM_{M^\pp}):\atop{MWM^\pp \text{ contains a cycle}}}\LR_{A_{M^\pp M^\p},\Phe}\RLR_{\e}(M W M^\pp)h(\bM_{M},\e,\om)\\
&\ +\sum_{M^\pp\in \bM_{1}}\sum_{i=1}^{m_{1}-1}\sum_{W\in W_{i}(\mM_{M^\pp})}\LR_{A_{M^\pp M^\p},\Phe}\RLR_{\e}(M^\pp W M^\pp)h(\bM_{M},\e,\om)
}
from (\ref{eq:speed_tsc_2}) by replacing $\bM=\bM_{M}$,
$g(\bM,\e,\cd)=h(\bM_{M},\e,\cd)$.
To estimate $I(\e)$ and $II(\e)$, we note the equation
\ali
{
\te(M,M^\p:\bM_{1})=\te(MM^\p)+\sum_{M^\pp\in \bM_{1}}\te(M,M^\pp:\bM_{1})\te(M^\pp M^\p).
}
We have that for $M^\pp\in \bM_{1}$ and $W\in W_{i}(\mM_{M^\pp})$ with $1\leq i\leq m_{1}$, if $W M^\pp$ contains a cycle, $\te(M,M^\pp:\bM_{1})>0$ and $\te(M^\pp M^\p)>0$, then $\bG_{\e}(\bM_{M},M^\pp)>0$ is satisfied and
\ali
{
&\frac{\|\LR_{A_{M^\pp M^\p},\Phe}\RLR_{\e}(W M^\pp)h(\bM_{M},\e,\cd)\|_{\infty}}{\te(M,M^\p:\bM_{1})}\\
\leq &\frac{\|\LR_{A_{M^\pp M^\p},\Phe}\RLR_{\e}(W M^\pp)h(\bM_{M},\e,\cd)\|_{\infty}}{\te(M,M^\pp:\bM_{1})\te(M^\pp M^\p)}\\
=&\frac{\|\LRLR_{\e}(W M^\pp M^\p)h(\bM_{M},\e,\cd)\|_{\infty}}{\te(M,M^\pp:\bM_{1})\te(M^\pp M^\p)}\\
\leq &c \frac{\te(W M^\pp M^\p)\te(M,W_{1}:\bM_{1})}{\te(M,M^\pp:\bM_{1})\te(M^\pp M^\p)}\asymp c \frac{\te(WM^\pp) \bG_{e}(\bM_{M},W_{1})}{\bG_{\e}(\bM_{M},M^\pp)}\to 0
}
for some positive constant $c$ with Proposition \ref{prop:prop_mr_LAMPhe3}(5). Therefore we obtain $|II(\e)|/\te(M,M^\p:\bM_{1})\to 0$. On the other hand, we have
\ali
{
\frac{I(\e)}{\te(M,M^\p:\bM_{1})}\to& c_{MM^\p}(MM^\p)\LR_{MM^\p}h(\bM_{M},\cd)\\
&+\sum_{W\in \SP(M,M^\p:\bM_{1})\,:\,|W|\geq 3}c_{MM^\p}(W)\LRLR(W)h(\bM_{M},\cd)\\
=&\sum_{W\in \SP(M,M^\p:\bM_{1})}c_{MM^\p}(W)\LRLR(W)h(\bM_{M},\cd).
}
Thus 
\ali
{
\frac{\mV_{\bbM,\e}(MM^\p)}{\te(M,M^\p:\bM_{1})}&\to \frac{\nu(M^\p,\sum_{W\in \SP(M,M^\p:\bM_{1})}c_{MM^\p}(W)\LRLR(W)h(\bM_{M},\cd))}{\nu(M^\p,h(\bM_{M^\p},\cd))}\\
&=\frac{\nu(M^\p,\sum_{W\in \SP(M,M^\p:\bM_{1})}c_{MM^\p}(W)\LRLR(W)h(\bM_{M},\cd))}{\nu(M^\p,h(M^\p,\cd))/\nu(\bM_{M^\p},\SiMp)}\\
&=\sum_{W\in \SP(M,M^\p:\bM_{1})}c_{MM^\p}(W)\nu(\bM_{M^\p},\LRLR(W)h(\bM_{M},\cd)).
}
\item[(2)] 
We have
\ali
{
\mG_{\bbM,\e}&(MM^\p)=\frac{\nu(\bM_{M^\p},\e,\sum_{M^\pp\in \bM_{M^\p}}\LR_{A_{MM^\pp},\Phe}g(\bM,\e,\cd))}{\nu(\bM_{M},\e,g(\bM,\e,\cd))}\\
&=\frac{\nu(\bM_{M^\p},\e,\sum_{M^\pp\in \bM_{M^\p}}\LR_{A_{MM^\pp},\Phe}g(\bM,\e,\cd))}{\bG_{\e}(\bM,M)}\times\\
&\ \left(\frac{\nu(\bM_{M},\e,\chi_{\SiM}g(\bM,\e,\cd))}{\bG_{\e}(\bM,M)}+\frac{\sum_{M^\pp\in \bM_{1}}\nu(\bM_{M},\e,\chi_{\SiMpp}g(\bM,\e,\cd))}{\bG_{\e}(\bM,M)}\right)^{-1}.
}
By a similar argument above (1), the last expression converges to the number $\nu(\bM_{M},h(M,\cd))^{-1}$ by using Proposition \ref{prop:prop_mr_LAMPhe3}(5). Let $f_{\e}=g(\bM,\e,\cd)/\bG_{\e}(\bM,M)$. Moreover, we have the following from using the expansion like (\ref{eq:prop:geneform_prop2_1}):
\ali
{
\sum_{M^\pp\in \bM_{M^\p}}\nu(\bM_{M},\e,\LR_{A_{MM^\pp},\Phe}f_{\e})
=&III(\e)+IV(\e)
}
with
\ali
{
&III(\e)=\nu(\bM_{M},\e,\LR_{A_{MM^\p},\Phe}f_{\e})\\
&\ +\sum_{M^\pp\in \bM_{1}}\sum_{i=1}^{m_{1}-1}\sum_{W\in W_{i}(\mM_{M^\pp}):\atop{M^\pp WM^\p \text{ is a simple cycle}}}\nu(\bM_{M^\p},\e,\LRR_{\e}(M^\pp WM^\p)\LR_{A_{MM^\pp},\Phe}f_{\e})\\
&\ +\sum_{M^\pp\in \bM_{1}}\nu(\bM_{M},\e,\LRR_{\e}(M^\pp M^\p)\LR_{A_{MM^\pp},\Phe}f_{\e}),\\
&IV(\e)=\sum_{M^\pp\in \bM_{1}}\sum_{W\in W_{m_{1}}(\mM_{M^\pp})}\nu(\bM_{M^\p},\e,\LRR_{\e}(M^\pp W)\LR_{A_{MM^\pp},\Phe}f_{\e})\\
&\ +\sum_{M^\pp\in \bM_{1}}\sum_{i=1}^{m_{1}-1}\sum_{W\in W_{i}(\mM_{M^\pp}):\atop{M^\pp WM^\p \text{ contains a cycle}}}\nu(\bM_{M^\p},\e,\LRR_{\e}(M^\pp WM^\p)\LR_{A_{MM^\pp},\Phe}f_{\e})\\
&\ +\sum_{M^\pp\in \bM_{1}}\sum_{i=1}^{m_{1}-1}\sum_{W\in W_{i}(\mM_{M^\pp})}\nu(\bM_{M^\p},\e,\LRR_{\e}(M^\pp WM^\pp)\LR_{A_{MM^\pp},\Phe}f_{\e}).
}
Note that
\ali
{
\te(M,M^\p:\bM_{1})=\te(MM^\p)+\sum_{M^\pp\in \bM_{1}}\te(MM^\pp)\te(M^\pp,M^\p:\bM_{1}).
}
For $M^\pp\in \bM_{1}$ and $W\in W_{i}(\mM_{M^\pp})$ with $1\leq i\leq m_{1}$, if$M^\pp W$ contains a cycle, $\te(MM^\pp)\te(M^\pp,M^\p:\bM_{1})>0$ and $\te(MM^\pp)>0$, then we have
\ali
{
&\frac{|\nu(\bM_{M},\e,\LRR_{\e}(M^\pp W)\LR_{A_{MM^\pp},\Phe}f_{\e})|}{\te(M,M^\p:\bM_{1})}\\
\leq &\frac{|\nu(\bM_{M},\e,\LRR_{\e}(M^\pp W)\LR_{A_{MM^\pp},\Phe}f_{\e})|}{\te(MM^\pp)\te(M^\pp,M^\p:\bM_{1})}\\
=&\frac{|\nu(\bM_{M},\e,\LRLR_{\e}(MM^\pp W)f_{\e})|}{\te(MM^\pp)\te(M^\pp,M^\p:\bM_{1})}\\
\leq &c\frac{\te(M^\pp W)\te(W_{|W|},M^\p:\bM_{1})}{\te(M^\pp,M^\p:\bM_{1})}\asymp c\frac{\te(M^\pp W)\bV_{\e}(\bM_{M},W_{|W|})}{\bV_{\e}(\bM_{M},M^\pp)}\to 0
}
from $WM^\pp$ contains a cycle for some constant $c$. Therefore we have convergence $|IV(\e)|/\te(M,M^\p:\bM_{1})\to 0$. On the other hand, we obtain
\ali
{
\frac{\mG_{\bbM,\e}(MM^\p)}{\te(M,M^\p:\bM_{1})}
&\to \frac{c_{MM^\p}(MM^\p)\nu(\bM_{M},\LR_{MM^\p}h(M,\cd))}{\nu(\bM_{M},h(M,\cd))}\\
&\ +\frac{\sum_{W\in \SP(M,M^\p:\bM_{1}\,:\,|W|\geq 3}c_{MM^\p}(W)\nu(\bM_{M},\LRLR(W)h(M,\cd))}{\nu(\bM_{M},h(M,\cd))}\\
&=\sum_{W\in \SP(M,M^\p:\bM_{1})}c_{MM^\p}(W)\nu(\bM_{M},\LRLR(W)h(\bM_{M},\cd))
}
by using Proposition \ref{prop:tsc_LAMPhe_T0=1}(3).
\proe
\cor
{\label{cor:geneform_prop2}
Assume the conditions $(\Si.1)$-$(\Si.3)$,  $(\Ph.1)$-$(\Ph.3)$ and $(\Si.\bM)$ with $\bM\cap \bT_{0}\neq \emptyset$. Then for all $M,M^\p\in \bM\cap\bT_{0}$ with $M\neq M^\p$ and $\te(M,M^\p:\bM\cap\bT_{1})>0$,  $\tmV_{\bbM,\e}(MM^\p)/\tmG_{\bbM,\e}(MM^\p)$ converges to $1$.
}
\pros
By the proof of Proposition \ref{cor:geneform_prop2}, we see $\mV_{\bbM,\e}(MM^\p)\asymp \te(M,M^\p:\bM\cap\bT_{1})\asymp \mG_{\bbM,\e}(MM^\p)$. Therefore, $\te(M,M^\p:\bM\cap\bT_{1})>0$ implies $\mV_{\bbM,\e}(MM^\p)>0$ and $\mG_{\bbM,\e}(MM^\p)>0$. Choose any positive sequence $(\e_{n})$ with $\inf_{n}\e_{n}=0$. We can take a subsequence $(\e_{n}^{\p})$ of $(\e_{n})$ satisfying that $\LR_{A_{MM^\p},\Phe}/\te(MM^\p)$ converges for all $M,M^\p\in \bM\cap\bT$ with $\te(MM^\p)>0$, and $\te(W)/\te(M,M^\p:\bM\cap\bT_{1})$ converges for all $W\in \SP(M,M^\p : \bM\cap\bT_{1})$. 
By virtue of Proposition \ref{cor:geneform_prop2} again and by the equations (\ref{eq:tmV_Me=}) and (\ref{eq:tmG_Me=}), we obtain $\tmV_{\bbM,\e}(MM^\p)/\tmG_{\bbM,\e}(MM^\p)=\mV_{\bbM,\e}(MM^\p)/\mG_{\bbM,\e}(MM^\p)\to 1$ as $\e\to 0$ running through $(\e_{n}^{\p})$. By arbitrary choosing $(\e_{n})$, we see the assertion.
\proe
\section{Proof of main theorems}\label{sec:proof}
This section is devoted to proofs of main results given in Section \ref{sec:intro}.
\subsection{Proof of Theorem \ref{th:conv_pressure}}\label{sec:lim_pre_Gibbs}
The assertion directly follows from Proposition \ref{prop:conv_pre_M} by putting $\bM=\bT$.\qed
\subsection{Proof of Theorem \ref{th:limpoint_Gibbs_entropy}}\label{sec:proof_limpoint_Gibbs_entropy}
By virtue of Theorem \ref{th:conv_gibbs_iff_mVM} with $\bM=\bT$, if $\mu(\e,\cd)$ converges to a measure $\mu$ weakly, then we obtain the form $\mu=\sum_{M\in \bT_{0}}\delta(M)\mu(M,\cd)$ for some constants $\delta(M)$.
It is sufficient to show that the last assertion holds. Recall the equation
\ali
{
h(\si_{A},\mu(\e,\cd))&=P(\si_{A},\Phe)-\int \Phe d\mu(\e,\cd)\\
&=P(\si_{A},\Phe)-\int \phe d\mu(\e,\cd)-\int \psi(\e,\cd)\chi_{N} d\mu(\e,\cd)
}
Now we show $\int \psi(\e,\cd)\chi_{N} d\mu(\e,\cd)\to 0$. Let $u(\e,\cd)=\psi(\e,\cd)\chi_{N}-\log g(\e,\cd)\circ \si_{A}+\log g(\e,\cd)$. We have
\ali
{
\int \psi(\e,\cd)\chi_{N} d\mu(\e,\cd)=&\int u(\e,\cd) d\mu(\e,\cd)\\
=&\lam(\e)^{-1}\int \LR_{A,\tPhe}u(\e,\cd) d\mu(\e,\cd)\\
=&\lam(\e)^{-1}\int \LR_{A,\phe}(e^{u(\e,\cd)}\log(e^{u(\e,\cd)})) d\mu(\e,\cd)\\
=&\lam(\e)^{-1}\sum_{M,M^\p\in \bT}\int \LR_{A_{MM^\p},\phe}(e^{u(\e,\cd)}\log(e^{u(\e,\cd)})) d\mu(\e,\cd)
}
from $\mu(\e,\cd)$ is $\si_{A}$-invariant.
Note $e^{u(\e,\cd)}\asymp e^{\tPhe}\asymp \tte(MM^\p)$ on $\Si_{MM^\p}$.
Therefore, 
\ali
{
&e^{u(\e,\cd)}\log(e^{u(\e,\cd)})\\
\in&\Big[\min(c \tte(MM^\p)\log (c^{-1}\tte(MM^\p)),c^{-1}\tte(MM^\p)\log (c^{-1}\tte(MM^\p))),\\
&\qquad\qquad\max(c \tte(MM^\p)\log (c\tte(MM^\p)),c^{-1}\tte(MM^\p)\log (c\tte(MM^\p)))\Big]
}
on $\Si_{MM^\p}$ for some constant $c\geq 1$. 
For $M^\p\in \bT_{0}$ and $M\in \bT$ with $M\neq M^\p$, $e^{u(\e,\cd)}\log(e^{u(\e,\cd)})$ converges to $0$ uniformly on $\Si_{MM^\p}$ by using $\tte(MM^\p)\to 0$ ( Proposition \ref{prop:prop_mr_LAMPhe3}(5) ). For $M^\p\in \bT_{1}$, $e^{u(\e,\cd)}\log(e^{u(\e,\cd)})$ is bounded on $\Si_{MM^\p}$ by $\tte(MM^\p)=O(1)$. Furthermore, when $M=M^\p\in \bT_{0}$, $\e^{u(\e,\cd)}$ converges to $0$ on $\Si_{MM}\cap N$ by the condition $(\Ph.2)$.
By using $\mu(\e,\sum_{M\in \bT_{1}}\chi_{\SiM})\to 0$, we obtain
\ali
{
&\left|\int \psi(\e,\cd)\chi_{N} d\mu(\e,\cd)\right|\\
\leq& \lam(\e)^{-1}\sum_{M,M^\p\in \bT}\mu(\e,\SiMp)\|\LR_{A_{MM^\p},\phe}1\|_{\infty}\|\chi_{\Si_{MM^\p}}e^{u(\e,\cd)}\log(e^{u(\e,\cd)})\|_{\infty}\to 0.
}
Thus $\int \psi(\e,\cd)\chi_{N} d\mu(\e,\cd)\to 0$. Consequently, we give
\ali
{
&h(\si_{A},\mu(\e,\cd))\\
\to& P(\si_{B},\ph_{B})-\sum_{M\in \bT_{0}}\delta(M)\int \ph_{M} d\mu(M,\cd)\\
=&\sum_{M\in \bT_{0}}\delta(M)\left(P(\si_{M},\ph_{M})-\int \ph_{M} d\mu(M,\cd)\right)=\sum_{M\in \bT_{0}}\delta(M)h(\si_{M},\mu(M,\cd))
}
by $P(\si_{M},\ph_{M})=P(\si_{B},\ph_{B})$ for $M\in \bT_{0}$.
\qed
\subsection{Proof of Theorem \ref{th:m0=2}}
Assume $(\Si.1)$-$(\Si.3)$ and $(\Ph.1)$-$(\Ph.3)$. We give a set $\bbM=\bbM(B,\bT,\ph)$ and recall the notation $\mV_{\bbM,\e}$, $\tmV_{\bbM,\e}$, $\mG_{\bbM,\e}$, $\tmG_{\bbM,\e}$, $\bV_{\bbM,\e}$, $\tbV_{\bbM,\e}$, $\bG_{\bbM,\e}$ and $\tbG_{\bbM,\e}$ defined in Section \ref{sec:matrep_LAi} with $\bM=\bT$. We see $\tbG_{\bbM,\e}(M)=1$ for all $M\in \bT_{0}$ by this definition. We also notice that the matrices $\mV_{\bbM,\e}$ and $\mG_{\bbM,\e}$ are both irreducible nonnegative matrices. We take the left Perron eigenvector of  $\mV_{\bbM,\e}$ and the left Perron eigenvector of $\tmV_{\bbM,\e}$ which are denoted by
\alil{
\bU_{\bbM,\e}=(\bU_{\bbM,\e}(M)),\quad \tbU_{\bbM,\e}=(\tbU_{\bbM,\e}(M))\label{eq:bU=}
}
respectively, satisfying $\bU_{\bbM,\e}\cdot \bV_{\bbM,\e}=\tbU_{\bbM,\e}\cdot \tbV_{\bbM,\e}=1$.  Here the dot $\cdot$ means the inner product.

Assume also $\sharp\bT_{0}=2$. By Proposition \ref{prop:cm=cb}(2)(3), we have the equations
\begin{align*}
\frac{\tbU_{\bbM,\e}(M_{1})}{\tbU_{\bbM,\e}(M_{2})}&=\frac{(\tmV_{\bbM,\e})_{\lam(\e)}(M_{1}M_{1})}{(\tmV_{\bbM,\e})_{\lam(\e)}(M_{2}M_{1})}=\frac{\lam(\e)-\lam(\bM_{M_{2}},\e)}{\tmV_{\bbM,\e}(M_{1}M_{2})},\\
1=\frac{\tbG_{\bbM,\e}(M_{1})}{\tbG_{\bbM,\e}(M_{2})}&=\frac{(\tmG_{\bbM,\e})_{\lam(\e)}(M_{1}M_{2})}{(\tmG_{\bbM,\e})_{\lam(\e)}(M_{2}M_{1})}=\frac{\lam(\e)-\lam(\bM_{M_{2}},\e)}{\tmG_{\bbM,\e}(M_{1}M_{2})}.
\end{align*}
Therefore, $\tbU_{\bbM,\e}(M_{1})/\tbU_{\bbM,\e}(M_{2})=\tmG_{\bbM,\e}(M_{1}M_{2})/\tmV_{\bbM,\e}(M_{1}M_{2})\to 1$ follows from Corollary \ref{cor:geneform_prop2}. Since $\sum_{k=1}^{2}\tbV_{\bbM,\e}(M_{k})\tbU_{\bbM,\e}(M_{k})=1$ holds,
\begin{align*}
\frac{1}{\tbU_{\bbM,\e}(M_{2})}=\tbV_{\bbM,\e}(M_{1})\frac{\tbU_{\bbM,\e}(M_{1})}{\tbU_{\bbM,\e}(M_{2})}+\tbV_{\bbM,\e}(M_{2})\to m(\Sigma_{1})+m(\Sigma_{2})=1
\end{align*}
is satisfied under the assumption $\mu(\e,\cd)\to m$. Thus $\tbU_{\bbM,\e}(M_{k})\to 1$ for $k=1,2$.
On the other hand,  Remark \ref{rem:cb=}(1) yields the form
\ali
{
\tbV_{\bbM,\e}(M_{k})\tbU_{\bbM,\e}(M_{k})=\frac{\lam(\e)-\lam(\bM_{M_{k^\p}},\e)}{\sum_{l=1}^{2}(\lam(\e)-\lam(\bM_{M_{l}},\e))}
}
for $\{k,k^\p\}=\{1,2\}$. Theorem \ref{th:conv_gibbs_iff_mVM} implies that the left hand side converges to a number for each $k=1,2$ if and only if $\mu(\e,\cd)$ converges to a measure. Hence the assertion is complete.
\qed
\subsection{Proof of Theorem \ref{th:m0=3}}
To prove this theorem, we need some lemmas as follows.
\lem
{\label{lem:m0=3_tU->1}
Assume that the conditions $(\Si.1)$-$(\Si.3)$, $(\Ph.1)$-$(\Ph.3)$ and $\sharp \bT_{0}=3$ are satisfied. Then Then $\tbU_{\bbM,\e}(M)\to 1$ for each $M\in \bT_{0}$.
}
\pros
We write $\bT_{0}=\{M_{1},M_{2},M_{3}\}$.
In this case,
\ali
{
\frac{\tbU_{\bbM,\e}(M_{1})}{\tbU_{\bbM,\e}(M_{2})}&=\frac{(\tmV_{\bbM,\e})_{\lam(\e)}(M_{1}M_{3})}{(\tmV_{\bbM,\e})_{\lam(\e)}(M_{2}M_{3})}\\
&=\frac{\tmV_{\bbM,\e}(M_{2}M_{1})\tmV_{\bbM,\e}(M_{3}M_{2})+(\lam(\e)-\lam(\bM_{M_{2}},\e))\tmV_{\bbM,\e}(M_{3}M_{1})}{(\lam(\e)-\lam(\bM_{M_{1}},\e))\tmV_{\bbM,\e}(M_{3}M_{2})+\tmV_{\bbM,\e}(M_{1}M_{2})\tmV_{\bbM,\e}(M_{3}M_{1})}
}
and
\ali
{
1=\frac{\tbG_{\bbM,\e}(M_{1})}{\tbG_{\bbM,\e}(M_{2})}&=\frac{(\tmG_{\bbM,\e})_{\lam(\e)}(M_{1}M_{3})}{(\tmG_{\bbM,\e})_{\lam(\e)}(M_{2}M_{3})}\\
=&\frac{\tmG_{\bbM,\e}(M_{2}M_{1})\tmG_{\bbM,\e}(M_{3}M_{2})+(\lam(\e)-\lam(\bM_{M_{2}},\e))\tmG_{\bbM,\e}(M_{3}M_{1})}{(\lam(\e)-\lam(\bM_{M_{1}},\e))\tmG_{\bbM,\e}(M_{3}M_{2})+\tmG_{\bbM,\e}(M_{1}M_{2})\tmG_{\bbM,\e}(M_{3}M_{1})}
}
are satisfied from Proposition \ref{prop:cm=cb}(2)(3). By $\tmV_{\bbM,\e}(M_{k}M_{l})/\tmG_{\bbM,\e}(M_{k}M_{l})\to 1$ for each $k\neq l$, we obtain convergence $\tbU_{\bbM,\e}(M_{1})/\tbU_{\bbM,\e}(M_{2})\to 1$ by using Corollary \ref{cor:geneform_prop2} and Lemma \ref{prop:a+b/c+d}. By a similar argument above, $\tbU_{\bbM,\e}(M_{i})/\tbU_{\bbM,\e}(M_{j})\to 1$ is fulfilled for $(i,j)=(1,3),(2,3)$. By $1=\sum_{k=1}^{3}\tbV_{\bbM,\e}(M_{k})\tbU_{\bbM,\e}(M_{k})$, we have
\ali
{
\frac{1}{\tbU_{\bbM,\e}(M_{k})}&=\sum_{M\in \bT_{0}}\tbV_{\bbM,\e}(M)\frac{\tbU_{\bbM,\e}(M)}{\tbU_{\bbM,\e}(M_{k})}\to \sum_{M\in \bT_{0}}\mu(\SiM)=1
}
when $\mu(\e,\cd)$ converges to a measure $\mu$ weakly by using Theorem \ref{th:conv_gibbs_iff_mVM}. Hence the assertion follows.
\proe
Let $\bi\subset \bT_{0}$ be a non-empty subset with $\sharp \bi=2$. Denoted by $\lam^{v}(\bbM,\bi,\e)$ the Perron eigenvalue of the submatrix $\mV_{\bbM,\e}(\bi)$ of $\mV_{\bbM,\e}$, and by $\bV_{\bbM,\e}^{v}(\bi,\cd)=(\bV_{\bbM,\e}^{v}(\bi,M))_{M\in \bi}$ the corresponding right Perron eigenvector of $\mV_{\bbM,\e}(\bi)$ satisfying $\sum_{M\in \bi}\bV_{\bbM,\e}^{v}(\bi,M)=1$. Note that the Perron eigenvalue of $\tmV_{\bbM,\e}(\bi)$ coincides with $\lam^{v}(\bbM,\bi,\e)$.
We also take the corresponding right Perron eigenvector $\tbV_{\bbM,\e}^{v}(\bi,\cd)$ of $\tmV_{\bbM,\e}(\bi)$  with $\sum_{M\in \bi}\tbV_{\bbM,\e}^{v}(\bi,M)=1$. 
\lem
{
\label{lem:lam_v-l/lam-l->1_m0>=3}
Assume that the conditions $(\Si.1)$-$(\Si.3)$, $(\Ph.1)$-$(\Ph.3)$ and $\sharp \bT_{0}\geq 2$ are satisfied. Let $\bi=\{M_{1},M_{2}\}\subset \bT_{0}$ with $M_{1}\neq M_{2}$, $\bM_{\bi}=\bi\cup \bT_{1}$ and $\bbJ=\bbM(B,\bM_{\bi},\ph)=\{\bM_{M_{1}},\bM_{M_{2}}\}$.Then we have
\ite{
\item $\tmV_{\bbM,\e}(M M^\p)>0$ if and only if $\tmV_{\bbJ,\e}(M M^\p)>0$ for $M,M^\p\in \bi$ with $M\neq M^\p$. In these cases, $\tmV_{\bbM,\e}(M M^\p)/\tmV_{\bbJ,\e}(M M^\p)\to 1$.
\item If $\tmV_{\bbM,\e}(M_{1} M_{2})>0$ and $\tmV_{\bbM,\e}(M_{2}M_{1})>0$, then we have convergence $(\lam^{v}(\bbM,\bi,\e)-\lam(\bM_{M},\e))/(\lam(\bbJ,\e)-\lam(\bM_{M},\e))\to 1$ for $M\in \bi$.
\item If $\tmV_{\bbM,\e}(M M^\p)>0$ then $\tbV^{v}_{\bbM,\e}(\bi,M)/\tbV_{\bbJ,\e}(M)\to 1$ for $\{M,M^\p\}=\bi$.
}
}
\pros
\item[(1)] By the proof of Proposition \ref{prop:geneform_prop2}, we see $\mV_{\bbM,\e}(M^\pp M^\ppp)\asymp \te(M^\pp,M^\ppp:T_{1})\asymp \mV_{\bbJ,\e}(M^\pp M^\ppp)$. In particular, the assertion in the first half is valid.
By virtue of Proposition \ref{prop:geneform_prop2}(1) again,
\ali
{
&\frac{\mV_{\bbM,\e}(MM^\p)}{\te(M,M^\p:\bT_{1})},\ \frac{\mV_{\bbJ,\e}(MM^\p)}{\te(M,M^\p:\bT_{1})}\\
\to& \sum_{W\in \SP(M,M^\p:\bT_{1})}c_{MM^\p}(W)\nu(\bM_{M^\p},\LRLR(W)h(\bM_{M},\cd))
}
are satisfied as $\e \to 0$ running through some subsequence of any sequence $(\e_{n})$.
Arbitrariness of $(\e_{n})$ yields
\ali{ \tmV_{\bbM,\e}(MM^\p)/\tmV_{\bbJ,\e}(MM^\p)=\mV_{\bbM,\e}(MM^\p)/\mV_{\bbJ,\e}(MM^\p)\to 1
}
as $\e\to 0$.
\item[(2)] Denoted by
\ali{
\tmV_{\bbM,\e}(\bi)=\matII{
\eta_{1}(\e) &a_{12}(\e) \\
a_{21}(\e) &\eta_{2}(\e) \\
} \text{ and } \tmV_{\bbJ,\e}=\matII{
\eta_{1}(\e) &b_{12}(\e) \\
b_{21}(\e) &\eta_{2}(\e) \\
}
}
with $\eta_{i}(\e)=\lam(\bM_{M_{i}},\e)$ for $i=1,2$, and $a_{ij}(\e)=\tmV_{\bbM,\e}(M_{i}M_{j})$ and $b_{ij}(\e)=\tmV_{\bbJ,\e}(M_{i}M_{j})$ for $i\neq j$. Then
$a_{12}(\e)/b_{12}(\e)\to 1$ and $a_{21}(\e)/b_{21}(\e)\to 1$ by (1). By the form
\ali
{
&\lam^{v}(\bbM,\bi,\e)-\min(\eta_{1}(\e),\eta_{2}(\e))\\
=&\left(|\eta_{1}(\e)-\eta_{2}(\e)|+\sqrt{(\eta_{1}(\e)-\eta_{2}(\e))^{2}+4a_{12}(\e)a_{21}(\e)}\right)/2\\
&\lam(\bbJ,\e)-\min(\eta_{1}(\e),\eta_{2}(\e))\\
=&\left(|\eta_{1}(\e)-\eta_{2}(\e)|+\sqrt{(\eta_{1}(\e)-\eta_{2}(\e))^{2}+4b_{12}(\e)b_{21}(\e)}\right)/2.
}
and by Lemma \ref{prop:a+b/c+d}, we obtain
\ali
{
\frac{\lam^{v}(\bbM,\bi,\e)-\min(\eta_{1}(\e),\eta_{2}(\e))}{\lam(\bbJ,\e)-\min(\eta_{1}(\e),\eta_{2}(\e))}\to 1.
}
On the other hand, we see
\ali
{
\frac{a_{12}(\e)a_{21}(\e)}{b_{12}(\e)b_{21}(\e)}=\frac{(\lam^{v}(\bbM,\bi,\e)-\min(\eta_{1}(\e),\eta_{2}(\e)))(\lam^{v}(\bbM,\bi,\e)-\max(\eta_{1}(\e),\eta_{2}(\e)))}{(\lam(\bbJ,\e)-\min(\eta_{1}(\e),\eta_{2}(\e)))(\lam(\bbJ,\e)-\max(\eta_{1}(\e),\eta_{2}(\e)))},
}
Thus
\ali
{
\frac{\lam^{v}(\bbM,\bi,\e)-\max(\eta_{1}(\e),\eta_{2}(\e))}{\lam(\bbJ,\e)-\max(\eta_{1}(\e),\eta_{2}(\e))}\to 1.
}
Since $\eta_{1}(\e),\eta_{2}(\e)$ are in $\{\min(\eta_{1}(\e),\eta_{2}(\e)),\max(\eta_{1}(\e),\eta_{2}(\e))\}$, the assertion is fulfilled.
\item[(3)] We have the form
\ali
{
\left(\begin{array}{c}
\tbV_{\bbJ,\e}(M^\pp)\\
\tbV_{\bbJ,\e}(M^\ppp)
\end{array}\right)=
&\frac{c(\e)}{c_{1}(\e)}
\left(\begin{array}{c}
\lam(\bbJ,\e)-\lam(\bM_{M^\pp},\e)\\
\tmV_{\bbJ,\e}(M^\ppp M^\pp)
\end{array}\right)\\
\left(\begin{array}{c}
\tbV_{\bbM,\e}^{v}(\bi,M^\pp)\\\tbV_{\bbM,\e}^{v}(\bi,M^\ppp)
\end{array}\right)
=&\frac{1}{c_{2}(\e)}
\left(\begin{array}{c}
\lam(\bbM,\bi,\e)-\lam(\bM_{M^\pp},\e)\\\tmV_{\bbM,\e}(M^\ppp M^\pp)
\end{array}\right)
}
with $c(\e)=\tbV_{\bbJ,\e}(M_{1})+\tbV_{\bbJ,\e}(M_{2})$, $c_{1}(\e)=\lam(\bbJ,\e)-\lam(\bM_{M^\pp},\e)+\tmV_{\bbJ,\e}(M^\ppp M^\pp)$ and $c_{2}(\e)=\lam(\bbM,\bi,\e)-\lam(\bM_{M^\pp},\e)+\tmV_{\bbM,\e}(M^\ppp M^\pp)$. We have $c(\e)\to 1$ from Theorem \ref{th:conv_gibbs_iff_mVM} by replacing $\bbM=\bbJ$. 
The assumption $\tmV_{\bbM,\e}(M^\pp M^\ppp)>0$ implies 
$\tbV_{\bbJ,\e}(M^\pp)>0$ and $\tbV_{\bbM,\e}^{v}(\bi,M^\pp)>0$. Hence the assertion follows from (1), (2) and Proposition \ref{prop:a+b/c+d}.
\proe
\lem
{\label{lem:le-lvie/le-lje->1_m0=3}
Assume that the conditions $(\Si.1)$-$(\Si.3)$, $(\Ph.1)$-$(\Ph.3)$ and $\sharp \bT_{0}=3$ are satisfied. We use the notation $\bi,\bM_{\bi},\bbJ, \lam(\bbJ.\e)$ in Lemma \ref{lem:lam_v-l/lam-l->1_m0>=3}. Then we have
\ali
{
\frac{\lam(\e)-\lam(\bbJ,\e)}{\lam(\e)-\lam^{v}(\bbM,\bi,\e)}\to 1.
}
}
\pros
$\bT_{0}$ denotes $\{M,M^\p,M^\pp\}$.
In the case when either $\tmV_{\bbM,\e}(MM^\p)=0$ or $\tmV_{\bbM,\e}(M^\p M)=0$, we notice
\ali
{
\lam(\bbJ,\e)=\lam^{v}(\bbM,\bi,\e)=\max(\lam(\bM_{M},\e),\lam(\bM_{M^\p},\e)).
}
Thus we obtain
$(\lam(\e)-\lam(\bbJ,\e))/(\lam(\e)-\lam^{v}(\bbM,\bi,\e))=1$. 

Assume $\tmV_{\bbM,\e}(M M^\p)>0$ and $\tmV_{\bbM,\e}(M^\p M)>0$. In this case, we notice that the $\bM=\bM_{\bi}$ satisfies the condition $(\Si.\bM)$. First we estimate the difference $\lam(\e)-\lam^{v}(\bbM,\bi,\e)$.
We have
\ali
{
&\lam(\e)-\lam^{v}(\bbM,\bi,\e)=\frac{\tbU_{\bbM,\e}(\tmV_{\bbM,\e}-\tmV_{\bbM,\e}(\bi))\tbV_{\bbM,\e}^{v}(\bi,\cd)}{\tbU_{\bbM,\e}\cd\tbV_{\bbM,\e}^{v}(\bi,\cd)}\\
=&\frac{\tbU_{\bbM,\e}(M^\pp)(\tmV_{\bbM,\e}(M^\pp M)\tbV^{v}_{\bbM,\e}(M)+\tmV_{\bbM,\e}(M^\pp M^\p)\tbV^{v}_{\bbM,\e}(M^\p))}{\tbU_{\bbM,\e}(M)\tbV_{\bbM,\e}^{v}(\bi,M)+\tbU_{\bbM,\e}(M^\p)\tbV_{\bbM,\e}^{v}(\bi,M^\p)}\\
=&I_{1}(\e)(\tmV_{\bbM,\e}(M^\pp M)I_{2}(\e)\tbV_{\bbJ,\e}(M)+\tmV_{\bbM,\e}(M^\pp M^\p)I_{3}(\e)\tbV_{\bbJ,\e}(M^\p))\\
=&I_{1}(\e)(\frac{\nu(\bM_{M^\pp},\e,g(\e,\cd))}{\nu(\bM_{M},\e,g(\e,\cd))}\mV_{\bbM,\e}(M^\pp M)I_{2}(\e)\frac{\nu(\bM_{M},\e,g(\e,\cd))}{\nu(\bM_{\bi},\e,g(\e,\cd))}\bV_{\bbJ,\e}(M)\\
&+\frac{\nu(\bM_{M^\pp},\e,g(\e,\cd))}{\nu(\bM_{M^\p},\e,g(\e,\cd))}\mV_{\bbM,\e}(M^\pp M^\p)I_{3}(\e)\frac{\nu(\bM_{M^\p},\e,g(\e,\cd))}{\nu(\bM_{\bi},\e,g(\e,\cd))}\bV_{\bbJ,\e}(M^\p))\\
=&I_{1}(\e)\left(\mV_{\bbM,\e}(M^\pp M)I_{2}(\e)\bV_{\bbJ,\e}(M)+\mV_{\bbM,\e}(M^\pp M^\p)I_{3}(\e)\bV_{\bbJ,\e}(M^\p)\right)\times\\
&\quad\frac{\nu(\bM_{M^\pp},\e,g(\e,\cd))}{\nu(\bM_{\bi},\e,g(\e,\cd))}
}
and
\ali
{
I_{1}(\e)&=\frac{1}{\frac{\tbU_{\bbM,\e}(M)}{\tbU_{\bbM,\e}(M^\pp)}\tbV_{\bbM,\e}^{v}(\bi,M)+\frac{\tbU_{\bbM,\e}(M^\p)}{\tbU_{\bbM,\e}(M^\pp)}\tbV_{\bbM,\e}^{v}(\bi,M^\p)}\to 1\\
I_{2}(\e)&=\frac{\tbV^{v}_{\bbM,\e}(\bi,M)}{\tbV_{\bbJ,\e}(M)}\to 1,\quad I_{3}(\e)=\frac{\tbV^{v}_{\bbM,\e}(\bi,M^\p)}{\tbV_{\bbJ,\e}(M^\p)}\to 1\\
}
by Lemma \ref{lem:m0=3_tU->1}, Lemma \ref{lem:lam_v-l/lam-l->1_m0>=3}(3) and the fact $\tbV_{\bbM,\e}^{v}(\bi,M)+\tbV_{\bbM,\e}^{v}(\bi,M^\p)=1$.
Let 
\ali
{\te(\bi)=\sum_{W\in \bi}\te(M^\pp ,W:\bT_{1})\bV_{\bbJ,\e}(W).
}
We take a sequence $(\e_{n})$ so that for any $W=W_{1}\cdots W_{|W|}\in \bigcup_{W^\p\in \bi}\SP(M^\pp,W^\p:\bT_{1})$,
$\te(W)\bV_{\bbJ,\e}(W_{|W|})/\te(\bi)$ converges to a number $c(W)$ when $\e=\e_{n}\to 0$. Thus, 
\alil
{
&\frac{\mV_{\bbM,\e}(M^\pp M)I_{2}(\e)\bV_{\bbJ,\e}(M)+\mV_{\bbM,\e}(M^\pp M^\p)I_{3}(\e)\bV_{\bbJ,\e}(M^\p)}{\te(\bi)}\nonumber\\
\to &\sum_{W^\p\in \bi}\sum_{W\in \SP(M^\pp,W^\p:\bT_{1})}c(W)\nu(\bM_{W^\p},\LRLR(W)h(\bM_{M^\pp},\cd))\label{eq:le-lvie/le-lje->1_m0=3_1}
}
is satisfied as $\e \to 0$ running through $(\e_{n})$.

Next, we will evaluate the difference $\lam(\e)-\lam(\bbJ,\e)$. By using the equations $\LR_{A,\tPhe}1=\lam(\e)$ and $\LR_{A(\bM_{\bi},\tPhe}^{*}\ti{\nu}(\bM_{\bi},\e,\cd)=\lam(\bbJ,\e)\ti{\nu}(\bM_{\bi},\e,\cd)$, 
we find
\ali
{
\lam(\e)-\lam(\bbJ,\e)=&\ti{\nu}(\bM_{\bi},\e,\LR_{A-A(\bM_{\bi}),\tPhe}1)
=\ti{\nu}(\bM_{\bi},\e,\sum_{W\in \bM_{\bi}}\LR_{A_{M^\pp W},\tPhe}1)\\
=&\frac{\nu(\bM_{\bi},\e,\sum_{W\in \bM_{\bi}}\LR_{A_{M^\pp W},\Phe}g(\e,\cd))}{\nu(\bM_{\bi},\e,g(\e,\cd))}\\
=&\nu(\bM_{\bi},\e,\sum_{W\in \bM_{\bi}}\LR_{A_{M^\pp W},\Phe}f(M^\pp,\e,\cd))\frac{\nu(\bM_{M^\pp},\e,g(\e,\cd))}{\nu(\bM_{\bi},\e,g(\e,\cd))}
}
with
\ali
{ f(M^\pp,\e,\cd)&=\frac{g(\e,\cd)\chi_{\SiMpp}}{\nu(\bM_{M^\pp},\e,g(\e,\cd))}\to \frac{h(M^\pp,\cd)}{\nu(\bM_{M^\pp},\e,h(M^\pp,\cd))}=
h(\bM_{M^\pp},\cd)
}
from $\nu(\bM_{M^\pp},\e,g(\e,\cd)\chi_{\Si_{W}})/\bG_{\e}(M^\pp)\to 0$ for $W\in \bT_{1}$ by using (\ref{eq:lem:thMM->_1}) and Proposition \ref{prop:tsc_LAMPhe_T0=1}(2). Here
we have
\ali
{
&\frac{\nu(\bM_{\bi},\e,\LR_{A_{M^\pp W},\Phe}f(M^\pp,\e,\cd))}{\te(\bi)}\\
=&\frac{\bV_{\bbJ,\e}(W)\te(M^\pp W)}{\te(\bi)}\frac{\nu(\bM_{\bi},\e,\LR_{A_{M^\pp W},\Phe}f(M^\pp,\e,\cd))/\nu(\bM_{\bi},\e,\chi_{\Si_{W}})}{\nu(\bM_{\bi},\e,h(\bM_{W},\e,\cd))\te(M^\pp W)/\nu(\bM_{\bi},\e,\chi_{\Si_{W}})}\\
\to &c(M^\pp W)\frac{\nu(W,\LR(M^\pp W) h(\bM_{M^\pp},\cd))}{\nu(W,h(\bM_{W},\cd))}\\
=&c(M^\pp W)\nu(\bM_{W},\LR(M^\pp W)h(\bM_{M^\pp},\cd))
}
as $\e\to 0$ running through $(\e_{n})$ by Proposition \ref{prop:conv_nor_TSC_M}(2). Furthermore, the proof of Proposition \ref{prop:speed_tscM}(2) implies that
for $W\in \bT_{1}$
\ali
{
&\nu(\bM_{\bi},\e,\LR_{A_{M^\pp W},\Phe}f(M^\pp,\e,\cd))/\te(\bi)\\
=&\sum_{W^\p\in \bi}\sum_{W^\pp\in \SP(W,W^\p:\bT_{1})}\nu(\bM_{\bi},\e,\LRR_{\e}(W^\pp)\LR_{A_{M^\pp W},\Phe}f(M^\pp,\e,\cd))/\te(\bi)\\
&+o(\VTe(\bM_{\bi},W))/\te(\bi)\\
=&\sum_{W^\p\in \bi}\sum_{W^\pp\in \SP(W,W^\p:\bT_{1})}\bV_{\bbJ,\e}(W)\te(M^\pp W^\pp)/\te(\bi)\times\\
&\nu(\bM_{\bi},\e,\LRR_{\e}(W^\pp)\LR_{A_{M^\pp W},\Phe}f(M^\pp,\e,\cd))/(\bV_{\bbJ,\e}(W)\te(M^\pp W^\pp))+o(1)\\
\to&\sum_{W^\p\in \bi}\sum_{W^\pp\in \SP(W,W^\p:\bT_{1})}c(M^\pp W^\pp)\nu(\bM_{W^\p},\LRLR(M^\pp W^\pp)h(\bM_{M^\pp},\cd))
}
by using the limit
\ali
{
\frac{\nu(\bM_{\bi},\e,f\chi_{\Si_{W^\p}})}{\bV_{\bbJ,\e}(W^\p)}=\frac{\nu(\bM_{\bi},\e,f\chi_{\Si_{W^\p}})}{\nu(\bM_{\bi},\e,h(\bM_{W^\p},\e,\cd))}\to \frac{\nu(W^\p,f)}{\nu(W^\p,h(\bM_{W^\p},\cd))}=\nu(\bM_{W^\p},f)
}
with Proposition \ref{prop:conv_nor_TSC_M}(2), where $\LRLR$ is defined in Proposition \ref{prop:geneform_prop2}. Thus we have
\alil
{
&\frac{\nu(\bM_{\bi},\e,\sum_{W\in \bM_{\bi}}\LR_{A_{M^\pp W},\Phe}f(M^\pp,\e,\cd))}{\te(\bi)}\nonumber\\
\to& \sum_{W^\p\in \bi}\sum_{W\in \SP(M^\pp,W^\p:\bT_{1})}c(W)\nu(\bM_{W^\p},\LRLR(W)h(\bM_{M^\pp},\cd))\label{eq:le-lvie/le-lje->1_m0=3_2}
}
as $\e\to 0$ running over $(\e_{n})$. Consequently, the limits (\ref{eq:le-lvie/le-lje->1_m0=3_1}) and (\ref{eq:le-lvie/le-lje->1_m0=3_2}) yield
\ali
{
&\frac{\lam(\e)-\lam(\bbJ,\e)}{\lam(\e)-\lam^{v}(\bbM,\bi,\e)}\\
=&\frac{I_{1}(\e)(\mV_{\bbM,\e}(M^\pp M)I_{2}(\e)\bV_{\bbJ,\e}(M)+\mV_{\bbM,\e}(M^\pp M^\p)I_{3}(\e)\bV_{\bbJ,\e}(M^\p))/\te(\bi)}{\nu(\bM_{\bi},\e,\sum_{W\in \bM_{\bi}}\LR_{A_{M^\pp W},\Phe}f(M^\pp,\e,\cd))/\te(\bi)}\to 1
}
as $\e\to 0$ running over $(\e_{n})$. This converging does not depend on how to choose a sequence $(\e_{n})$.
Hence the proof is complete.
\proe
\ \\(Proof of Theorem \ref{th:m0=3})\\
Recall the form
\ali
{
&\tbV_{\e}(M)\tbU_{\e}(M)\ti{c}(\e)=(\lam(\e)-\lam^{v}(\bbM,\{M^\p,M^\pp\},\e))\times\\
&\qquad\qquad\qquad(\lam(\e)-\lam(\bM_{M^\p},\e)+\lam^{v}(\bbM, \{M^\p,M^\pp\},\e)-\lam(\bM_{M^\pp},\e))
}
by Remark \ref{rem:cb=}(2), where $\ti{c}(\e)$ is a normalizing constant. By virtue of Lemma \ref{lem:m0=3_tU->1}, Lemma \ref{lem:lam_v-l/lam-l->1_m0>=3}(2) and Lemma \ref{lem:le-lvie/le-lje->1_m0=3} in addition to Lemma \ref{prop:a+b/c+d}, we obtain $\tbV_{\e}(M)/\delta_{\e}(M)\to 1$.
By using Theorem \ref{th:conv_gibbs_iff_mVM}, the proof is complete.\qed
\section{Examples}\label{sec:ex}
\subsection{A relation between potentials and convergence of Gibbs measures}\label{sec:ex_asymp}
In Theorem \ref{th:m0=2} and Theorem \ref{th:m0=3}, we found the relation between convergence of the Gibbs measure $\mu(\e,\cd)$ and convergence of the expression $\delta_{\e}(k)$ using eigenvalues of generalized Ruelle operators. However, the relation between convergence of $\mu(\e,\cd)$ and the potential $\Phe$ is not immediately clear. In this section, we illustrate this relation by using asymptotic expansion techniques for eigenvalues of Ruelle operators under the case $\sharp \bT_{0}=2$.

Assume that $(\Sigma.1)$-$(\Sigma.3)$, $(\varPhi.1)$-$(\varPhi.3)$ and $\sharp \bT_{0}=2$ are satisfied. We use the notation $\lam(\e)$ and $\lam(\bM,\e)$ defined in Theorem \ref{th:m0=2}.
Let
\ali
{ c_{1}(\e)=&\frac{\lam(\e)-\lam(\bM(1),\e)}{\lam(\e)-\lam(\bM(2),\e)},\\ c_{2}(\e)=&\frac{|\lam(\bM(2),\e)-\lam(\bM(1),\e)|(\lam(\bM(2),\e)-\lam(\bM(1),\e))}{\mV_{\bbM,\e}(M_{1}M_{2})\mV_{\bbM,\e}(M_{2}M_{1})}.
}
By direct calculation, we have the expression
\ali
{
\frac{|1-c_{1}(\e)|(1-c_{1}(\e))}{c_{1}(\e)}=\frac{|\lam(\bM(2),\e)-\lam(\bM(1),\e)|(\lam(\bM(2),\e)-\lam(\bM(1),\e))}{(\lam(\e)-\lam(\bM(1),\e))(\lam(\e)-\lam(\bM(2),\e))}.
}
By expanding the equation $\det(\lam(\e)I-\mV_{\bbM,\e})=0$, 
we see the relation
\alil{
\frac{|1-c_{1}(\e)|(1-c_{1}(\e))}{c_{1}(\e)}=c_{2}(\e).
}
Thus Theorem \ref{th:m0=2} implies that $\mu(\e,\cd)$ converges iff $c_{1}(\e)$ converges in $\R\cup\{\infty\}$ iff $c_{2}(\e)$ converges in $\R\cup\{\pm\infty\}$. If $c_{2}(\e)\to c_{2}$ then $\mu(\e,\cd)$ converges to
\ali
{
\case
{
\frac{1+\mathrm{sign}(c_{2})\sqrt{|c_{2}|/(|c_{2}|+4)}}{2}\mu(M_{1},\cd)+\frac{1-\mathrm{sign}(c_{2})\sqrt{|c_{2}|/(|c_{2}|+4)}}{2}\mu(M_{2},\cd),& c_{2}\in \R\\
\mu(M_{1},\cd),& c_{2}=\infty\\
\mu(M_{2},\cd),& c_{2}=-\infty.\\
}
}
\par
For simplicity we consider only the case:
\ite
{
\item[(a)] $\bT=\bT_{0}$ and $A_{MM}=M$ for each $M\in \bT_{0}$.
}
Under this condition, $\lam(\bM_{M},\e)$ equals $\lam(M,\e)$ and is the Perron eigenvalue of the operator $\LR_{A_{MM},\Phe}=\LR_{M,\phe}$ for $M\in \bT_{0}$.
Moreover, for fixed integers $n(1)\geq 0$ and $n(2)\geq 1$, we introduce the following conditions of asymptotic expansions:
\ite
{
\item[(b)] $\phe\in F_{\theta}(\SiAp,\R)$ has the $n(1)$-order asymptotic expansion $\phe=\ph+\ph_{1}\e+\cdots+\ph_{n(1)}\e^{n(1)}+\ti{\ph}_{n(1)}(\e,\cd)\e^{n(1)}$ and $\|\ti{\ph}_{n(1)}(\e,\cd)\|_{\infty}\to 0$ as $\e\to 0$, where $\ph,\ph_{1},\dots, \ph_{n(1)}$ are in $F_{\theta}(\SiAp,\R)$.
\item[(c)] $\psi(\e,\cd)\in F_{\theta}(\SiAp,\R)$ has the $n(2)$-order asymptotic expansion $e^{\psi(\e,\cd)}\chi_{N}=\psi_{1}\e+\cdots+\psi_{n(2)}\e^{n(2)}+\ti{\psi}_{n(2)}(\e,\cd)\e^{n(2)}$ and $\|\ti{\psi}_{n(2)}(\e,\cd)\|_{\infty}\to 0$ as $\e\to 0$, where $\psi_{1},\dots,\psi_{n(2)}$ are in $F_{\theta}(\SiAp,\R)$.
\item[(d)] $\sup_{\e>0}[[\ti{\ph}_{n(1)}(\e,\cd)]]_{\theta}<\infty$ and $\sup_{\e> 0}[[\ti{\psi}_{n(2)}(\e,\cd)]]_{\theta}<\infty$.
}
We notice that the assumptions $(b)-(d)$ yield the conditions $(\Ph.1)$-$(\Ph.3)$. The condition (b) implies the following:
\prop{[\cite{T2011}]\label{prop:asymp_lam}
Assume that $(\Sigma.1)$-$(\Sigma.3)$, (a) and (b). Then for each $M\in \bT_{0}$, there exist unique numbers $\lam_{M,j}\in \R$ $(1\leq j\leq n)$ such that the asymptotic expansion
\ali
{
\lam(\bM_{M},\e)=\lam+\lam_{M,1}\e+\cdots+\lam_{M,n}\e^{n}+\ti{\lam}_{M,n}(\e)\e^{n}
}
is satisfied with $|\ti{\lam}_{M,n}(\e)|\to 0$ as $\e\to 0$.
}
\pros
By applying Lemma 4.1 in \cite{T2011} to the operator $\LR_{M,\phe_{M}}\in \LR(C(\Si_{M}^{+}))$, this assertion follows immediately.
\proe
For a convenience, we put $\psi_{0}=\chi_{\SiAp\setminus N}$ and write $\bT_{0}=\{1,2\}$ by regarding as $k=B_{kk}$ for $k=1,2$. We introduce some numbers below.
\ali{
s&=
\case
{
l,&\lam_{1,1}=\lam_{2,1},\dots,\lam_{1,l-1}=\lam_{2,l-1} \text{ and }\\
&\qquad\lam_{1,l}\neq \lam_{2,l} \text{ for some }1\leq l\leq n(1) \\
\infty,&\lam_{1,l}= \lam_{2,l} \text{ for any }1\leq l\leq n(1)
}\\
s(kk^\p)&=
\case
{
l,&\psi_{0}=\cdots=\psi_{l-1}=0 \text{ on } \Si_{kk^\p} \text{ and }\\
&\qquad \psi_{l}\neq 0 \text{ on } \Si_{kk^\p} \text{ for some }0\leq l\leq n(2) \\
\infty,&\psi_{l}= 0 \text{ on } \Si_{kk^\p} \text{ for any }0\leq l\leq n(2)
}\\
d(kk^\p)&=\nu(k^\p,\LR_{B_{kk^\p},\ph}(h(k,\cd)\psi_{s(kk^\p)})) \text{ if }s(kk^\p)<\infty.
}
for $k,k^\p\in \bT_{0}$ with $k\neq k^\p$. Under these notation, we obtain the following:
\prop{
\label{prop:c1toc2}
Assume that $(\Sigma.1)$-$(\Sigma.3)$, $(a)$-$(d)$ and $\sharp \bT_{0}=2$ are satisfied. Then $c_{2}(\e)$ converges to
\alil
{
c_{2}=\case
{
\displaystyle \frac{|\lam_{1,s}-\lam_{2,s}|(\lam_{2,s}-\lam_{2,s})}{d(12)d(21)},&\text{ if }2s=s(12)+s(21)<\infty\\
0,&\text{ if }2s>s(12)+s(21)\\
\mathrm{sign}(\lam_{1,s}-\lam_{2,s})\infty,&\text{ if }2s<\min(n(2)+1,s(12))\\
&\qquad\qquad+\min(n(2)+1,s(21)).\\
}\label{eq:prop:c1toc2}
}
}
\pros
First we show that if $s(kk^\p)<\infty$ then $\mV_{\bbM,\e}(MM^\p)/\e^{s(kk^\p)}$ converges to $d(kk^\p)$. To do this, we will claim $e^{\Phe}/\e^{s(kk^\p)}\to e^{\ph}\psi_{s(kk^\p)}$ on $\Si_{kk^\p}$.
 Indeed, in the case when $s(kk^\p)=0$ i.e. $\psi_{0}=\chi_{\Si_{kk^\p}\setminus N}\neq 0$ on $\Si_{kk^\p}$, we have
\ali
{
e^{\Phe}=e^{\phe}\psi_{0}+e^{\phe}e^{\psi(\e,\cd)}\chi_{N}\to e^{\ph}\psi_{0} \text{ on } \Si_{kk^\p}.
}
In the case when $0<s(kk^\p)<\infty$, we see
\ali
{
e^{\Phe}/e^{s(kk^\p)}=&e^{\phe}e^{\psi(\e,\cd)}\chi_{N}/e^{s(kk^\p)}\to e^{\ph}\psi_{s(kk^\p)} \text{ on }\Si_{kk^\p}.
}
Therefore, we obtain
\ali
{
\frac{\mV_{\bbM,\e}(MM^\p)}{\e^{s(kk^\p)}}&=\frac{\nu(\bM,\e,\LR_{A_{MM^\p},\Phe}h(M,\e,\cd))}{\nu(\bM,\e,h(M^\p,\e,\cd))\e^{s(kk^\p)}}\\
&=\frac{\nu(\bM,\e,\LR_{A_{MM^\p},0}(h(M,\e,\cd)e^{\Phe}/\e^{s(kk^\p)})/\nu(\bM,\e,\Si_{M^\p})}{\nu(\bM,\e,h(M^\p,\e,\cd))/\nu(\bM,\e,\Si_{M^\p})}\to d(kk^\p).
}

Next we will show the assertion by considering the following cases:\\
The case $s<\infty$, $s(12)<\infty$ and $s(21)<\infty$: In this case, the expression
\ali
{
c_{2}(\e)=&\frac{|\lam_{1,s}-\lam_{2,s}|(\lam_{1,s}-\lam_{2,s})\e^{2s}+o(\e^{2s})}{d(12)d(21)\e^{s(12)+s(21)}+o(\e^{s(12)+s(21)})}.
}
is satisfied. This implies the assertion (\ref{eq:prop:c1toc2}).\\
The case $s=\infty$, $s(12)<\infty$ and $s(21)<\infty$: We see
\ali
{
c_{2}(\e)=&\frac{o(\e^{2n})}{d(12)d(21)\e^{s(12)+s(21)}+o(\e^{s(12)+s(21)})}\to 0.
}
The case $s<\infty$, $s(12)=\infty$ and $s(21)=\infty$: In this case, we have
\ali
{
c_{2}(\e)=&\frac{|\lam_{1,s}-\lam_{2,s}|(\lam_{1,s}-\lam_{2,s})\e^{2s}+o(\e^{2s})}{o(\e^{2n})}\to \mathrm{sign}(\lam_{1,s}-\lam_{2,s})\infty.
}
The case $s(12)=\infty$, $s(21)<\infty$ and $2s<n(2)+1+s(21)$: 
\ali
{
c_{2}(\e)=&\frac{|\lam_{1,s}-\lam_{2,s}|(\lam_{1,s}-\lam_{2,s})\e^{2s}+o(\e^{2s})}{o(\e^{n+s(21)})}\to \mathrm{sign}(\lam_{1,s}-\lam_{2,s})\infty.
}
The case $s(21)=\infty$, $s(12)<\infty$ and $2s<n(2)+1+s(12)$: By a similar argument in above case, we obtain the assertion.
\proe
If $s$, $s(12)$ and $s(21)$ do not satisfy in the conditions (\ref{eq:prop:c1toc2}), then $c_{2}(\e)$ might not converge. For example, we assume
\ali
{
\ti{\varphi}_{n(1)}(\e,\cd)&=\sin(1/\e)\chi_{\Si_{11}}\e^{n(1)+1}+o(\e^{n(1)+1}),\\
\ti{\psi}_{n(2)}(\e,\cd)&=\chi_{\Si_{12}\cup\Si_{21}}\e^{n(2)+1}+o(\e^{n(2)+1})
}
in $C(\SiAp)$. Choose any positive sequence $(\e(k))$ with $\inf_{k}\e(k)=0$ so that $\sin(1/\e(k))$ converges to a number $e(j)$ as $k\to \infty$.
Assume $s=s(12)=s(21)=\infty$ and $n=n(1)=n(2)$. Lemma 4.1 in \cite{T2011} again implies that
\alil
{
c_{2}(\e)\to c_{2}=\frac{|\lam_{1,n+1}-\lam_{2,n+1}|(\lam_{1,n+1}-\lam_{2,n+1})}{\nu(2,\LR_{B_{12},\ph}(h(1,\cd)))\nu(2,\LR_{B_{12},\ph}(h(1,\cd)))}\label{eq:ex_asymp2}
}
running through $(\e(k))$, where $\lam_{k,n+1}$ has the form 
\alil
{\lam_{k,n+1}=&\sum_{j=1}^{n+1}\nu(k,\mathcal{P}_{k,n+1-j}\LR_{k,j}h(k,\cd))),\label{eq:ex_asymp1}\\
\LR_{k,j}=&\LR_{B_{kk},\ph}(F_{j}\cdot),\qquad F_{j}=\sum_{0\leq i_{1},\cdots,i_{j}\leq j:\atop{i_{1}+2i_{2}+\cdots+j\cdot i_{j}=j}}\ph_{1}^{i_{1}}\cdots \ph_{j}^{i_{j}}/(i_{1}!\cdots i_{j}!)\nonumber
}
and $\mathcal{P}_{k,0}$ is the eigenprojection of the eigenvalue $\lam$ of $\LR_{B_{kk},\ph}$ defined by $\mathcal{P}_{k,0}f=\nu(k,h(k,\cd)f)$,  $\mathcal{P}_{k,j}=\sum_{i=1}^{j}\mathcal{P}_{k,j-i}(\lam_{k,j}\mathcal{I}-\LR_{k,j})\mathcal{S}_{k}$ and $\mathcal{S}_{k}=(\LR_{B_{kk},\ph}-\mathcal{P}_{k,0}-\lam\mathcal{I})^{-1}(\mathcal{I}-\mathcal{P}_{k,0})$. Notice that terms depending on $e(1)$ in the summation (\ref{eq:ex_asymp1}) are only $\nu(1,\mathcal{P}_{1,0}\LR_{1,n+1}h(1,\cd)))$, and this term is equal to $\lam e(1)$. Therefore the number $c_{2}$ given by (\ref{eq:ex_asymp2}) depends on $e(1)$. This implies that $c_{2}(\e)$ does not converge as $\e\to 0$ and hence so is for the Gibbs measure $\mu(\e,\cd)$.
\subsection{An example in the case when $\sharp \bT_{0}=3$}\label{sec:ex_m0=3}
In this section, we demonstrate limits of Gibbs measures by using Theorem \ref{th:m0=3} under the case $\sharp \bT_{0}=3$. Let $A$ and $B$ be $6\times 6$ matrices so that
\ali
{
A=\matVI
{
I&I&O\\
O&I&I\\
I&I&I\\
},\quad
B=\matVI
{
I&O&O\\
O&I&O\\
O&O&I\\
} \text{ with }
I=\matII
{
1&1\\
1&1\\
},\\
\varPhi(\e,\om)=\case
{
\log \e &\om\in \Si_{12}\cup \Si_{23}\\
s\log \e &\om\in \Si_{32}\\
(\sin(1/\e)/3+1)\log \e &\om\in \Si_{31}\\
0 &\text{otherwise},
}}
where $s$ is a positive number and $\Si_{kk^\p}=\{\om\in \SiAp\,:\,\om_{0}\om_{1}\in \{2k-1,2k\}\times \{2k^\p-1,2k^\p\}\}$. In this case, we see that the six conditions $(\Si.1)$-$(\Si.3)$, $(\Ph.1)$-$(\Ph.3)$ are satisfied, and $\sharp \bT_{0}=3$ holds. For a simple, we write $\bT_{0}=\{1,2,3\}$ by regarding as $k=B_{kk}$ for $k=1,2,3$. We find
\ali
{
\lam(1,\e)&=\lam(2,\e)=\lam(3,\e)=\lam(\{1,2\},\e)=\lam(\{1,3\},\e)=2 \text{ and}\\
\lam(\{2,3\},\e)&=2+2e^{(1+s)/2}.
}
Therefore, the expression (\ref{eq:dk_m=3}) have the forms
\ali
{
\delta^{0}_{\e}(1)&=(\lam(\e)-\lam(\{2,3\},\e))(\lam(\e)+\lam(\{2,3\},\e)-\lam(2,\e)-\lam(3,\e))\\
&=(\lam(\e)-2)^2-4e^{1+s}\\
\delta^{0}_{\e}(2)&=(\lam(\e)-\lam(\{1,3\},\e))(\lam(\e)+\lam(\{1,3\},\e)-\lam(1,\e)-\lam(3,\e))=(\lam(\e)-2)^{2}\\
\delta^{0}_{\e}(3)&=(\lam(\e)-\lam(\{1,2\},\e))(\lam(\e)+\lam(\{1,2\},\e)-\lam(1,\e)-\lam(2,\e))=(\lam(\e)-2)^{2}.
}
We see
\ali
{
\delta_{\e}(1)=c(\e)\left(1-\frac{4e^{1+s}}{(\lam(\e)-2)^2}\right),\quad \delta_{\e}(2)=\delta_{\e}(3)=c(\e)
}
with $c(\e)=(3-4e^{1+s}/(\lam(\e)-2)^2)^{-1}$. Moreover, 
\ali
{
\lam(\e)=\frac{1}{3}K(\e)^{1/3}+\frac{4\e^{s+1}}{K(\e)^{1/3}}+2
}
with $K(\e)=108\e^{2}\e^{(1/3)\sin(1/\e)+1}+12\sqrt{-12\e^{3(s+1)}+81\e^{4}\e^{(2/3)\sin(1/\e)+2}}$.
For any sequence $a(\e)$ on $\R$, denoted by $A(a(\e))$ the totally of the accumulation points of this squence as $\e\to 0$.
Now, we consider the five cases as follows:
\smallskip
\par
The case $0<s<7/9$ : We have $K(\e)=\e^{(3s+3)/2}(12\sqrt{-12}+o(1))$, $\lam(\e)=\e^{(s+1)/2}(2+o(1))+2$ and therefore $4\e^{s+1}/(\lam(\e)-2)^2\to 1$. Thus we see $\mu(\e,\cd)\to \sum_{k=2}^{3}\mu(B_{kk},\cd)/2$.
\smallskip
\par
The case $s=7/9$ : In this case, $K(\e)=c_{1}(\e)\e^{(3s+3)/2}$ with  $c_{1}(\e)=108\e^{-2/3}\e^{\sin(1/\e)/3+1}+12\sqrt{-12+81\e^{-4/3}\e^{2\sin(1/\e)/3+2}}$ and $4\e^{s+1}/(\lam(\e)-2)^2=4/(c_{1}(\e)^{1/3}/3+4/c_{1}(\e)^{1/3})^{2}$ are satisfied. Since $A(c_{1}(\e))=[0,c_{1}]$ with $c_{1}=108+12\sqrt{69}$ and the map $x\mapsto x^{1/3}/3+4/x^{1/3}$ has the lower bound $4/\sqrt{3}$ at $x=24\sqrt{3}<c_{1}$, we see $A(\mu(\e,\Si_{1}))=[1/9,1/3]$ and $A(\mu(\e,\Si_{k}))=[1/3,4/9]$ for $k=2,3$.
In particular, $\mu(\e,\cd)$ does not converge.
\smallskip
\par
The case $7/9<s<11/9$ : If $\sin(1/\e)=-1$ then we notice $K(\e)=\e^2\e^{\sin(1/\e)/3+1}(c_{1}+o(1))$,  $\lam(\e)=\e^{8/9}(c_{1}^{1/3}/3+o(1))+2$ and therefore we obtain $4\e^{s+1}/(\lam(\e)-2)^2\to 0$. On the other hand, if $\sin(1/\e)=1$ then $4\e^{s+1}/(\lam(\e)-2)^2\to 1$ as well as the case $(II)$. Thus $A(\mu(\e,\Si_{1}))=[0,1/3]$ and $A(\mu(\e,\Si_{k})=[1/3,1/2]$ for $k=2,3$. 
\smallskip
\par
The case $s=11/9$ : We have $K(\e)=c_{3}(\e)\e^2 \e^{\sin(1/\e)/3+1}$ with the map $c_{3}(\e)=108+12\sqrt{-12\e^{8/2}\e^{-2\sin(1/\e)/3-2}+81}$. Therefore $A(4\e^{s+1}/(\lam(\e)-2)^2)=[0,c_{2}]$, where $c_{2}=4/(c_{1}^{1/3}/3+4/c_{1}^{1/3})$. In particular, $\mu(\e,\cd)$ does not converge. In particular, $A(\mu(\e,\Si_{1}))=[(1-c_{2})/(3-c_{2}),1/3]$.
\smallskip
\par
The case $s>11/9$ : $K(\e)\in [\e^{10/3}(108+12\sqrt{-12\e^{8/3}+81}), \e^{8/3}(108+12\sqrt{-12\e^{10/3}+81})]$ is satisfied and thus $\lam(\e)\in [2^{2/3}\e^{10/9}, 3\e^{8/9}]$ for a small $\e>0$.  We obtain $e^{1+s}/(\lam(\e)-2)^2\to 0$ and hence $\mu(\e,\cd)\to \sum_{k=1}^{3}\mu(B_{kk},\cd)/3$.
\subsection{An example in the case when $\sharp \bT_{0}= 4$}\label{sec:m0>=4}
Assume $(\Si.1)$-$(\Si.3)$, $(\Ph.1)$-$(\Ph.3)$ and $\sharp \bT_{0}= 4$.
By Proposition \ref{prop:decom_eigen_Mmat}, we have the forms $\tbV_{\bbM,e}(k)\tbU_{\bbM,e}(k)=(\tmV_{\bbM,\e})_{\lam(\e)}(kk)/\sum_{l=1}^{4}(\tmV_{\bbM,\e})_{\lam(\e)}(ll)$ and
\ali{
&(\tmV_{\bbM,\e})_{\lam(\e)}(i_{1}i_{1})=\Big(\lam(\e)-\lam^{v}(\{i_{2},i_3,i_4\},\e)\Big)\times\\
&\qquad\Big((\lam(\e)-\lam^{v}(\{i_3,i_4\},\e))(\lam(\e)-\lam^{v}(i_3,\e)+\lam^{v}(\{i_3,i_4\},\e)-\lam^{v}(i_4,\e))\\
&\qquad+(\lam^{v}(\{i_2,i_3,i_4\},\e)-\lam^{v}(\{i_2,i_4\},\e))(\lam(\e)-\lam^{v}(i_4,\e))\\
&\qquad+(\lam(\e)-\lam^{v}(\{i_2,i_4\},\e))(\lam^{v}(\{i_2,i_4\},\e)-\lam^{v}(i_2,\e))\\
&\qquad+(\lam^{v}(\{i_2,i_3,i_4\},\e)-\lam^{v}(\{i_2,i_3\},\e))\times\\
&\qquad\quad(\lam^{v}(\{i_2,i_3,i_4\},\e)-\lam^{v}(i_2,\e)+\lam^{v}(\{i_2,i_3\},\e)-\lam^{v}(i_3,\e))\Big)
}
for $\{i_{1},i_{2},i_{3},i_{4}\}=\bT_{0}$.
If
\alil
{
\frac{\lam^{v}(\{j_{1},\dots,j_{k}\},\e)-\lam^{v}(\{j_{1}^{\p},\dots,j_{k^\p}^{\p}\},\e)}{\lam(\{j_{1},\dots,j_{k}\},\e)-\lam(\{j_{1}^{\p},\dots,j_{k^\p}^{\p}\},\e)}\to 1\label{eq:lav/lato1}
}
as $\e\to 0$ for all distinct elements $j_{1},\dots,j_{k}\in \bT_{0}$ and $j_{1}^{\p},\dots,j_{k^\p}^{\p}\in \{j_{1},\dots,j_{k}\}$ with $k^\p<k$, and if $\tbU_{\bbM,\e}(i)\to 1$ for each $i\in \bT_{0}$, then $\tbV_{\bbM,\e}(i)/\delta_{\e}(i)$ converges to $1$ for each $i\in \bT_{0}$, where $\delta_{\e}$ is defined by $\delta_{\e}(k)=\delta_{\e}^{0}(k)/\sum_{l=1}^{4}\delta_{\e}^{0}(l)$ and
\ali
{
&\delta_{\e}^{0}(i_1)=\Big(\lam(\e)-\lam(\{i_2,i_3,i_4\},\e)\Big)\times\\
\quad&\Big((\lam(\e)-\lam(\{i_3,i_4\},\e))(\lam(\e)-\lam(i_3,\e)+\lam(\{i_3,i_4\},\e)-\lam(i_4,\e))\\
\quad&+(\lam(\{i_2,i_3,i_4\},\e)-\lam(\{i_2,i_4\},\e))(\lam(\e)-\lam(i_4,\e))\\
\quad&+(\lam(\e)-\lam(\{i_2,i_4\},\e))(\lam(\{i_2,i_4\},\e)-\lam(i_2,\e))\\
\quad&+(\lam(\{i_2,i_3,i_4\},\e)-\lam(\{i_2,i_3\},\e))\times\\
&\quad(\lam(\{i_2,i_3,i_4\},\e)-\lam(i_2,\e)+\lam(\{i_2,i_3\},\e)-\lam(i_3,\e))\Big)
}
for $\{i_{1},i_{2},i_{3},i_{4}\}=\bT_{0}$. In fact, under the case when $\sharp \bT_{0}=3$,
Lemma \ref{lem:lam_v-l/lam-l->1_m0>=3}(2) and Lemma \ref{lem:le-lvie/le-lje->1_m0=3} provide that the ratio between the difference of the eigenvalues of generalized Ruelle operators and the difference of the eigenvalues of the submatrices of $\tmV_{\bbM,\e}$ converges to $1$ as $\e\to 0$. This fact plays an important role in the proof of Theorem \ref{th:m0=3}. When $\sharp \bT_{0}=4$, unfortunately (\ref{eq:lav/lato1}) are not valid in general.

For example, we take $8\times 8$ matrices $A$ and $B$, and a potential $\varPhi(\e,\cd)$ as follows:
\ali
{
A&=\matVIII
{
1&1&1&1&1&1&1&1\\
1&1&1&1&1&1&1&1\\
1&1&1&1&0&0&0&0\\
1&1&1&1&1&1&0&0\\
1&1&0&0&1&1&1&1\\
1&1&1&1&1&1&1&1\\
1&1&1&1&1&1&1&1\\
1&1&1&1&1&1&1&1\\
},\quad
B=\matVIII
{
1&1&0&0&0&0&0&0\\
1&1&0&0&0&0&0&0\\
0&0&1&1&0&0&0&0\\
0&0&1&1&0&0&0&0\\
0&0&0&0&1&1&0&0\\
0&0&0&0&1&1&0&0\\
0&0&0&0&0&0&1&1\\
0&0&0&0&0&0&1&1\\
},\\
\Ph(\e,\omega)&=
\case{
\log(\frac{11}{10}\e+1), &\omega_{0}\omega_{1}=43\\
\log(\frac{11}{10}\e+2), &\omega_{0}\omega_{1}=44\\
\log 2, &\omega_{0}\omega_{1}\in \Si\setminus\{43,44\} \text{ and }\om_{0}=\om_{1}\\
0, &\omega_{0}\omega_{1}\in \Si\setminus\{43,44\} \text{ and }\om_{0}\neq \om_{1}\\
4\log \e, &\omega_{0}\omega_{1}\in \{45,46\}\\
\log(\frac{1}{10}\e), &\omega_{0}\omega_{1}\in \{57,67,76,86\}\\
\log \e,&\text{otherwise},
}\\
\text{where } \Si&=\bigcup_{i\in \{1,3,5,7\}}\{ii,i(i+1),(i+1)i,(i+1)(i+1)\}.
}
In these setting, we easily see that the conditions $(\Si.1)$-$(\Si.3)$, $(\Ph.1)$-$(\Ph.3)$ and $\sharp \bT_{0}=4$ are satisfied. In particular, we see $\lam(2,\e)=\lam(\{3,4\},\e)$.  Since $\Ph(\e,\om)$ depends only on two coordinates $\om_{0}\om_{1}$, the eigenfunction $g(\e,\om)=g(\e,\om_{0})$ depends only on one coordinate $\om_{0}$. We write $\bT_{0}=\{1,2,3,4\}$ by regarding as $k=B_{kk}$ for $k=1,2,3,4$.
We have the form
\ali
{
\mV_{\bbM,\e}=
\matIV
{
3&2\e&2\e&2\e\\
2\e&\frac{11}{10}\e+3&\e^{4}&0\\
2\e&\e&3&\frac{\e}{5}\frac{g(\e,7)+10g(\e,8)}{g(\e,7)+g(\e,8)}\\
2\e&2\e&\frac{\e}{5}\frac{10g(\e,5)+g(\e,6)}{g(\e,5)+g(\e,6)}&3\\
}.}
This yields the following:
\prop{\label{prop:lam234-lam2toinf}
Under the above notation, 
\ali
{
\frac{\lam(\{2,3,4\},\e)-\lam(2,\e)}{\lam^{v}(\{2,3,4\},\e)-\lam^{v}(2,\e)}\to \infty.
}
}
\pros
Indeed,
\ali
{
\lam(2,\e)&=\lam^{v}(2,\e)=3+\frac{11}{10}\e\\
\lam(\{2,3,4\},\e)&=3+\frac{11}{10}\e+\frac{\sqrt{6}}{2} e^{5/2}+o(e^{3})\\
\lam^{v}(\{2,3,4\},\e)&=3+\frac{11}{10}\e+\frac{20}{3}\frac{14884 c^2+9028 c+2035}{14884 c^2+5368 c-605}\e^{3}+o(\e^{3})
}
with $c=\cos(\frac{1}{3}\arctan((180\sqrt{1273610})/101269))$ are satisfied 
\proe
By a similar argument, we obtain
\ali
{
&\frac{\lam(\{2,3,4\},\e)-\lam(\{2,3\},\e)}{\lam^{v}(\{2,3,4\},\e)-\lam(\{2,3\},\e)}\to \infty,\quad \frac{\lam(\{2,3,4\},\e)-\lam(\{2,4\},\e)}{\lam^{v}(\{2,3,4\},\e)-\lam(\{2,4\},\e)}\to \infty,\\
&\frac{\lam(\{2,3,4\},\e)-\lam(\{3,4\},\e)}{\lam^{v}(\{2,3,4\},\e)-\lam(\{3,4\},\e)}\to 0.
}

\endthebibliography
\end{document}